\documentclass[11pt,a4paper]{article}
\usepackage{amsmath,amssymb,amsthm}
\usepackage{color}
\usepackage{ascmac}
\usepackage{dsfont}

\theoremstyle{definition}
 \newtheorem{dfn}{Definition}[section]
 
  \newtheorem{rmk}[dfn]{Remark}

\theoremstyle{plain}
 \newtheorem{thm}[dfn]{Theorem}
 \newtheorem{prop}[dfn]{Proposition}
 \newtheorem{lem}[dfn]{Lemma}
 \newtheorem{cor}[dfn]{Corollary}
 
\newtheorem*{Rmk}{Remark}

\newcommand{\bn}{{\mathbf n}}

\newcommand{\bk}{{\mathbf k}}
\newcommand{\bp}{{\mathbf p}}
\newcommand{\bu}{{\mathbf u}}

\newcommand{\bv}{{\mathbf v}}
\newcommand{\bw}{{\mathbf w}}

\newcommand{\bq}{{\mathbf q}}
\newcommand{\bff}{{\mathbf f}}
\newcommand{\bA}{{\mathbf A}}
\newcommand{\bB}{{\mathbf B}}

\newcommand{\bD}{{\mathbf D}}

\newcommand{\bF}{{\mathbf F}}
\newcommand{\bG}{{\mathbf G}}
\newcommand{\bH}{{\mathbf H}}
\newcommand{\bI}{{\mathbf I}}

\newcommand{\bK}{{\mathbf K}}
\newcommand{\bL}{{\mathbf L}}
\newcommand{\bM}{{\mathbf M}}

\newcommand{\bV}{{\mathbf V}}
\newcommand{\bS}{{\mathbf S}}

\newcommand{\bU}{{\mathbf U}}

\newcommand{\DV}{{\rm Div}\,}
\newcommand{\dv}{{\rm div}\,}

\newcommand{\BR}{{\mathbb R}}
\newcommand{\BC}{{\mathbb C}}
\newcommand{\BM}{{\mathbb M}}
\newcommand{\BN}{{\mathbb N}}

\newcommand{\BF}{{\mathbb F}}

\newcommand{\BV}{{\mathbb V}}

\newcommand{\CA}{{\mathcal A}}
\newcommand{\CB}{{\mathcal B}}

\newcommand{\CD}{{\mathcal D}}

\newcommand{\CF}{{\mathcal F}}
\newcommand{\CI}{{\mathcal I}}
\newcommand{\CJ}{{\mathcal J}}

\newcommand{\CL}{{\mathcal L}}
\newcommand{\CM}{{\mathcal M}}

\newcommand{\CG}{{\mathcal G}}
\newcommand{\CR}{{\mathcal R}}

\newcommand{\CT}{{\mathcal T}}
\newcommand{\CH}{{\mathcal H}}
\newcommand{\CO}{{\mathcal O}}
\newcommand{\CP}{{\mathcal P}}
\newcommand{\CQ}{{\mathcal Q}}

\newcommand{\CV}{{\mathcal V}}
\newcommand{\CW}{{\mathcal W}}

\newcommand{\fd}{{\mathfrak d}}

\newcommand{\fQ}{{\mathfrak Q}}

\newcommand{\be}{{\mathbf e}}
\newcommand{\bg}{{\mathbf g}}
\newcommand{\bh}{{\mathbf h}}
\newcommand{\pd}{\partial}

\usepackage{xcolor}

\newcommand{\ep}{\varepsilon}
\newcommand{\pa}{\partial}
\newcommand{\wt}{\widetilde}
\newcommand{\wh}{\widehat}
\newcommand{\ol}{\overline}

\newcommand{\br}{{\mathbf r}}

\newcommand{\BBC}{{\mathbb C}}
\newcommand{\BBN}{{\mathbb N}}

\newcommand{\Bf}{{\mathbf f}}
\newcommand{\Bg}{{\mathbf g}}
\newcommand{\Bn}{{\mathbf n}}
\newcommand{\Bu}{{\mathbf u}}
\newcommand{\Bv}{{\mathbf v}}
\newcommand{\Bw}{{\mathbf w}}
\newcommand{\Fs}{{\mathfrak s}}

\newcommand{\Di}{{\rm Div}\,}
\newcommand{\di}{{\rm div}\,}
\newcommand{\loc}{{\rm loc\,}}
\newcommand{\Hol}{{\rm Hol\,}}
\newcommand{\supp}{{\rm supp\,}}
\newcommand{\Ker}{{\rm Ker\,}}
\newcommand{\Span}{{\rm Span\,}}

\newcommand{\BPsi}{{\mathbf \Psi}}

\numberwithin{equation}{section} 

\usepackage{scalerel,stackengine}
\stackMath
\newcommand\reallywidehat[1]{%
\savestack{\tmpbox}{\stretchto{%
  \scaleto{%
    \scalerel*[\widthof{\ensuremath{#1}}]{\kern-.6pt\bigwedge\kern-.6pt}%
    {\rule[-\textheight/2]{1ex}{\textheight}}
  }{\textheight}%
}{0.5ex}}%
\stackon[1pt]{#1}{\tmpbox}%
}

\newcommand{\vertiii}[1]{{\left\vert\kern-0.25ex\left\vert\kern-0.25ex\left\vert #1 
    \right\vert\kern-0.25ex\right\vert\kern-0.25ex\right\vert}}

\begin{document}

\title{\bf The $L_p$-$L_q$ decay estimate for the multidimensional compressible
flow with free surface in the exterior domain}
\author{Yoshihiro Shibata
\thanks{Department of Mathematics,  Waseda University, 
\endgraf
Department of Mechanical Engineering and Materials Science,
University of Pittsburgh, USA \endgraf
mailing address: 
Department of Mathematics, Waseda University,
Ohkubo 3-4-1, Shinjuku-ku, Tokyo 169-8555, Japan. \endgraf
e-mail address: yshibata325@gmail.com \endgraf
Partially supported by Top Global University Project, 
JSPS Grant-in-aid for Scientific Research (A) 17H0109, and Toyota
 Central Research Institute Joint Research Fund }
and Xin Zhang \thanks{School of Mathematical Sciences,
Tongji University, 
No.1239, Siping Road, Shanghai (200092), China. \endgraf
e-mail address: xinzhang2020@tongji.edu.cn}}

\maketitle

\begin{abstract}
The aim of this paper is to develop the general $L_p$ theory for the barotropic compressible Navier-Stokes equations with the free boundary condition in the exterior domain in $\mathbb{R}^N$ ($N\geq 3$).
By the spectral analysis, we obtain the classical $L_p$-$L_q$ decay estimate for the linearized model problem (with variable coefficients) in view of the partial Lagrangian transformation.
\end{abstract}

\tableofcontents

\section{Introduction}
\subsection{Model}
 The motion of viscous gases in some moving exterior domain 
 $\Omega_t \subset \BR^N$ ($N\geq 3$) is described by the following barotropic compressible Navier-Stokes equations with the free boundary condition:
\begin{equation}\label{1.1}
	\left\{ \begin{aligned}
&\pd_t\rho+ \dv\big((\rho_e+\rho)\bv\big) =0
&&\quad&\text{in} \quad \bigcup_{0<t<T} \Omega_t \times\{t\}, \\
&(\rho_e+\rho)(\pd_t\bv +\bv\cdot\nabla\bv)
-\DV\big(\bS(\bv) - P(\rho_e+\rho)\bI\big) = 0 
&&\quad&\text{in} \quad \bigcup_{0<t<T} \Omega_t \times\{t\}, \\
&\big(\bS(\bv) - P(\rho_e+\rho)\bI\big)\bn_{\Gamma_t} = P(\rho_e)\bn_{\Gamma_t},\,\,\,
V_{\Gamma_t} = \bv\cdot\bn_{\Gamma_t} 
&&\quad&\text{on} \quad \bigcup_{0<t<T} \Gamma_t \times\{t\}. \\
\end{aligned}
\right.
\end{equation}
Given the reference mass density $\rho_e>0$ and the initial condition
\begin{equation}\label{1.2}
(\rho, \bv,\Omega_t)|_{t=0} = (\rho_0, \bv_0,\Omega),
\end{equation}
we need to determine the unknown mass density $\rho+\rho_e,$ velocity 
field $\bv=(v_1, \ldots, v_N)^\top,$ and the moving exterior domain 
$\Omega_t.$ 
In \eqref{1.1}, the Cauchy stress tensor 
$$\bS(\bv) = \mu \bD(\bv) + (\nu-\mu)\dv \bv\bI 
\,\,\,\text{for constants}\,\, \mu,\nu>0,$$
and the doubled deformation tensor
 $\bD(\bv) = \nabla\bv + (\nabla\bv)^\top.$
Here, the $(i,j)$th entry of the matrix $\nabla \bv$ is $\pd_i v_j,$ 
$\bI$ is the $N\times N$ identity matrix,
and $\bM^\top$ is the transposed the matrix $\bM=[M_{ij}].$
In addition, $\DV \bM$ denotes an $N$-vector of functions whose 
$i$-th component is $\sum_{j=1}^N \pd_j M_{ij},$
$\dv \bv = \sum_{j=1}^N \pd_jv_j,$ and $\bv\cdot\nabla = \sum_{j=1}^N v_j\pd_j$ with $\pd_j = \pd/\pd x_j.$

 For the last equation in \eqref{1.1},
 $\bn_{\Gamma_t}$ is the outer unit normal vector to the boundary $\Gamma_t$ of $\Omega_t,$ 
and $V_{\Gamma_t}$ stands for the normal velocity of the moving surface $\Gamma_t.$ 
Moreover, the pressure law $P(\cdot)$ is a smooth function defined on $\BR_+.$ In fact, $P(\rho_e)$ coincides with the pressure of the atmosphere.
\medskip

The long time issue for the compressible Navier-Stokes equations in $\BR^3$ is first studied by Matsumura and Nishida in \cite{MN1980} provided the initial data with the $H^3(\BR^3)$ regularity.
Moreover, the authors  in \cite{MN1979} obtained the decay properties of the solutions. For instance, they proved the following $L_2$-$L_1$ type estimates for the perturbation $(\rho,\bv)$ near the equilibrium $(\rho_e, 0)$,
\begin{equation}\label{eq:MN80}
\|(\rho, \bv)\|_{L_2(\BR^3)} \leq C_0 (1+t)^{-3/4},
\end{equation} 
for some constant $C_0$ depending on the small quantity $\|(\rho_0,\bu_0)\|_{L_1(\BR^3)\cap H^3(\BR^3)}.$ 
The decay rate in \eqref{eq:MN80} was recently improved in \cite{LZ2011} by assuming negative Besov regularity for the initial data.
Beyond the classical works \cite{MN1979,MN1980}, Danchin \cite{Dan2000} constructed the global solution with the hyrid Besov regularity in the $L_2$ framework. 
Furthermore, the extension to the general $L_p$ framework was done in the works \cite{CD2010,CMZ2010,Has2011}. Very recently, the decay property of the solution of the compressible Navier-Stokes equations in the Besov norms was investigated by \cite{DanX2017,Okita2014}.
\medskip

A natural question to \eqref{eq:MN80} is whether the general $L_p$-$L_q$ type estimates hold, especially for the boundary value problem in the exterior domain. Namely, the $L_p$ norm of the solution decays provided the initial states in some $L_q$ space. 
For example, recalling the Cauchy problem of the heat flow in $\BR^N$ ($N\geq 2$), we have 
\begin{equation}\label{eq:heat_pq}
\|\pa_{x}^{\alpha}e^{t\Delta} \bu_0\|_{L_p(\BR^N)} \leq C t^{-N(1/q-1/p)/2-|\alpha|/2} \|\bu_0\|_{L_q(\BR^N)}
\end{equation}
for any $1\leq q\leq p\leq \infty,$ $\alpha \in \BN_0^N$ and $t>0.$ Here $\BN_0$ denotes the set of all nonnegative integers.
The extension of \eqref{eq:heat_pq} for the incompressible flow in the exterior domain was first done by Iwashita \cite{Iwa1989}. In \cite{Iwa1989}, as well as the later contribution \cite{MS1997} by Maremonti and Solonnikov within the framework of the potential theory,  the $L_p$-$L_q$ decay property of the Stokes operator $A_{_S}$ associated to the Dirichlet boundary condition in the smooth exterior domain $\Omega \subset \BR^N$ ($N\geq 3$) is established:
\begin{equation}\label{es:Stokes}
\begin{aligned}
\|e^{t A_{_S}} \bu_0\|_{L_p(\Omega)} 
&\leq  C t^{-N(1/q-1/p)/2} \|\bu_0\|_{L_q(\Omega)},\\
 \|\nabla e^{t A_{_S}} \bu_0\|_{L_p(\Omega)} 
 &\leq C t^{-\sigma_1(p,q,N)} \|\bu_0\|_{L_q(\Omega)},
\end{aligned}
\end{equation}
for $t>1,$ $1<q\leq p<\infty$ and 
\begin{equation*}
\sigma_1(p,q,N)= 
\begin{cases}
(N\slash q-N\slash p)\slash 2 + 1\slash 2 & \text{for}\,\,\, 1< p\leq N,\\
N\slash (2q) & \text{for}\,\,\,  N<p<\infty.
\end{cases}
\end{equation*}
Moreover, \cite{MS1997} proves that the gradient estimate of $e^{t A_{_S}}$ in \eqref{es:Stokes} is also sharp for $p>N.$
If one considers the Stokes operator $A_{_S}$ in the 2-D exterior domain, the theory in \cite{DS1999,MS1997} implies that the $L_p$-$L_q$ decay rate becomes worse in the plane than the higher dimensional case. For other discussion on the decay property of the semigroup generated by the Stokes operator, we refer to \cite{HS2009,Ku2007,Shi2018,ShiShi2007a} and the references therein.
\smallbreak

However, concerning the compressible system, one can not expect such $L_p$-$L_q$ estimates for all indices $(p,q)$ as in \eqref{eq:heat_pq} or \eqref{es:Stokes} even 
for the Cauchy problem 
in view of the results \cite{HZ1995,HZ1997,Ponce1985,KS2002}. 
Roughly speaking, the linearized equations of \eqref{1.1} are no longer purely parabolic (see $\eqref{eq:Lame}_1$ for instance).
For the compressible system with the Dirichlet boundary condition, the $L_p$-$L_q$ type estimate was studied in \cite{KS1999} in the exterior domain, and the optimal decay rate like \eqref{eq:MN80} in the half space $\BR^3_+$ was verified in \cite{KK2005}. 
The aim of this paper is to prove the $L_p$-$L_q$ type estimate for the linearized model of \eqref{1.1} with the free boundary condition whenever $\Omega$ is a smooth exterior domain in $\BR^N$ ($N\geq 3$).

\subsection{Main result}
To analysis \eqref{1.1}, it is necessary to transfer \eqref{1.1} to the model in some fixed reference domain. For the local well-posedness issue of \eqref{1.1}, it is sufficient to use the (full) Lagrangian transformation by choosing the initial domain $\Omega$ as the reference domain (see \cite{EvBS2014} for more details). However, to construct the long time solution of \eqref{1.1}, the so-called \emph{partial} Lagrangian coordinates are more helpful due to the hyperbolicity of the mass equation and the free boundary condition $\eqref{1.1}_3.$ 
For the seek of the generality of our mathematical theory, we will treat some system (i.e.\eqref{LL:Lame_1} below) with the slightly variable coefficients, while our main result also applies for the system with constant coefficients (see the comment after Theorem \ref{thm:main_LpLq}). 
\smallbreak

To define the partial Lagrangian coordinates, we assume that $\Omega \subset \BR^N$ with the boundary $\Gamma$ is an exterior domain such that $\CO = \BR^N\setminus \Omega$ is a subset of the ball $B_R$, centred at origin with radius $R>1.$
Let $\kappa$ be a $C^\infty$ functions which equals to one for $x \in B_{R}$ and vanishes outside of $B_{2R}.$
For such $\kappa$ and some fixed $T>0,$ we define the partial Lagrangian mapping 
\begin{equation}\label{def:LL_w}
x =X_{\bw}(y,T)= y + \int^T_0 \kappa(y) \bw(y, s)\,ds
\quad (\forall \, y\in \Omega),
\end{equation}
with $\bw=\bw(\cdot,s)$ defined in $\Omega.$ 
Moreover, we may take $\bw$ in some maximal regularity class 
\begin{equation*}
\bw \in H^1_{\bar p}\big(0,T;L_{\bar q}(\Omega)^N\big) \cap  
L_{\bar p}\big(0,T;H^2_{\bar q}(\Omega)^N\big)
\end{equation*}
for $1<{\bar p}<\infty$ and $N<{\bar q}<\infty.$
To guarantee the invertibility of $X_{\bw},$ we assume that 
\begin{equation}\label{cdt:small_1}
\int^T_0\|\kappa(\cdot)\bw(\cdot, s)\|_{H^1_\infty(\Omega)}\,ds \leq\delta <1
\end{equation}
for some small constant $\delta.$
By \eqref{def:LL_w}, \eqref{1.1} can be reduced to the following linearized model problem with ignoring all nonlinear terms
\footnote{The details of the linearization of \eqref{1.1} will be present in Appendix \ref{sec:PL}. Although the convection term $\bv \cdot \nabla \rho$ in the mass conservation law can not be simply regarded as a perturbation term for the global well-posedness issue of \eqref{1.1}, it is fundamental to derive the decay properties of \eqref{LL:Lame_1}.},
\begin{equation}\label{LL:Lame_1}
	\left\{ \begin{aligned}
&\pd_t\rho 
+ \gamma_1\, \ol\dv \bv  =0
&&\quad&\text{in}& \quad \Omega \times \BR_+, \\
&\gamma_1\pd_t\bv 
-\ol\DV\big(\ol\bS(\bv) - \gamma_2\rho\bI\big) =0
&&\quad&\text{in}& \quad \Omega \times \BR_+, \\
&\big(\ol\bS(\bv) - \gamma_2 \rho\bI\big)\ol\bn_{\Gamma} = 0
&&\quad&\text{on}& \quad \Gamma \times \BR_+, \\
&(\rho, \bv)|_{t=0} = (\rho_0, \bv_0)
&&\quad&\text{in}& \quad \Omega, \\
\end{aligned}
\right.
\end{equation}
where we have set that 
\begin{gather}
\gamma_1=\rho_e, \quad \gamma_2= P'(\rho_e), \quad 
\bI+\bV=\big(\nabla_y X_{\bw}\big)^{-1}(\cdot,T),\quad 
 J= \big(\det(\bI+\bV)\big)^{-1},\nonumber\\ 
\ol\dv\bv = (\bI+\bV):\nabla_y \bv
=J^{-1}\dv_y \big( J (\bI+ \bV)^\top \bv\big), \nonumber\\
\ol \DV\bM=  J^{-1} \DV_{y} \big(J \bM (\bI + \bV) \big)
=\DV_y \bM+(\bV \nabla_y|\bM), \label{LL:sym_0}\\
\ol \bD (\bv) =(\bI + \bV)\nabla_{y} \bv 
+(\nabla_{y} \bv)^{\top}\big(\bI + \bV\big)^{\top},\nonumber\\
\ol\bS(\bv)=\mu \ol\bD(\bv)+(\nu-\mu) (\ol\dv \bv) \bI, \quad 
 \ol\bn_{\Gamma}= (\bI + \bV)\bn_{\Gamma}.\nonumber
\end{gather}
Above, for any matrices $\bA=[A_{jk}]$ and $\bB=[B_{jk}],$ 
$\bA:\bB$ denotes the trace of 
the product $\bA\bB,$ that is,
\begin{equation*}
\bA:\bB=\sum_{j,k=1}^N A_{jk} B_{kj},
\end{equation*}
and $(\bB \nabla_y|\bA)$ stands for an $N$-vector with the $i$th component $\sum_{j,k} B_{jk} \pd_k A_{ij}.$ 
In particular, it is easy to see that 
\begin{equation*}
\ol \DV (f\bI)=\nabla_y f + (\bV \nabla_y \mid f\bI) = (\bI +\bV) \nabla_y f,
\end{equation*}
for any smooth function $f.$
For simplicity, we denote $\ol\nabla =(\bI +\bV) \nabla_y$ in what follows.
Hereafter, we assume that 
\begin{equation}\label{assump:1}
\|(J^{\pm 1}-1, \bV)\|_{L_\infty(\Omega)}
+\|\nabla_y (J^{\pm 1}, \bV)\|_{L_r(\Omega)}\leq \sigma(<<1)
\end{equation}
with some small constant $\sigma=\sigma(\delta)$ and $N < r < \infty$. 
By the assumption \eqref{assump:1}, \eqref{LL:Lame_1} can be regarded as the perturbation of the standard Lam\'e operator, which is the motivation of our main result.
\smallbreak

Following \cite{EvBS2014}, it is not hard to see that there exists a $C_0$-semigroup $\{T(t)\}_{t\geq 0}$ generated by the operator
\begin{equation*}
\CA_{\Omega}(\rho,\bv)=
\Big(\gamma_1 \ol\di \bv,-\gamma_1^{-1}\ol\Di \big(\ol\bS(\bv)-\gamma_2\rho\bI\big)\Big)
\end{equation*}
in the space $H^{1,0}_{p}(\Omega)=H^1_p(\Omega)\times L_p(\Omega)^N$ for $1<p<\infty$
(see Theorem \ref{thm:sg} in section \ref{sec:large}). 
Sometimes,  we also denote the solution of \eqref{LL:Lame_1} by 
$(\rho,\bv)=T(t)(\rho_0,\bv_0)$ and 
$\bv=\CP_v T(t)(\rho_0,\bv_0).$ 
Now our main result for the $L_p$-$L_q$ decay estimate of $\{T(t)\}_{t\geq 0}$ reads as follows.
\begin{thm}{($L_p$-$L_q$ decay estimate)}
\label{thm:main_LpLq}
Let $\Omega$ be a $C^{3}$ exterior domain in $\BR^N$ with $N\geq 3,$
and let \eqref{assump:1} be satisfied for some $N<r<\infty.$
Assume that 
$(\rho_0,\bv_0) \in L_q(\Omega)^{1+N}\cap H^{1,0}_p(\Omega)$
with $H^{1,0}_{p}(\Omega)=H^1_p(\Omega)\times L_p(\Omega)^N$
for $1\leq q \leq 2 \leq p \leq r<\infty,$ 
and $\{T(t)\}_{t\geq 0}$ is the semigroup associated to \eqref{LL:Lame_1} 
in $H^{1,0}_{p}(\Omega).$ 
For convenience, we set
\begin{equation*}
\vertiii{(\rho_0, \bv_0)}_{p,q} =\|(\rho_0, \bv_0)\|_{L_q(\Omega)} 
+ \|(\rho_0, \bv_0)\|_{H^{1,0}_p(\Omega)}.
\end{equation*}
Then for $t\geq 1,$ there exists a positive constant $C$ such that
\begin{align*}
\|T(t)(\rho_0, \bv_0)\|_{L_p(\Omega)} &\leq C 
t^{-(N\slash q-N\slash p)\slash 2} \vertiii{(\rho_0, \bv_0)}_{p,q}, \\
\|\nabla T(t)(\rho_0, \bv_0)\|_{L_p(\Omega)} &\leq C 
t^{-\sigma_1(p,q,N)} \vertiii{(\rho_0, \bv_0)}_{p,q}, \\
\|\nabla^2 \CP_v T(t)(\rho_0, \bv_0)\|_{L_p(\Omega)} &\leq C 
t^{-\sigma_2(p,q,N)} \vertiii{(\rho_0, \bv_0)}_{p,q}, 
\end{align*}
where the indices $\sigma_1(p,q,N)$ and $\sigma_2(p,q,N)$ 
are given by
\begin{align*}
\sigma_1(p,q,N) &= 
\begin{cases}
(N\slash q-N\slash p)\slash 2 + 1\slash 2 & \text{for}\,\,\, 2\leq p\leq N,\\
N\slash (2q) & \text{for}\,\,\,  N<p<\infty,\\
\end{cases}\\
\sigma_2(p,q,N) &= 
\begin{cases}
3\slash (2q) & \text{for}\,\,\, N=3,\\
(N\slash q-N\slash p)\slash 2 + 1 & \text{for}\,\,\, N\geq 4\,\,\, \text{and}\,\,\, 
2\leq p\leq N\slash 2,\\
N\slash (2q) & \text{for}\,\,\,   N\geq 4\,\,\, \text{and}\,\,\, 
N\slash 2<p<\infty.
\end{cases}
\end{align*}
\end{thm}

\begin{Rmk}
Let us give some comments on Theorem \ref{thm:main_LpLq} above.
\begin{enumerate}
\item \eqref{LL:Lame_1} is a system with variable coefficients only near the boundary $\Gamma.$
Notice that 
\begin{equation}\label{support_0}
\bV=0, \quad J = 1, \quad \overline{\dv}\bv= \dv\bv, 
\quad \ol\bD(\bv) = \bD(\bv), \quad 
\ol\bS(\bv) = \bS(\bv) 
\end{equation}
outside of the ball $B_{2R}$ by the choice of $\kappa$ in \eqref{def:LL_w}. 

\item Taking $\bV=0$ and $J=1$ in Theorem \ref{thm:main_LpLq}, we have in particular the same decay properties for the system with constant coefficients, 
\begin{equation}\label{eq:Lame}
	\left\{ \begin{aligned}
&\pd_t\rho
+ \gamma_1\, \dv \bv  =0
&&\quad&\text{in}& \quad \Omega \times \BR_+, \\
&\gamma_1\pd_t\bv 
-\DV\big(\bS(\bv) - \gamma_2\rho\bI\big) =0
&&\quad&\text{in}& \quad \Omega \times \BR_+, \\
&\big(\bS(\bv) - \gamma_2\rho\bI\big)\bn_{\Gamma} =0
&&\quad&\text{on}& \quad \Gamma \times \BR_+, \\
&(\rho, \bv)|_{t=0} = (\rho_0, \bv_0)
&&\quad&\text{in}& \quad \Omega.
\end{aligned}
\right.
\end{equation}
Of course, the assumption \eqref{assump:1} is not necessary for \eqref{eq:Lame} any longer.

\item For simplicity, we will not seek for the optimal assumption on the regularity of $\Omega.$ The assumption of $\Omega$ in $C^3$ class is for the technical reason, as the property of the Bogovskii operator will be used in the later proof.
\end{enumerate}
\end{Rmk}

To establish the $L_p$-$L_q$ estimates in Theorem \ref{thm:main_LpLq}, we use the so-called 
\emph{local energy approach}. More precisely, we shall establish the following theorem.
\begin{thm}{(local energy estimate)}
\label{thm:local_energy}
Let $\Omega$ be a $C^3$ exterior domain in $\BR^N$ for $N\geq 3,$ $N<r<\infty,$ $1<p\leq r,$ and $L> 2R.$ 
Denote that \footnote{$\overline{E}$ stands for the closesure of $E$ for any subset 
$E \subset \BR^N.$} 
\begin{gather*}
\Omega_L = \Omega \cap B_L, \quad
H^{1,2}_p(\Omega_L)=H^{1}_p(\Omega_L)
\times H^{2}_p(\Omega_L)^N,\\
X_{p,L}(\Omega)= \{(d,\bff)\in H^{1,0}_p(\Omega) 
: \supp d, \,\supp \bff \subset \overline{\Omega_L}\,\}.
\end{gather*}
Then for any $(\rho_0, \bv_0) \in X_{p,L}(\Omega)$ and $k\in \BBN_0=\BBN\cup \{0\},$ there exists a positive constant $C_{p,k,L}$ such that
\begin{equation*}
\| \pa^k_tT(t)(\rho_0,\bv_0)\|_{H^{1,2}_p(\Omega_L)} 
\leq C_{p,k,L}\, t^{-N\slash 2-k}\|(\rho_0, \Bv_0)\|_{H^{1,0}_p(\Omega)}, 
\,\,\,\forall \,\, t\geq 1.
\end{equation*}
\end{thm}

To prove Theorem \ref{thm:local_energy}, we consider the resolvent problem of \eqref{LL:Lame_1}:
\begin{equation}\label{resolvent_0}
\left\{\begin{aligned}
&\lambda \eta + \gamma_1\overline{\dv}\bu = d 
&&\quad &\text{in}& \quad \Omega, \\
&\gamma_1\lambda \bu 
- \ol\DV\big(\ol\bS(\bu)-\gamma_2\eta \bI \big)= \bff 
&&\quad &\text{in}& \quad \Omega, \\
&\big(\ol\bS(\bu)-\gamma_2\eta \bI \big)\ol\bn_{\Gamma}=0 
&&\quad &\text{on}& \quad \Gamma.
\end{aligned}\right.
\end{equation}
In fact, it is easy to study \eqref{resolvent_0} whenever $\lambda$ is sufficient large (see Theorem \ref{thm:large_res}). The case $\lambda$ is uniformly bounded from above is more involved, which is the core concern of this paper. Especially, one difficulty appears when we handle \eqref{resolvent_0} for $\lambda=0.$ Namely, there always exists non-trivial stationary solution for the stationary Lam\'e system (see \cite[chapter 6]{Lions1998} for more details). 
Here, one can first consider the following simplified model problem in the whole space $\BR^3$:
\begin{equation}\label{eq:model_3}
- \Delta \bu - \nabla \di \bu + \nabla \eta = \Bf.
\end{equation}
Obviously, $(\bu, \eta)=(0,F)$ always solves \eqref{eq:model_3} as long as $\Bf =\nabla F$ for some $ F \in C_0^{\infty}(\BR^3),$ which is quite annoying for our study of \eqref{resolvent_0} with $\lambda=0.$ 
To handle such trouble, we introduce some weighted inner product structure for the rigid motion space (see subsection \ref{subsec:para}) to eliminate the homogeneity effect from the stationary Lam\'e system.
\smallbreak

This paper is folded as follows.
In sections \ref{sec:important} and \ref{sec:near}, we will discuss the case where $\lambda$ in \eqref{resolvent_0} is closed to zero.
Afterwards, via the sections \ref{sec:aux} and \ref{sec:large},  we study \eqref{resolvent_0} when $\lambda$ is uniformly bounded from below and also from above.
Finally, we combine the previous results and prove Theorems \ref{thm:local_energy} and \ref{thm:main_LpLq} in the last section.
In Appendix \ref{sec:PL}, we give the derivation of \eqref{LL:Lame_1} from \eqref{1.1}.

\subsection*{Notation}
For convenience, we introduce some useful notation.
For any domain $G$ in $\BR^N,$ $1\leq p\leq\infty$ and $k\in \BN,$   $L_p(G)$ ($L_{p,\loc}(G)$) stands for the (local) Lebesgue space, 
and $H^k_p(G)$ ($H^k_{p,\loc}(G)$)  for the (local) Sobolev space. 
Moreover, we write 
\begin{equation*}
H^{k,\ell}_p(G) = H^{k}_p(G) \times  H^{\ell}_p(G)^N, \quad
H^{k,\ell}_{p,\loc}(G) = H^{k}_{p,\loc}(G) \times  H^{\ell}_{p,\loc}(G)^N.
\end{equation*}
In addition, the Besov space $B^{s}_{q,p}(G)$ 
for $k-1<s\leq k$ and $1<p,q<\infty$
is defined by the real interpolation functor
\begin{equation*}
B_{q,p}^{s}(G)= \big(L_q(G),H^{k}_q(G)\big)_{s\slash k,\,p}.
\end{equation*} 
Sometimes, we write $W^{s}_q(G)=B^s_{q,q}(G)$ for simplicity.
\smallbreak

For any Banach spaces $X,Y,$ the total of the bounded linear transformations from $X$ to $Y$ is denoted by $\CL(X;Y).$ We also write $\CL(X)$ for short if $X=Y.$ 
In addition, $\Hol (\Lambda;X)$ denotes the set of $X$-valued analytic mappings defined on the domain $\Lambda \subset \BBC.$
\smallbreak

To study the resolvent problem \eqref{resolvent_0}, we introduce that 
\begin{gather}
\Sigma_{\ep} = \{\lambda \in \BC \backslash \{ 0\}: 
|\arg \lambda|\leq \pi-\ep\}, \quad
\Sigma_{\ep,b}=\{\lambda \in \Sigma_{\ep}: |\lambda| \geq b\}, 
\nonumber\\
 K= \Big\{\lambda \in \BBC: 
\big(\Re \lambda +\frac{\gamma_1\gamma_2}{\mu+\nu} \big)^2 
+\Im \lambda^2 > \big(\frac{\gamma_1\gamma_2}{\mu+\nu} \big)^2
\Big\},
\label{def:domain} \\ \nonumber
V_{\ep,b}=\Sigma_{\ep,b} \cap K, \quad
\dot{U}_{b}=\{\lambda \in 
\BC \backslash (-\infty, 0]: |\lambda| <b\}
\end{gather}
for any $0<\ep<\pi/2$ and $b>0.$

\section{Important propositions for local energy decay}
\label{sec:important}
In this section, we will consider several model problems in the bounded domain $\Omega_{5R}=\Omega \cap B_{5R},$ which play a vital role in the construction of the solution of \eqref{resolvent_0} in the next section. 

\subsection{The divergence equation with variable coefficients}
First, we consider the modified divergence operator in \eqref{LL:sym_0}, that is,
\begin{equation}\label{eq:div_1}
\ol\dv \bu=
\dv \bu + \sum_{j,k=1}^N V_{jk}\frac{\pd u_j}{\pd y_k} = f
\quad \text{in}\,\,\, \Omega_{5R}=\Omega \cap B_{5R},
\end{equation}
with $\bV=[V_{jk}]$ satisfying 
\eqref{assump:1} for some $0<\sigma<1$ and $N<r<\infty.$

\begin{prop}\label{prop:div} 
Let $1 < p \leq r$.  Then, there exist positive
constants $\sigma_0=\sigma_{0}(p,r)$ and $M=M(p,R)$ such that 
the equation \eqref{eq:div_1} admits a solution $\bu \in H^2_p(\Omega_{5R})^N$ 
possessing the estimate:
\begin{equation*}
\|\bu\|_{H^2_p(\Omega_{5R})} \leq M\|f\|_{H^1_p(\Omega_{5R})}
\end{equation*}
for any $0 < \sigma \leq \sigma_0.$ 
\end{prop}

Hereafter, the following estimate in any uniform $C^1$ domain $G$ will be constantly used,
\begin{equation}\label{eq:sobolev}
\|gh\|_{L_p(G)} \leq C \|g\|_{L_r(G)} \|h\|_{H^1_p(G)}
\end{equation}
for any $r>N$ and $1<p\leq r.$ 
In fact, when $p=r > N$, we have 
\begin{equation*}
\|gh\|_{L_p(G)} \leq C \|g\|_{L_p(G)} \|h\|_{L_\infty(G)}
\leq C \|g\|_{L_r(G)} \|h\|_{H^1_p(G)}.
\end{equation*}
For $1 < p < r$, using the H\"older's inequality for 
$1/r + 1/s = 1/p$ and 
the Sobolev embedding theorem yield that
\begin{equation*}
\|gh\|_{L_p(G)} \leq C \|g\|_{L_r(G)} \|h\|_{L_s(G)}
\leq C \|g\|_{L_r(G)} \|h\|_{H^1_p(G)}.
\end{equation*}
\begin{proof}
To construct the solution $\bu$ of \eqref{eq:div_1}, we solve the following equation with Dirichlet boundary condition
\begin{equation}\label{eq:div_2}
\Delta \theta + \sum_{j,k=1}^N V_{jk}\frac{\pd^2 \theta}{\pd y_j\pd y_k}=f
\quad \text{in $\Omega_{5R}$}, \quad \theta|_{\pd \Omega_{5R}}=0,
\end{equation}
where $\pd\Omega_{5R}= \Gamma \cup S_{5R}$ denotes the boundary of $\Omega_{5R}.$
As $\bV$ satisfies \eqref{assump:1}, the left-hand side of \eqref{eq:div_2} can be regarded as the perturbation of the Laplace operator. Thus we use the standard Banach fixed point argument to find $\theta.$

Given $\omega \in H^3_p(\Omega_{5R})$ with $\omega|_{\pd\Omega_{5R}}=0,$ we consider the following equation
\begin{equation}\label{eq:div_3}
\Delta \theta =
f -  \sum_{j,k=1}^N V_{jk}\frac{\pd^2 \omega}{\pd y_j\pd y_k}
\quad \text{in $\Omega_{5R}$}, \quad \theta|_{\pd \Omega_{5R}}=0.
\end{equation}
Then it is not hard to see from \eqref{eq:sobolev} that 
\begin{equation*}
\Big\|V_{jk}\frac{\pd^2 \omega}{\pd y_j\pd y_k}\Big\|_{H^1_p(\Omega_{5R})}
\leq C\big(\|\bV\|_{L_\infty(\Omega_{5R})} 
+ \|\nabla\bV\|_{L_r(\Omega_{5R})} \big)
\|\omega\|_{H^3_p(\Omega_{5R})},
\end{equation*}
with the constant $C$ depending solely on $p,$ $r$ and $\Omega_{5R}.$ 
\smallbreak 

Next, in view of \eqref{assump:1}, there exists a solution $\theta \in H^3_p(\Omega_{5R})$ of \eqref{eq:div_3} fulfilling
\begin{equation}\label{eq:4.6}
\|\theta\|_{H^3_p(\Omega_{5R})} 
\leq \frac{M}{2} \big(\|f\|_{H^1_p(\Omega_{5R})} + 
C\sigma\|\omega\|_{H^3_p(\Omega_{5R})} \big),
\end{equation}
where $M$ is a positive constant depending on $\Omega_{5R}$ and $p.$ 
Thus, if $\omega$ satisfies the estimate
$$\|\omega\|_{H^3_p(\Omega_{5R})}
\leq M\|f\|_{H^1_p(\Omega_{5R})},$$
then \eqref{eq:4.6} yields that
$$\|\theta\|_{H^3_p(\Omega_{5R})} 
\leq (M/2)\|f\|_{H^1_p(\Omega_{5R})} 
+(C\sigma M^2/2)\|f\|_{H^1_p(\Omega_{5R})}.
$$
Choosing $\sigma_0>0$ small enough such that $CM\sigma_0 \leq 1$, we have
$$\|\theta\|_{H^3_p(\Omega_{5R})} 
\leq M\|f\|_{H^1_p(\Omega_{5R})}. $$
For simplicity, we write the solution mapping of \eqref{eq:div_3} by $\theta =\Phi(\omega)$ for $\omega \in H^3_p(\Omega_{5R})$ with $\omega|_{\pd\Omega_{5R}}=0.$ 
\smallbreak

To verify the contraction property of $\Phi,$ we take $\theta_i=\Phi(\omega_i) \in H^3_p(\Omega_{5R})$ for $i=1,2.$
Then the difference $\theta_1-\theta_2$ satisfies the following Poisson's equation
$$\Delta (\theta_1 - \theta_2) 
= -\sum_{j,k=1}^N V_{jk}
\frac{\pd^2 (\omega_1- \omega_2)}{\pd y_j\pd y_k}
\quad \text{in $\Omega_{5R}$}, \quad 
(\theta_1-\theta_2)|_{\pd \Omega_{5R}}=0.$$
Thus, by \eqref{eq:4.6} we have
$$\|\theta_1-\theta_2\|_{H^3_p(\Omega_{5R})}
\leq (M/2)C\sigma\|\omega_1-\omega_2\|_{H^3_p(\Omega_{5R})}
\leq (1/2)\|\omega_1-\omega_2\|_{H^3_p(\Omega_{5R})},
$$
for $CM\sigma_0 \leq 1.$
Therefore, the Banach's fixed point theorem implies that there exists a unique
solution  $\theta \in H^3_p(\Omega_{5R})$ of \eqref{eq:div_2} with
$\|\theta\|_{H^3_p(\Omega_{5R})} \leq M\|f\|_{H^1_p(\Omega)}.$
\smallbreak

At last, setting $\bu= \nabla\theta$, we have the desired result.  
This completes the proof of Proposition \ref{prop:div}. 
\end{proof}

\subsection{Modified Stokes and reduced Stokes operators}
Recall the notion in \eqref{LL:sym_0} and 
the assumption \eqref{assump:1},
we first consider the resolvent problem of the modified Stokes equations as follows
\begin{equation}\label{st:1}
\left\{  \begin{aligned}
&\gamma_1 \lambda \bv
-\ol\DV \big( \mu \ol\bD(\bv)
-\gamma_2\omega\bI\big)= \bff
&&\quad &\text{in}& \quad \Omega_{5R}, \\
&\overline{\dv }\bv = 0
&&\quad &\text{in}& \quad \Omega_{5R}, \\
&\big(\mu \ol\bD(\bv)- \gamma_2\omega\bI\big)
\ol\bn_{_{\pd\Omega_{5R}}}= 0
&& \quad &\text{on}&\quad  \pd \Omega_{5R},
\end{aligned} \right.
\end{equation}
where $\ol\bn_{_{\pd\Omega_{5R}}}=(\bI+\bV)\bn_{_{\pd\Omega_{5R}}}$ and $\bn_{_{\pd\Omega_{5R}}}$ denotes the unit outward normal vector to 
$\pd\Omega_{5R} = \Gamma \cup S_{5R}.$ 

\medskip

\subsubsection{Analysis of \eqref{st:1} for $\lambda \not=0$} 
We will establish the following theorem 
by assuming $\lambda\not=0$ in \eqref{st:1}.
\begin{thm}\label{thm:ms1} 
Assume that $\Omega$ is a $C^2$ exterior domain.
Let $1 < p \leq r,$ $0 < \varepsilon < \pi/2,$ and $\lambda_0 > 0.$
Then, there exists a $\sigma_1$ depending on 
$\mu$, $\varepsilon,$ $\lambda_0,$ $p$, $r$ and $\Omega_{5R}$ 
such that
if $0 < \sigma \leq \sigma_1$,
then for any $\lambda \in \Sigma_{\varepsilon, \lambda_0}$,
 $\bff \in L_p(\Omega_{5R})^N,$ problem \eqref{st:1} admits a
unique solution $(\bv,\omega) \in H^2_p(\Omega_{5R})^N 
\times  H^1_p(\Omega_{5R})$
possessing the estimate:
\begin{align*}
|\lambda|\|\bv\|_{L_p(\Omega_{5R})} 
+|\lambda|^{1/2}\|\bv\|_{H^1_p(\Omega_{5R})} 
+ \|\bv\|_{H^2_p(\Omega_{5R})}
+ \|\omega\|_{H^1_p(\Omega_{5R})}
 \leq C \|\bff\|_{L_p(\Omega_{5R})},
\end{align*}
with some constant $C$ depending on 
$\mu$, $\varepsilon,$ $\lambda_0,$ $p,$ $r$ and $\Omega_{5R}.$ 
\end{thm}

The proof of Theorem \ref{thm:ms1} is similar to Proposition \ref{prop:div}.
By the notion in \eqref{LL:sym_0}, \eqref{st:1} is rewritten as 
\begin{equation}\label{st:1_1}
\left\{  \begin{aligned}
&\gamma_1 \lambda \bv
-\DV \big( \mu \bD(\bv)
-\gamma_2\omega\bI\big)= \bff+\bF(\Bu,\omega)
&&\quad &\text{in}& \quad \Omega_{5R}, \\
&\dv \bv = g(\bv)=\dv\bg(\bv)
&&\quad &\text{in}& \quad \Omega_{5R}, \\
&\big(\mu\bD(\bv)- \gamma_2\omega\bI\big)
\ol\bn_{_{\pd\Omega_{5R}}}= \bH(\bv,\omega)
&& \quad &\text{on}&\quad  \pd \Omega_{5R},
\end{aligned} \right.
\end{equation}
with 
\begin{align*}
\bF(\Bv,\omega) &= \DV\big(\mu \bV \nabla \bv 
+ \mu (\bV \nabla \bv)^\top \big) 
+ \big(\bV \nabla \mid \mu \ol \bD(\bv)
-\gamma_2 \omega \bI \big),\\
g(\bv)&=-\bV:\nabla_y \bu, \quad 
\bg(\bv)= -\bV^{\top} \bv-(J-1)(\bI+\bV)^{\top}\bv,\\
\bH(\bv,\omega)&=-\mu \big( \bV \nabla \bv 
+ (\bV \nabla \bv)^\top\big)\bn_{_{\pd\Omega_{5R}}}
-\big(\mu \ol\bD(\bv)- \gamma_2\omega\bI\big) \bV
\bn_{_{\pd\Omega_{5R}}}.
\end{align*}
Thanks to \eqref{assump:1}, \eqref{st:1} can be regarded as the perturbation of the following linear equations:
\begin{equation}\label{ms:4}
\left\{  \begin{aligned}
&\gamma_1 \lambda \bv
-\DV \big( \mu \bD(\bv)
-\gamma_2\omega\bI\big)= \bff
&&\quad &\text{in}& \quad \Omega_{5R}, \\
&\dv\bv = g= \dv\bg
&&\quad &\text{in}& \quad \Omega_{5R}, \\
&(\mu \bD(\bv)- \gamma_2\omega\bI)
\bn_{_{\pd\Omega_{5R}}}= \bh
&& \quad &\text{on}&\quad  \pd \Omega_{5R}.
\end{aligned} \right.
\end{equation}
For \eqref{ms:4}, we know the following theorem from \cite[Sec. 4]{Shi2020}.
\begin{thm}\label{Stokes:1} 
Assume that $\Omega$ is a $C^2$ exterior domain.
Let $1 < p < \infty$, $0 < \varepsilon < \pi/2$, and $\lambda_0 > 0$.
Then, for any $\lambda \in \Sigma_{\varepsilon, \lambda_0},$
$\bff \in L_p(\Omega_{5R})^N,$ 
$g \in H^1_p(\Omega_{5R}),$ 
$\bg \in L_p(\Omega_{5R})^N$ and 
$\bh \in H^1_p(\Omega_{5R})^N,$ problem \eqref{ms:4} admits a
unique solution $(\bv,\omega) \in 
H^2_p(\Omega_{5R})^N \times H^1_p(\Omega_{5R})$
possessing the estimate:
\begin{align*}
&|\lambda|\|\bv\|_{L_p(\Omega_{5R})}
 + |\lambda|^{1/2}\big(\|\bv\|_{H^1_p(\Omega_{5R})}
 +\|\omega\|_{L_p(\Omega_{5R})} \big)
+\|\bv\|_{H^2_p(\Omega_{5R})}
+ \|\omega\|_{H^1_p(\Omega_{5R})} \\
&\quad \leq C(\|\bff\|_{L_p(\Omega_{5R})}
+|\lambda|^{1/2} \|(g, \bh)\|_{L_p(\Omega_{5R})}
+ \|(g, \bh)\|_{H^1_p(\Omega_{5R})} 
+ |\lambda|\|\bg\|_{L_p(\Omega_{5R})})
\end{align*}
with some constant $C$ depending on 
$\mu$, $\varepsilon,$ $\lambda_0,$ $p$ and $\Omega_{5R}.$ 
\end{thm}

Now, let us give the proof of Theorem \ref{thm:ms1}
by the standard Banach fixed point theorem and Theorem \ref{Stokes:1}.
\begin{proof}[Proof of Theorem \ref{thm:ms1}]

Given $(\bu,\theta) \in H^2_p(\Omega_{5R})^N 
\times H^1_p(\Omega_{5R}),$ we consider
\begin{equation}\label{ms:5}
\left\{  \begin{aligned}
&\gamma_1 \lambda \bv
-\DV \big( \mu \bD(\bv)
-\gamma_2\omega\bI\big)= \bff+\bF(\Bu,\theta)
&&\quad &\text{in}& \quad \Omega_{5R}, \\
&\dv \bv = g(\bu)=\dv\bg(\bu)
&&\quad &\text{in}& \quad \Omega_{5R}, \\
&\big(\mu\bD(\bv)- \gamma_2\omega\bI\big)
\bn_{_{\pd\Omega_{5R}}}= \bH(\bu,\theta)
&& \quad &\text{on}&\quad  \pd \Omega_{5R}.
\end{aligned} \right.
\end{equation}
By \eqref{assump:1} and \eqref{eq:sobolev}, it is not hard to find that
\begin{equation}\label{ms:6}
\begin{aligned}
\|\bF(\bu,\theta)\|_{L_p(\Omega_{5R})} 
& \leq C\sigma \big(\|\bu\|_{H^2_p(\Omega_{5R})}
+\|\theta\|_{H^1_p(\Omega_{5R})} \big), \\
\|g(\bu)\|_{H^k_p(\Omega_{5R})} 
& \leq C\sigma\|\bu\|_{H^{k+1}_p(\Omega_{5R})},  
\quad \forall\,\,\,k=0,1,\\
\|\bg(\bu)\|_{L_p(\Omega_{5R})} & 
\leq C\sigma\|\bu\|_{L_p(\Omega_{5R})},\\
\| \bH(\bu,\theta)\|_{H^k_p(\Omega_{5R})}
& \leq  C\sigma \big(\|\bu\|_{H^{k+1}_p(\Omega_{5R})}
+\|\theta\|_{H^{k}_p(\Omega_{5R})} \big), 
\quad \forall\,\,\,k=0,1,
\end{aligned}
\end{equation}
for some constant $C>0.$
For simplicity, we set 
\begin{equation*}
N_{\lambda}(\bu,\theta;\Omega_{5R})
=|\lambda| \|\bu\|_{L_p(\Omega_{5R})}+
|\lambda|^{1/2}\big(\|\bu\|_{H^1_p(\Omega_{5R})}
 +\|\omega\|_{L_p(\Omega_{5R})} \big)
+ \|\bu\|_{H^2_p(\Omega_{5R})}
+ \|\theta\|_{H^1_p(\Omega_{5R})}.
\end{equation*}
Then, for any fixed $\lambda \in \Sigma_{\ep,\lambda_0},$
we infer from \eqref{ms:6} that
\begin{multline}\label{ms:7}
 \|\bF(\bu,\theta)\|_{L_p(\Omega_{5R})}
+|\lambda|^{1/2}\big\|\big(g(\bu),
 \bH(\bu,\theta)\big)\big\|_{L_p(\Omega_{5R})}
+ \big\|\big(g(\bu), \bH(\bu,\theta)\big)\big\|_{H^1_p(\Omega_{5R})}\\
+|\lambda|\|\bg(\bu)\|_{L_p(\Omega_{5R}))}
\leq C\sigma N_{\lambda}(\bw,\theta;\Omega_{5R}).
\end{multline}
Thus, by Theorem \ref{Stokes:1} and \eqref{ms:7}, there exists a solution $(\bv,\omega)$ of \eqref{ms:5} in 
$H^2_p(\Omega_{5R})^N \times H^1_p(\Omega_{5R})$ 
satisfying
\begin{equation}\label{ms:8}
N_{\lambda}(\bv,\omega;\Omega_{5R})
\leq (M/2)\big( \|\bff\|_{L_p(\Omega_{5R})} 
+ C\sigma N_{\lambda}(\bu,\theta;\Omega_{5R})\big).
\end{equation}
Thus, if $N_{\lambda}(\bu,\theta;\Omega_{5R}) 
\leq M \|\bff\|_{L_p(\Omega_{5R})},$
then choosing $\sigma > 0$ so small that $C\sigma M \leq 1$ in \eqref{ms:8} we have
\begin{equation}\label{ms:9}
N_{\lambda}(\bv,\omega;\Omega_{5R}) 
\leq M\|\bff\|_{L_p(\Omega_{5R})}.
\end{equation}

Next, for $(\bu_i, \theta_i) \in H^2_p(\Omega_{5R})^N 
\times H^1_p(\Omega_{5R})$ ($i = 1,2$), 
let $(\bv_i,\omega_i)$ be the solution of \eqref{ms:5} with respect to 
$(\bu_i, \theta_i).$
Noting that the functionals $\bF$, $g$, $\bg$ and $\bH$ are linear, we have
\begin{equation*}
\left\{  \begin{aligned}
&\gamma_1 \lambda (\bv_1-\bv_2)
-\DV \big( \mu \bD(\bv_1-\bv_2)
-\gamma_2(\omega_1-\omega_2)\bI\big)
= \bF(\bu_1-\bu_2,\theta_1-\theta_2)
&&\quad &\text{in}& \quad \Omega_{5R}, \\
&\dv (\bv_1-\bv_2) = g(\bu_1-\bu_2)=\dv\bg(\bu_1-\bu_2)
&&\quad &\text{in}& \quad \Omega_{5R}, \\
&\big(\mu\bD(\bv_1-\bv_2)
- \gamma_2(\omega_1-\omega_2)\bI\big)
\bn_{_{\pd\Omega_{5R}}}= \bH(\bu_1-\bu_2,\theta_1-\theta_2)
&& \quad &\text{on}&\quad  \pd \Omega_{5R}.
\end{aligned} \right.
\end{equation*}
Applying Theorem \ref{Stokes:1} and using \eqref{ms:6} yield that
\begin{equation}\label{ms:10}
N_{\lambda}(\bv_1-\bv_2,\omega_1-\omega_2;\Omega_{5R})
 \leq C(M/2) \sigma 
 N_{\lambda}(\bu_1-\bu_2,\theta_1-\theta_2;\Omega_{5R}).
\end{equation}
Since the quantity $N_{\lambda}(\bv,\omega;\Omega_{5R})$ gives an
equivalent norm of 
$H^2_p(\Omega_{5R})^N\times H^1_p(\Omega_{5R})$ 
for $\lambda\not=0,$
the solution maping of \eqref{ms:5}: $(\bu,\theta) \to (\bv, \eta)$ is contraction by choosing $CM \sigma \leq 1.$
Thus, there exists a unique 
$(\bv, \omega) \in H^2_p(\Omega_{5R})^N
\times H^1_p(\Omega_{5R})$ 
satisfying equations \eqref{st:1}.
Moreover, by \eqref{ms:8}, this $(\bv, \omega)$ satisfies \eqref{ms:9}.
This completes the proof of Theorem \ref{thm:ms1}.
\end{proof}

In fact, the proof of Theorem \ref{thm:ms1} also implies the following result.
\begin{cor}\label{cor:ms1}
Let $\bh\in W^{1-1/p}_p(\pd \Omega_{5R}).$ Under the assumptions in Theorem \ref{thm:ms1}, the system
\begin{equation*}
\left\{  \begin{aligned}
&\gamma_1 \lambda \bv
-\ol\DV \big( \mu \ol\bD(\bv)
-\gamma_2\omega\bI\big)= \bff
&&\quad &\text{in}& \quad \Omega_{5R}, \\
&\overline{\dv }\bv = 0
&&\quad &\text{in}& \quad \Omega_{5R}, \\
&\big(\mu \ol\bD(\bv)- \gamma_2\omega\bI\big)
\ol\bn_{_{\pd\Omega_{5R}}}= \bh
&& \quad &\text{on}&\quad  \pd \Omega_{5R},
\end{aligned} \right.
\end{equation*}
admits a unique solution 
$(\bv,\omega) \in H^2_p(\Omega_{5R})^N 
\times  H^1_p(\Omega_{5R})$
possessing the estimate:
\begin{align*}
|\lambda|\|\bv\|_{L_p(\Omega_{5R})} 
+|\lambda|^{1/2}\|\bv\|_{H^1_p(\Omega_{5R})}  
+ \|\bv\|_{H^2_p(\Omega_{5R})}
+ \|\omega\|_{H^1_p(\Omega_{5R})}
 \leq C\big( \|\bff\|_{L_p(\Omega_{5R})} 
 +\|\bh\|_{W^{1-1/p}_p(\pd \Omega_{5R})}\big),
\end{align*}
with some constant $C$ depending on 
$\mu$, $\varepsilon$, $\lambda_0$, $p$ and $\Omega_{5R}.$ 
\end{cor}

\subsubsection{Reduced Stokes problem from \eqref{st:1}} 
To solve \eqref{st:1} with $\lambda=0$, we introduce some reduced Stokes problem associated to \eqref{st:1}.  
First, we review the reduced problem for the standard Stokes system.
By the transformation \eqref{def:LL_w}, we set that
\begin{gather*}
\Omega_{R,T} =\{x = X_{\bw}(y,T) 
\mid y\in\Omega_{5R} \}, \quad
\Gamma_{R,T} = \{x = X_{\bw}(y,T) \mid y\in \pd\Omega_{5R}\}, \\
\bu(x) = \bv\big(X_{\bw}^{-1}(x,T) \big), \quad 
\zeta(x) = \omega\big(X_{\bw}^{-1}(x,T) \big), \quad 
\wt\bff(x) = \bff \big(X_{\bw}^{-1}(x,T) \big).
\end{gather*}
Write $\wt\bn$ for the unit outward normal to $\Gamma_{R, T}.$
Then \eqref{st:1} turns to the resolvent problem of the standard Stokes operator in $\Omega_{R,T},$ 
\begin{equation}\label{ms:2}
\left\{  \begin{aligned}
&\gamma_1 \lambda \bu
-\DV_x \big( \mu \bD(\bu)
-\gamma_2\zeta\bI\big)= \wt\bff
&&\quad &\text{in}& \quad \Omega_{R,T}, \\
&\dv_x \bu = 0
&&\quad &\text{in}& \quad \Omega_{R,T}, \\
&\big(\mu\bD(\bu)- \gamma_2\zeta\bI\big)
\wt\bn= 0
&& \quad &\text{on}&\quad  \Gamma_{R,T}.
\end{aligned} \right.
\end{equation}
For \eqref{ms:2}, we introduce the weak Dirichlet problem:
\begin{equation}\label{wd:1}
(\nabla_x \wt u, \nabla_x \varphi)_{\Omega_{R,T}}
= (\wt\bff, \nabla_x\varphi)_{\Omega_{R,T}}
\quad\text{for any $\varphi \in \wh H^1_{p', 0}(\Omega_{R,T})$},
\end{equation}
where $\wt \bff \in L_p(\Omega_{R,T})^N,$ 
$1<p,p' = p/(p-1)<\infty,$ and 
$$\wh H^1_{p,0}(\Omega_{R,T}) 
= \{\psi \in L_{p, {\rm loc}}(\Omega_{R,T}) :
\nabla \psi \in L_p(\Omega_{R, T})^N, \,\,\, 
\psi|_{\Gamma_{R,T}}=0 \}.$$
Note that in view of Poincar\'e's inequality, 
$$\wh H^1_{p,0}(\Omega_{R,T})
= H^1_{p, 0}(\Omega_{R,T})
= \{\psi \in H^1_{p}(\Omega_{R,T})
:\psi|_{\Gamma_{R,T}}=0 \}.$$
\smallbreak

Given $\bu \in H^2_p(\Omega_{R,T})^N,$
let $K(\bu) \in H^1_{p}(\Omega_{R,T}) + 
\wh H^1_{p,0}(\Omega_{R,T})$ 
be a unique solution of the weak Dirichlet problem
\begin{equation}\label{wd:2}
\big(\nabla_x K(\bu), \nabla_x \varphi\big)_{\Omega_{R,T}} 
= \Big( \DV_x \big(\mu \bD(\bu)\big)-\nabla_x\dv_x\bu, 
\nabla_x\varphi \Big)_{\Omega_{R,T}}
\quad\text{for any $\varphi \in \wh H^1_{p',0}(\Omega_{R,T})$}
\end{equation}
subject to
\begin{equation}\label{wd:3}
K(\bu) = <\mu \bD(\bu)\wt\bn, \wt\bn> - \dv_x\bu \quad\text{on $\Gamma_{R, T}$}.
\end{equation}
And then, the reduced Stokes equations corresponding to equations \eqref{ms:2} is the following:
\begin{equation}\label{ms:11}
\left\{  \begin{aligned}
&\gamma_1 \lambda \bu
-\DV_x \big( \mu \bD(\bu)
-K(\bu) \bI\big)= \wt\bff
&&\quad &\text{in}& \quad \Omega_{R,T}, \\
&\big(\mu \bD(\bu)- K(\bu)\bI\big)
\wt\bn= 0
&& \quad &\text{on}&\quad  \Gamma_{R,T}.
\end{aligned} \right.
\end{equation}
For more discussion on \eqref{ms:11}, we refer to the lecture note \cite{Shi2020}.
\medskip

It is not hard to see that the weak Dirichlet problem
\eqref{wd:1} is transformed to the following variational equation:
\begin{equation}\label{wd:4}
(\ol \nabla u, J\ol\nabla \varphi)_{\Omega_{5R}} = (\bff, J\ol\nabla\varphi)_{\Omega_{5R}}
\quad\text{for any $\varphi \in \wh H^1_{p',0}(\Omega_{5R})$},
\end{equation}
with $\ol\nabla =(\bI + \bV)\nabla_y,$ 
$ J= \big(\det(\bI+\bV)\big)^{-1},$ and 
$$\wh H^1_{p,0}(\Omega_{5R}) 
= \{\psi \in L_{p, {\rm loc}}(\Omega_{5R}) :
\nabla \psi \in L_p (\Omega_{5R})^N, \,\,\, 
\psi|_{\pd\Omega_{5R}}=0 \}.$$ 
\eqref{wd:4} can be solved uniquely for any $\bff \in L_p(\Omega_{5R})^N$ ($1<p\leq r<\infty$).  
In fact, employing the standard Banach fixed point theorem as in the proof of Theorem \ref{thm:ms1}
and the unique solvability of the weak Dirichlet problem in $\Omega_{5R}$, we can show the unique existence of the solution to equations \eqref{wd:4}.  
Thus, for any $\bu \in H^2_p(\Omega_{5R})^N,$
let $\bar K(\bu)$ be a unique solution of 
the weak Dirichlet problem:
\begin{equation}\label{wd:5}
\big(\ol\nabla\, \bar K(\bu), J\,\ol\nabla\varphi \big)_{\Omega_{5R}}
= \Big(\ol\DV\big(\mu \ol\bD(\bu)\big) -\ol\nabla \, \overline{\dv}\bu,
 J\ol \nabla \varphi \Big)_{\Omega_{5R}}
\quad\text{for any $\varphi \in \wh H^1_{p',0}(\Omega_{5R})$},
\end{equation}
subject to
\begin{equation}\label{wd:6}
\bar K(\bu) = \fd^{-1}<\mu \ol\bD(\bu)\ol\bn_{_{\pd\Omega_{5R}}}, 
\ol\bn_{_{\pd\Omega_{5R}}}> -\overline{\dv}\bu
\quad\text{on $\pd\Omega_{5R}$},
\end{equation}
where we have set 
$\fd = |\ol\bn_{_{\pd\Omega_{5R}}}|^2=|(\bI+\bV)\bn_{_{\pd\Omega_{5R}}}|^2.$
By \eqref{def:LL_w}, the reduced Stokes equations
\eqref{ms:11} are transformed to the following equations: 
\begin{equation}\label{st:2}
\left\{  \begin{aligned}
&\gamma_1 \lambda \bv
-\ol\DV \big( \mu \ol\bD(\bv)
-\bar K(\bv) \bI\big)= \Bf
&&\quad &\text{in}& \quad \Omega_{5R}, \\
&\big(\mu\ol\bD(\bv)-\bar K(\bv)\bI\big)
\ol\bn_{_{\pd\Omega_{5R}}}= 0
&& \quad &\text{on}&\quad  \pd\Omega_{5R}.
\end{aligned} \right.
\end{equation}
Let $\bar J_p(\Omega_{5R})$ be a modified solenoidal space defined by
 \begin{equation}\label{solenoid:1}
 \bar J_p(\Omega_{5R}) = \{\bg\in L_p(\Omega_{5R})^N \mid 
 (\bg, J\ol\nabla\varphi)_{\Omega_{5R}}
 = 0 \quad\text{for any $\varphi \in \wh H^1_{p',0}(\Omega_{5R})$}\}. 
\end{equation}
Assume that $(\bv,\omega) \in H^2_p(\Omega_{5R})^N \times H^1_p(\Omega_{5R})$ is a solution of equations \eqref{st:1} for some $\bff \in \bar J_p(\Omega_{5R})$ by Theorem \ref{thm:ms1}.
Then $\bv$ also satisfies \eqref{st:2}. In fact, it suffices to verify
$ \bar K(\bv)=\gamma_2\omega.$
For any $\varphi \in \wh H^1_{p',0}(\Omega_{5R})$, by \eqref{st:1} and 
the solenoidality of $\bff$, we have
\begin{align*}
0 = (\bff, J\ol\nabla \varphi)_{\Omega_{5R}} 
= &\gamma_1\lambda(\bv, J\,\ol\nabla \varphi)_{\Omega_{5R}}
-\Big(\ol\DV\big(\mu\ol\bD(\bv)\big) 
- \ol\nabla\,\overline{\dv}\bv, J\,\ol\nabla\varphi)_{\Omega_{5R}}
\\
&-(\ol\nabla\,\overline{\dv}\bv, J\,\ol\nabla\varphi)_{\Omega_{5R}}
+ \big(\ol\nabla (\gamma_2 \omega), J\,\ol\nabla \varphi \big)_{\Omega_{5R}}\\
=&\gamma_1\lambda(\bv, J\,\ol\nabla \varphi)_{\Omega_{5R}}
+\Big(\ol\nabla \big(\gamma_2 \omega-\bar K(\bv)\big), J\,\ol\nabla \varphi \Big)_{\Omega_{5R}}.
\end{align*} 
By integration by parts, $\varphi|_{\pd\Omega_{5R}}=0$ and $\overline{\dv}\bv = 0,$  we obtain that
\begin{equation}\label{eq:div_solen}
(\bv, J\,\ol\nabla \varphi)_{\Omega_{5R}} 
= -\Big(\dv \big((\bI+\bV)^\top\bv\big),   \varphi\Big)_{\Omega_{5R}}
=(J\overline{\dv}\bv, \varphi)_{\Omega_{5R}}=0.
\end{equation}
Then we have
$$\Big(\ol\nabla \big(\gamma_2 \omega-\bar K(\bv)\big), J\,\ol\nabla \varphi \Big)_{\Omega_{5R}}=0
\quad\text{for any $\varphi \in \wh H^1_{p',0}(\Omega_{5R})$}.
$$
From the boundary condition in \eqref{st:1}, it follows that
$$\gamma_2\omega = \fd^{-1}<\ol\bD(\bv)\ol\bn_{_{\pd\Omega_{5R}}}, 
\ol\bn_{_{\pd\Omega_{5R}}}> - \overline{\dv}\bv = \bar K(\bv)
\quad\text{on $\pd\Omega_{5R}$},
$$
where we have used $\overline{\dv}\bv=0$.  Thus, the uniqueness of the 
weak Dirichlet problem yields that $\gamma_2\omega = \bar K(\bv)$, which shows
that $\bv \in H^2_p(\Omega_{
5R})^N$ solves equations \eqref{st:2}.  
Moreover, thanks to \eqref{eq:div_solen},
$\overline{\dv} \bv=0$ yields that
$\bv \in \bar J_p(\Omega_{5R}).$ 

From these observations and Theorem \ref{thm:ms1}, we have the following 
theorem.
\begin{thm}\label{thm:ms2} 
Let $1 < p \leq r,$ $0 < \varepsilon < \pi/2$, and $\lambda_0 > 0$.
Then, there exists a $\sigma_1$ depending on 
$\mu$, $\varepsilon$, $\lambda_0$, $p$, $r$
and $\Omega_{5R}$ such that if $0 < \sigma \leq \sigma_1$,
then for any $\lambda \in \Sigma_{\varepsilon, \lambda_0}$,
 $\bff \in \bar J_p(\Omega_{5R})$, problem \eqref{st:2} admits a
unique solution 
$\bv \in H^2_p(\Omega_{5R})^N\cap \bar J_p(\Omega_{5R})$ 
possessing the estimate:
\begin{align*}
|\lambda|\|\bv\|_{L_p(\Omega_{5R})} 
+|\lambda|^{1/2}\|\bv\|_{H^1_p(\Omega_{5R})}  
+ \|\bv\|_{H^2_p(\Omega_{5R})}
 \leq C\|\bff\|_{L_p(\Omega_{5R})}
\end{align*}
with some constant $C$ depending on $\mu$, $\varepsilon$, $\lambda_0$, 
$p,$ $r$ and $\Omega_{5R}$. 
\end{thm}

\subsubsection{Analysis of \eqref{st:1} for $\lambda=0$} 
To study \eqref{st:1} for the parameter $\lambda=0,$ we consider the following reduced model: 
\begin{equation}\label{st:3}
\left\{  \begin{aligned}
&-\ol\DV \big( \mu \ol\bD(\bv)
-\bar K(\bv) \bI\big)= \Bf
&&\quad &\text{in}& \quad \Omega_{5R}, \\
&\big(\mu\ol\bD(\bv)-\bar K(\bv)\bI\big)
\ol\bn_{_{\pd\Omega_{5R}}}= 0
&& \quad &\text{on}&\quad  \pd\Omega_{5R}.
\end{aligned} \right.
\end{equation}
First, we introduce the rigid motion associated to \eqref{ms:2}.  Set
\begin{equation}\label{eq:rid}
\br_{j}(x)=\begin{cases}
\be_{j}= (0, \ldots, 
\underbrace{1}_{\mbox{\tiny ${\rm j}$th component}},
\ldots, 0) & \mbox{for}\,\,\, j=1, \ldots, N,\\
x_k\be_{\ell}-x_{\ell}\be_k \,\,\,
(k,\ell=1,\ldots,N)& \mbox{for}\,\,\, j=N+1, \ldots, M.
\end{cases}
\end{equation}
Above, $M$ is a constant only depending of the dimension $N.$
For any vector $\bu$ satisfying
$\bD(\bu) = 0,$ $\bu$ is represented by a linear combination of
$\{\br_j\}_{j=1}^M$, namely $\bu = \sum_{j=1}^M a_j\br_j$ with some 
$a_j \in \BR.$ For convenience, we denote $\{\bp_j\}_{j=1}^M$ for the orthonormalization of $\{\br_j\}_{j=1}^M,$ that is,  $\{\bp_j\}_{j=1}^M$
satisfies
\begin{equation}\label{orth:3.9}
(\bp_j, \bp_k)_{\Omega_{R,T} }=\delta_{jk}
\quad\text{for $j, k = 1, \ldots, M$}.
\end{equation}
By the transformation \eqref{def:LL_w}, we define $\bar\bp_j(y) = \bp_j(x)$ for $j=1,\ldots,M,$  which, thanks to \eqref{orth:3.9}, satisfy the orthogonal conditions:
\begin{equation}\label{orth:4}
(\bar\bp_j, J\bar\bp_k)_{\Omega_{5R} }=\delta_{jk}
\quad\text{for $j, k = 1, \ldots, M$}.
\end{equation}
In addition,
$\ol\bD(\bar\bp_j) = \bD(\bp_j) = 0$ and $\overline{\dv} \bar\bp_j
=\dv\bp_j=0.$ 
For convenience, we introduce the spaces 
\begin{equation}\label{orth:2}
\begin{aligned}
\bar\CR_d(\Omega_{5R}) &= \Span \{\bar\bp_1,\ldots,\bar\bp_M\},\\
\bar\CR_d(\Omega_{5R})^\perp &= \{ \bv \in L_1(\Omega_{5R})^N \mid 
(\bv, J\bar\bp_j)_{\Omega_{5R}} = 0 
\,\,\,\text{for $j=1, \ldots, M$}\}.
\end{aligned}
\end{equation}
\smallbreak

Now, let $\bv\in H^2_p(\Omega_{5R})^N \cap \bar J_p(\Omega_{5R})$ be a solution of \eqref{st:2}
for some $\bff \in \bar J_p(\Omega_{5R}) \cap 
 \bar \CR_d(\Omega_{5R})^\perp$ and $\lambda \in \Sigma_{\ep,\lambda_0}$ by Theorem \ref{thm:ms2}.
Then, $\bv \in \bar \CR_d(\Omega_{5R})^\perp$ as well.  
In fact, 
\begin{equation}\label{eq:vp}
\begin{aligned}
0 = (\bff, J\bar\bp_j)_{\Omega_{5R}} &= 
\gamma_1\lambda(\bv, J\bar\bp_j)_{\Omega_{5R}} 
-\Big(\ol\DV\big(\mu \ol\bD(\bv)
-\bar K(\bv)\bI\big), J\bar\bp_j\Big)_{\Omega_{5R}}
\\
& = \gamma_1\lambda(\bv, J\bar\bp_j)_{\Omega_{5R}} 
-\Big(\big(\mu \ol\bD(\bv)-\bar K(\bv)\bI\big)\ol\bn_{_{\pd\Omega_{5R}}}, J\bar\bp_j\Big)_{\pd\Omega_{5R}}
\\
&\quad + \int_{\Omega_{5R}} 
\big(\mu \ol\bD(\bv)-\bar K(\bv)\bI\big)(\bI+\bV):
\nabla_y\bar\bp_j \,J\,dy 
 \\
& = \gamma_1\lambda(\bv, J\bar\bp_j)_{\Omega_{5R}} 
+\frac{\mu}{2}
\int_{\Omega_{5R}} \ol\bD(\bv):\ol \bD(\bar\bp_j) \,J\,dy
- (J\bar K(\bv), \overline{\dv} \bar \bp_j)_{\Omega_{5R}} 
\\
&= \gamma_1\lambda(\bv, J\bar\bp_j)_{\Omega_{5R}}. 
\end{aligned}
\end{equation}

Under these preparations, we have the following theorem. 
\begin{thm}\label{thm:ms3} 
Let $1 < p \leq r$. Then, for any 
$\bff \in \bar J_p(\Omega_{5R}) \cap \bar\CR_d(\Omega_{5R})^{\perp}$, 
problem \eqref{st:3} admits a unique solution 
$\bv \in H^2_p(\Omega_{5R})^N \cap \bar J_p(\Omega_{5R}) \cap 
\bar\CR_d(\Omega_{5R})^{\perp}$
possessing the estimate:
$$\|\bv\|_{H^2_p(\Omega_{5R})} \leq C\|\bff\|_{L_p(\Omega_{5R})}.
$$
\end{thm}
\begin{proof}
In view of Theorem \ref{thm:ms2} and the Fredholm alternative principle, to prove
Theorem \ref{thm:ms3}, it is enough to verify the uniqueness in the $L_2(\Omega_{5R})^N$
framework. Namely, if  
$\bv \in H^2_2(\Omega_{5R})^N \cap J_2(\Omega_{5R}) 
\cap \bar \CR_d(\Omega_{5R})^{\perp} $ 
satisfies the homogeneous equations:
 \begin{equation*}
\left\{  \begin{aligned}
&-\ol\DV \big( \mu \ol\bD(\bv)
-\bar K(\bv) \bI\big)= 0
&&\quad &\text{in}& \quad \Omega_{5R}, \\
&\big(\mu \ol\bD(\bv)-\bar K(\bv)\bI\big)
\ol\bn_{_{\pd\Omega_{5R}}}= 0
&& \quad &\text{on}&\quad  \pd\Omega_{5R}.
\end{aligned} \right.
\end{equation*} 
 then we shall show that $\bv = 0$.  By the similar argument to \eqref{eq:vp}, we have 
 \begin{align*}
 0 =  \Big( -\ol\DV \big(\mu \ol\bD(\bv) - \bar K(\bv)\bI \big), J \bv\Big)_{\Omega_{5R}}
  =  \frac{\mu}{2} \int_{\Omega_{5R}} \ol\bD(\bv): \ol\bD(\bv) \,Jdy.
 \end{align*}
 Therefore, $\ol\bD(\bv)=0$ in $\Omega_{5R}$, and then $\bv$ is a linear combination of $\{\bar\bp_j\}_{j=1}^M.$
 Since $\bv \in  \bar \CR_d(\Omega_{5R})^{\perp} $, 
 we have $\bv=0$. 
 This completes the proof of Theorem \ref{thm:ms3}. 
\end{proof}

We now consider the Stokes equation:
\begin{equation}\label{st:4}
\left\{  \begin{aligned}
&-\ol\DV \big( \mu \ol\bD(\bv)
-\gamma_2\omega\bI\big)= \bff
&&\quad &\text{in}& \quad \Omega_{5R}, \\
&\overline{\dv }\bv = 0
&&\quad &\text{in}& \quad \Omega_{5R}, \\
&\big(\mu \ol\bD(\bv)- \gamma_2\omega\bI\big)
\ol\bn_{_{\pd\Omega_{5R}}}= 0
&& \quad &\text{on}&\quad  \pd \Omega_{5R},
\end{aligned} \right.
\end{equation}
From Theorem \ref{thm:ms3}, we have the following corollary.
 \begin{cor}\label{cor:st1}
 Let $1 < p \leq r$. 
Then, for any $\bff \in \bar J_p(\Omega_{5R}) \cap 
 \bar\CR_d(\Omega_{5R})^{\perp},$
problem \eqref{st:4} admits a unique solution 
$\bv \in H^2_p(\Omega_{5R})^N \cap \bar J_p(\Omega_{5R}) 
\cap \bar \CR_d(\Omega_{5R})^{\perp} $
and $\omega \in H^1_p(\Omega_{5R})$ 
possessing the estimate:
$$\|\bv\|_{H^2_p(\Omega_{5R})} 
+ \|\omega\|_{H^1_p(\Omega_{5R})}
\leq C\|\bff\|_{L_p(\Omega_{5R})}.
$$
\end{cor}
\begin{proof}
Let $\bv \in H^2_p(\Omega_{5R})^N \cap \bar J_p(\Omega_{5R}) 
\cap \bar \CR_d(\Omega_{5R})^{\perp}$
be a solution of equations \eqref{st:3}, whose existence is guaranteed by
Theorem \ref{thm:ms3}. For any 
$\varphi \in \wh H^1_{p',0}(\Omega_{5R})$, 
\begin{align*}
0 = (\bff, J\,\ol\nabla \varphi)_{\Omega_{5R}} 
&= -\Big( \ol\DV \big( \mu \ol\bD(\bv)\big), 
J\,\ol\nabla \varphi \Big)_{\Omega_{5R}}
+\big( \ol\nabla \bar K(\bv), J\,\ol\nabla\varphi \big)_{\Omega_{5R}}\\
&= -(\ol\nabla\, \overline{\dv}\bv, J\,\ol\nabla\varphi)_{\Omega_{5R}}.
\end{align*}
Moreover, by the boundary condition in \eqref{st:3}, we have
\begin{align*}
\overline{\dv}\bv 
&= \fd^{-1}<\mu \ol\bD(\bv) \ol\bn_{_{\pd\Omega_{5R}}},
\ol\bn_{_{\pd\Omega_{5R}}}> -\bar K(\bv) \\
&= \fd^{-1}<\big(\mu \ol\bD(\bv) - \bar K(\bv)\big)
\ol\bn_{_{\pd\Omega_{5R}}},
\ol\bn_{_{\pd\Omega_{5R}}}> =0\quad\text{on $\pd\Omega_{5R}$}.
\end{align*}
Then, the uniqueness of the weak Dirichlet problem yields that 
$\overline{\dv}\bv=0$ in $\Omega_{5R}.$ 
Thus, $\bv$ automatically 
satisfies equations \eqref{st:4} with imposing $\omega =\gamma_2^{-1}\bar K(\bv),$
which completes the proof of Corollary \ref{cor:st1}. 
\end{proof}

\subsection{Modified compressible model problem}
In this subsection, we consider the following problem:
\begin{equation}\label{eq:4.7}
\left\{  \begin{aligned}
&\gamma_1\ol\dv \bv= d
&&\quad &\text{in}& \quad \Omega_{5R}, \\
&-\ol\DV \big( \ol\bS(\bv)
-\gamma_2\omega\bI\big)= \bff
&&\quad &\text{in}& \quad \Omega_{5R}, \\
&\big(\ol\bS(\bv)- \gamma_2\omega\bI\big)
\ol\bn_{_{\pd\Omega_{5R}}}= \bh
&& \quad &\text{on}&\quad  \pd \Omega_{5R},
\end{aligned} \right.
\end{equation}
with $\ol\bS(\bv)=\mu \ol\bD(\bv)+(\nu-\mu) (\ol\dv \bv) \bI.$
Concerning \eqref{eq:4.7}, we have the following theorem. 
\begin{thm}\label{thm:st2} 
Let $1 < p \leq r$. 
Then, for any 
$d \in H^1_p(\Omega_{5R})$,  $\bff \in L_p(\Omega_{5R})^N$ 
and $\bh \in H^1_p(\Omega_{5R})^N$
satisfying the orthogonal condition:
\begin{equation}\label{orth:3}
(\bff, J\bar\bp_j)_{\Omega_{5R}} 
+ (\bh, J\bar\bp_j)_{\pd \Omega_{5R}} = 0
\quad (j=1, \ldots, M), 
\end{equation}
problem \eqref{eq:4.7}
admits a unique solution 
$\bv \in H^2_p(\Omega_{5R})^N \cap \bar\CR_d(\Omega_{5R})^\perp$
and $\omega \in H^1_p(\Omega_{5R})$ possessing the estimate:
$$\|\bv\|_{H^2_p(\Omega_{5R})} 
+ \|\omega\|_{H^1_p(\Omega_{5R})}
\leq C\big(\|\bff\|_{L_p(\Omega_{5R})} 
+\|(d,\bh)\|_{H^1_p(\Omega_{5R})} \big). 
$$
\end{thm}
\begin{proof}
Let $\bv_1 \in H^2_p(\Omega_{5R})^N$ be a solution of the divergence equation:
$\gamma_1\overline{\dv}\bv_1 = d$ in $\Omega_{5R}$ possessing the estimate:
$$\|\bv_1\|_{H^2_p(\Omega_{5R})} 
\leq C\|d\|_{H^1_p(\Omega_{5R})}.$$ 
The existence of such $\bv_1$ is guaranteed by 
Proposition \ref{prop:div}. 
Therefore, to construct the solution of \eqref{eq:4.7}, it suffices to solve the following system:
\begin{equation}\label{eq:uomega}
\left\{  \begin{aligned}
&\ol\dv \bu= 0
&&\quad &\text{in}& \quad \Omega_{5R}, \\
&-\ol\DV \big( \mu \ol\bD(\bu)
-\gamma_2\omega\bI\big)= \bff+F_1
&&\quad &\text{in}& \quad \Omega_{5R}, \\
&\big(\mu \ol\bD(\bu)- \gamma_2\omega\bI\big)
\ol\bn_{_{\pd\Omega_{5R}}}= \bh+H_1
&& \quad &\text{on}&\quad  \pd \Omega_{5R},
\end{aligned} \right.
\end{equation}
with $F_1 = \ol\DV\big(\ol \bS(\bv_1)\big),$ and 
$H_1=-\ol \bS(\bv_1) \ol\bn_{_{\pd\Omega_{5R}}}.$
Then $\bv = \bv_1 + \bu$ and $\omega$ satisfy the equations \eqref{eq:4.7}. 
Moreover, the divergence theorem of Gau\ss\, yields that
\begin{equation}\label{orth:3.5}
(F_1, J\bar\bp_j)_{\Omega_{5R}} 
+ (H_1, J\bar\bp_j)_{\pd\Omega_{5R}}
=-\frac{1}{2} \int_{\Omega_{5R}} \ol\bS(\bv_1):
\ol\bD(\bar\bp_j) J\,dy=0.
\end{equation}

To obtain the solution of \eqref{eq:uomega}, we consider the following two linear systems:
\begin{equation}\label{eq:4.75}
\left\{  \begin{aligned}
&\ol\dv \bv_2= 0
&&\quad &\text{in}& \quad \Omega_{5R}, \\
&\gamma_1 \lambda_0 \bv_2-\ol\DV \big( \mu \ol\bD(\bv_2)
-\gamma_2\omega_2\bI\big)= \bff+F_1
&&\quad &\text{in}& \quad \Omega_{5R}, \\
&\big(\mu \ol\bD(\bv_2)- \gamma_2\omega_2\bI\big)
\ol\bn_{_{\pd\Omega_{5R}}}= \bh+H_1
&& \quad &\text{on}&\quad  \pd \Omega_{5R},
\end{aligned} \right.
\end{equation}
\begin{equation}\label{eq:4.8}
\left\{  \begin{aligned}
&\ol\dv \bv_3= 0
&&\quad &\text{in}& \quad \Omega_{5R}, \\
&-\ol\DV \big( \mu \ol\bD(\bv_3)
-\gamma_2\omega_3\bI\big)= \gamma_1 \lambda_0 \bv_2
&&\quad &\text{in}& \quad \Omega_{5R}, \\
&\big(\mu \ol\bD(\bv_3)- \gamma_2\omega_3\bI\big)
\ol\bn_{_{\pd\Omega_{5R}}}= 0
&& \quad &\text{on}&\quad  \pd \Omega_{5R},
\end{aligned} \right.
\end{equation}
where  $\lambda_0>0$ is fixed large number. At least formally, $(\bu,\omega)=(\bv_2,\omega_2)+(\bv_3,\omega_3)$ gives a solution of \eqref{eq:4.7}.
Thanks to Corollary \ref{cor:ms1}, there exists a solution of \eqref{eq:4.75} such that 
\begin{align*}
\|\bv_2\|_{H^2_p(\Omega_{5R})} + \|\omega_2\|_{H^1_p(\Omega_{5R})}
&\leq C(\|\bff + F_1\|_{L_p(\Omega_{5R})} 
+ \|\bh + H_1\|_{H^1_p(\Omega_{5R})})\\
&\leq C(\|\bff\|_{L_p(\Omega_{5R})} 
+ \|(d, \bh)\|_{H^1_p(\Omega_{5R})}).
\end{align*}
Additionally, $\bv_2 \in \bar\CR_d(\Omega_{5R})^\perp.$ 
In fact, \eqref{orth:3} and \eqref{orth:3.5} yield that
\begin{align*} 
\gamma_1\lambda_0 (\bv_2, J\bar\bp_j)_{\Omega_{5R}}
&= \Big(\ol\DV \big( \mu \ol\bD(\bv_2)
-\gamma_2\omega_2\bI\big), J\bar\bp_j\Big)_{\Omega_{5R}}
+(\bff + F_1, J\bar\bp_j)_{\Omega_{5R}}\\
 & = (\bh + H_1, J\bar\bp_j)_{\pd\Omega_{5R}} 
 + (\bff + F_1, J\bar\bp_j)_{\Omega_{5R}}=0.
\end{align*}

On the other hand, to solve \eqref{eq:4.8}, we notice that 
$\gamma_1 \lambda_0 \bv_2 \in\bar J_p(\Omega_{5R})$ from the divergence equation, as
\begin{equation*}
(\bv_2, J\,\ol\nabla \varphi )_{\Omega_{5R}} = (J\varphi \bv_2, \ol \bn_{\pd \Omega_{5R}})_{\pd \Omega_{5R}} 
- (\overline{\dv} \bv_2, J\varphi )_{\Omega_{5R}}=0.
\end{equation*}
Then, by Corollary \ref{cor:st1}, there exists a solution 
$\bv_3 \in H^2_p(\Omega_{5R})^N 
\cap \bar\CR_d(\Omega_{5R})^\perp$
and $\omega_3 \in H^1_p(\Omega_{5R})$
of the system \eqref{eq:4.8} satisfying
$$\|\bv_3\|_{H^2_p(\Omega_{5R})} 
+ \|\omega_3\|_{H^1_p(\Omega_{5R})}
\leq C\|\gamma_1 \lambda_0 \bv_2\|_{L_p(\Omega_{5R})}
\leq C(\|\bff\|_{L_p(\Omega_{5R})} 
+ \|(d, \bh)\|_{H^1_p(\Omega_{5R})}).$$
Combing the discussion above, we find one solution $(\bv,\omega)$ of \eqref{eq:4.7}  by defining
$$(\bv,\omega)=(\bv_1+\bv_2+\bv_3, \omega_2+\omega_3) 
\in H^2_p(\Omega_{5R})^N  \times H^1_p(\Omega_{5R}).$$ 
Moreover, we have
$$\|\bv\|_{H^2_p(\Omega_{5R})}
 + \|\omega\|_{H^1_p(\Omega_{5R})}
\leq C\big(\|(f,\bh)\|_{H^1_p(\Omega_{5R})} 
+ \|\bg\|_{L_p(\Omega_{5R})}\big).$$

To meet the orthogonal condition, i.e.
$\bv \in \CR_d(\Omega_{5R})^\perp,$
we refine the definition of the velocity field by
$$\wt \bv = \bv 
- \sum_{j=1}^M (\bv, J\bar\bp_j)_{\Omega_{5R}} 
\bar\bp_j  \in \bar\CR_d(\Omega_{5R})^\perp.
$$
Since $\ol\bD(\bar\bp_j) = 0$ and $\overline{\dv}\bar\bp_j = 0$,
we see that $\wt\bv$ and $\omega$ satisfy equations \eqref{eq:4.7} with the desired estimate.
\medskip

We now prove the uniqueness.  Let $\bv \in H^2_p(\Omega_{5R})^N 
\cap \bar\CR_d(\Omega_{5R})^\perp$ and $\omega \in H^1_p(\Omega_{5R})$ for $2\leq p\leq r$ 
satisfy the homogeneous equations: 
\begin{equation}\label{eq:4.9}
\left\{  \begin{aligned}
&\ol\dv \bv= 0
&&\quad &\text{in}& \quad \Omega_{5R}, \\
&-\ol\DV \big( \mu \ol\bD(\bv)
-\gamma_2\omega\bI\big)= 0
&&\quad &\text{in}& \quad \Omega_{5R}, \\
&\big(\mu \ol\bD(\bv)- \gamma_2\omega\bI\big)
\ol\bn_{_{\pd\Omega_{5R}}}= 0
&& \quad &\text{on}&\quad  \pd \Omega_{5R},
\end{aligned} \right.
\end{equation}
By the previous discussion on the existence issue of \eqref{eq:4.7},
given any $\bg \in L_{2}(\Omega_{5R})^N \cap \bar \CR_d(\Omega_{5R})^{\perp}$, 
there exists $(\bu,\theta) \in \big( H^2_{2}(\Omega_{5R})^N \cap \bar \CR_d(\Omega_{5R})^{\perp} \big) \times H^1_{2}(\Omega_{5R})$
fulfilling the equations:
\begin{equation*}
\left\{  \begin{aligned}
&\ol\dv \bu= 0
&&\quad &\text{in}& \quad \Omega_{5R}, \\
&-\ol\DV \big( \mu \ol\bD(\bu)
-\gamma_2\theta\bI\big)= \bg
&&\quad &\text{in}& \quad \Omega_{5R}, \\
&\big(\mu \ol\bD(\bu)- \gamma_2\theta\bI\big)
\ol\bn_{_{\pd\Omega_{5R}}}= 0
&& \quad &\text{on}&\quad  \pd \Omega_{5R}.
\end{aligned} \right.
\end{equation*}
By the divergence theorem of Gau\ss, we have
\begin{align*}
0&= \Big(\ol\DV\big(\mu \ol\bD(\bv)-\gamma_2\omega\bI\big), 
J\bu\Big)_{\Omega_{5R}} \\
&= \Big( \big(\mu \ol\bD(\bv)- \gamma_2\omega\bI\big)\ol\bn_{_{\Omega_{5R}}}, J\bu\Big)_{\pd\Omega_{5R}}
+\frac{\mu }{2}\int_{\Omega_{5R}} \ol\bD(\bv): \ol\bD(\bu) J dy
+ (J\gamma_2\omega, \overline{\dv}\bu)_{\Omega_{5R}} \\
&= \frac{\mu}{2}\int_{\Omega_{5R}} \ol\bD(\bv): \ol\bD(\bu) J dy\\
& = \Big( J\bv, \ol\DV\big(\mu \ol\bD(\bu)- \gamma_2\theta\bI\big) \Big)_{\Omega_{5R}}
= (J\bv, \bg)_{\Omega_{5R}}. 
\end{align*}

For any $\bff \in L_{2}(\Omega_{5R})^N,$ we set
$$\bff_{\perp} = \bff - \sum_{j=1}^M(\bff, J\bar\bp_{j})_{\Omega_{5R}}
\bar\bp_j\in L_{2}(\Omega_{5R})^N \cap \bar\CR_d(\Omega_{5R})^{\perp}.$$
Then we infer from $\bv \in \bar\CR_d(\Omega_{5R})^\perp$ that
\begin{equation*}
 (J\bv, \bff)_{\Omega_{5R}} = 
(J\bv, \bff_{\perp})_{\Omega_{5R}}
+ \sum_{j=1}^M(\bff, J\bar\bp_{j})_{\Omega_{5R}} (J\bv,\bar\bp_j)_{\Omega_{5R}}=(J\bv, \bff_{\perp})_{\Omega_{5R}}=0.
\end{equation*}
The arbitrary choice of $\bff$ yields that $J\bv=0$, and thus $\bv=0.$
\smallbreak

Now, the equations \eqref{eq:4.9} are reduced to
the form:
\begin{equation*}
\left\{  \begin{aligned}
&\ol\nabla \omega =(\bI+\bV) \nabla_y \omega =0
&&\quad &\text{in}& \quad \Omega_{5R}, \\
&\omega\ol\bn_{_{\pd\Omega_{5R}}} 
=\omega(\bI+\bV)\bn_{_{\pd\Omega_{5R}}}= 0
&& \quad &\text{on}&\quad  \pd \Omega_{5R}.
\end{aligned} \right.
\end{equation*}
Therefore, $\nabla_y \omega=0$ in $\Omega_{5R},$ and then $\omega=0$ by the boundary condition. 
\smallbreak

Next, we assume that 
$(\bv,\omega) \in \big( H^2_{p}(\Omega_{5R})^N \cap \bar \CR_d(\Omega_{5R})^{\perp} \big) \times H^1_{p}(\Omega_{5R})$
for some $1<p<2$ satisfies the equations \eqref{eq:4.9}.
Notice that the embedding $H^1_p(\Omega_{5R}) \hookrightarrow L_q(\Omega_{5R})$ with $0<N(1/p-1/q)<1.$
Then, by Theorem \ref{thm:ms1}, there exists 
$(\bu,\theta) \in  
H^2_{q}(\Omega_{5R})^N \times H^1_{q}(\Omega_{5R})$ 
satisfying 
\begin{equation*}
\left\{  \begin{aligned}
&\gamma_1 \lambda \bu
-\ol\DV \big( \mu \ol\bD(\bu)
-\gamma_2\theta \bI\big)= \gamma_1 \lambda \bv
&&\quad &\text{in}& \quad \Omega_{5R}, \\
&\overline{\dv }\bu = 0
&&\quad &\text{in}& \quad \Omega_{5R}, \\
&\big(\mu \ol\bD(\bu)- \gamma_2\theta\bI\big)
\ol\bn_{_{\pd\Omega_{5R}}}= 0
&& \quad &\text{on}&\quad  \pd \Omega_{5R},
\end{aligned} \right.
\end{equation*}
for some $\lambda>0.$
In fact, we have $(\bu,\theta)=(\bv,\omega) \in H^2_{q}(\Omega_{5R})^N \times H^1_{q}(\Omega_{5R})$ by the uniqueness in Theorem \ref{thm:ms1}.
If $q\geq 2,$ we have
$(\bv,\omega) \in H^2_{2}(\Omega_{5R})^N \times H^1_{2}(\Omega_{5R})$ and the uniqueness of \eqref{eq:4.9} from the previous discussion on the case $p \geq 2$.
Otherwise, we repeat such argument in finite times.
This completes the proof of Theorem \ref{thm:st2}.
\end{proof}


\section{Resolvent problem for $\lambda$ near zero}
\label{sec:near}

In this section, we will study the behaviour of the solution of the system \eqref{resolvent_0} whenever $\lambda$ lies near the origin.
The main result reads as follows.
\begin{thm}\label{thm:ext_1}
Let $(d,\bff) \in X_{p,L}(\Omega)$ for $1<p\leq r$ and $L>2R>0.$
Then there exist a constant $\lambda_1>0$ and two families of the operators $(\BM_{\lambda},\BV_{\lambda})$ for any 
$\lambda \in\dot{U}_{\lambda_1}=\{\lambda \in 
\BC \backslash (-\infty, 0]: |\lambda| <\lambda_1\}$ 
with 
\begin{align*}
\BM_{\lambda} \in \Hol\Big(\dot{U}_{\lambda_1}; 
\CL\big( X_{p,L}(\Omega) ;
H^{1}_{p,\loc}(\Omega)\big)\Big),\quad
\BV_{\lambda} \in \Hol\Big(\dot{U}_{\lambda_1};
\CL\big( X_{p,L}(\Omega) ;
H^{2}_{p,\loc}(\Omega)^N\big)\Big),
\end{align*}
so that $(\eta,\bu)=(\BM_{\lambda},\BV_{\lambda})(d,\bff)$ solves \eqref{resolvent_0}.
Moreover, there exist families of the operators
\begin{align*}
\BM_{\lambda}^i &\in \Hol\Big(\dot{U}_{\lambda_1}; 
\CL\big( X_{p,L}(\Omega) ;
H^{1}_{p,\loc}(\Omega)\big)\Big)\,\,\,(i=1,2),\\
\BV_{\lambda}^j &\in \Hol\Big(\dot{U}_{\lambda_1}; 
\CL\big( X_{p,L}(\Omega) ;
H^{2}_{p,\loc}(\Omega)^N\big)\Big)\,\,\,(j=0,1,2),
\end{align*}
such that 
\begin{align*}
\BM_{\lambda} &= (\lambda^{N-2} \log \lambda) \BM_{\lambda}^1  + \BM_{\lambda}^2, \\
\BV_{\lambda} &= \big(\lambda^{N\slash 2 -1} (\log \lambda)^{\sigma(N)} \big) \BV_{\lambda}^0
+ (\lambda^{N-2} \log \lambda) \BV_{\lambda}^1 +\BV_{\lambda}^2,
\end{align*}
for any $\lambda \in \dot{U}_{\lambda_1}$ and 
$\sigma (N)=\big((-1)^N+1 \big)\slash 2.$
\end{thm}

\subsection{Resolvent problem in $\BR^N$}
In this subsection, we review the result of some model problem in 
$\BR^N:$ 
\begin{equation}\label{eq:resolvent_whole}
	\left\{\begin{aligned}
&\lambda \eta  + \gamma_1 \di \bu = d  
	    &&\quad\hbox{in}\quad \BR^N, \\
&\gamma_1\lambda \bu-\Di\big( \bS(\bu) -\gamma_2 \eta \bI \big)=\bff
		&&\quad\hbox{in}\quad \BR^N,
	\end{aligned}\right.
\end{equation}
with the parameters in \eqref{eq:resolvent_whole} satisfy $\mu, \nu,\gamma_1,\gamma_2>0.$
Now, we recall the notion in \eqref{def:domain} and the results in \cite[Subsec. 3.1]{ShiE2018}.
\begin{thm}\label{thm:wh_1}
Let $1<p<\infty,$ $0<\ep<\pi \slash 2,$ and $\lambda_0>0.$
Then there exist two families of operators 
$$\big(\CP(\lambda),\CV(\lambda)\big) \in 
\Hol\Big(\Sigma_{\ep} \cap K; \CL\big( H^{1,0}_p(\BR^N) ;
H^{1,2}_p(\BR^N)\big)\Big),$$
such that $(\eta, \bu)= \big(\CP(\lambda),\CV(\lambda)\big)(d,\bff)$ is a solution of \eqref{eq:resolvent_whole}. 
Moreover, there exists a constant $C_{\ep,\lambda_0}$ so that 
\begin{align*}
|\lambda|\|\eta\|_{H^{1}_p(\BR^N)} 
+ \sum_{j=0}^2 |\lambda|^{j\slash 2}\| \bu\|_{H^{2-j}_p(\BR^N)}
\leq C_{\ep,\lambda_0}  \|(d,\bff)\|_{H^{1,0}_p(\BR^N)}
\end{align*}
for any $\lambda \in V_{\ep,\lambda_0}.$
\end{thm}

\begin{thm}\label{thm:wh_2}
Let $1<p<\infty,$ $0<\ep<\pi \slash 2,$ $L\geq R>0,$ and $N\geq 3.$
Set that 
\begin{equation*}
X_{p,L}(\BR^N)= \{(d,\bff)\in H^{1,0}_p(\BR^N) 
: \supp d, \,\supp \bff \subset \overline{B_L}\,\}.
\end{equation*}
Then the following assertions hold:
\begin{enumerate}
\item There exist a constant $\lambda_0>0$ and two families of operators 
\begin{align*}
\CM_{\lambda} &\in \Hol\Big(\dot{U}_{\lambda_0}; 
\CL\big( X_{p,L}(\BR^N) ;
H^{1}_{p,\loc}(\BR^N)\big)\Big),\\
\CV_{\lambda} &\in \Hol\Big(\dot{U}_{\lambda_0};
\CL\big( X_{p,L}(\BR^N) ;
H^{2}_{p,\loc}(\BR^N)^N\big)\Big),
\end{align*}
such that for any $(d,\bff)\in X_{p,L}(\BR^N)$ and $\lambda \in K \cap \dot{U}_{\lambda_0}$
\begin{equation*}
(\CM_{\lambda}, \CV_{\lambda} ) (d,\bff) 
= \big(\CP(\lambda),\CV(\lambda)\big) (d,\bff).
\end{equation*}
Moreover, there exist families of operators
\begin{align*}
\CM_{\lambda}^i &\in \Hol\Big(U_{\lambda_0}; 
\CL\big( X_{p,L}(\BR^N); H^{1}_{p}(B_L)\big)\Big)\,\,\,(i=1,2),\\
\CV_{\lambda}^j &\in \Hol\Big(U_{\lambda_0}; 
\CL\big( X_{p,L}(\BR^N) ;H^{2}_{p}(B_L)^N\big)\Big)\,\,\,(j=0,1,2),
\end{align*}
such that 
\begin{align*}
\CM_{\lambda} &= (\lambda^{N-2} \log \lambda) \CM_{\lambda}^1  + \CM_{\lambda}^2, \\
\CV_{\lambda} &= \big(\lambda^{N\slash 2-1} (\log \lambda)^{\sigma(N)} \big) \CV_{\lambda}^0
+ (\lambda^{N-2} \log \lambda) \CV_{\lambda}^1 +\CV_{\lambda}^2
\end{align*}
for any $\lambda \in \dot{U}_{\lambda_0}$ and 
$\sigma (N)=\big((-1)^N+1 \big)\slash 2.$

\item There exist operators
\begin{equation*}
\CM_0\in \CL\big( X_{p,L}(\BR^N) ;
H^{1}_{p,\loc}(\BR^N)\big),\quad 
\CV_0\in \CL\big( X_{p,L}(\BR^N) ;
H^{2}_{p,\loc}(\BR^N)^N\big)
\end{equation*}
such that for any $(d,\bff) \in X_{p,L}(\BR^N),$ 
$(\eta_0,\bu_0)=(\CM_0,\CV_0)(d,\bff)$ is a solution of 
\begin{equation*}
	\left\{\begin{aligned}
&\gamma_1\di \bu = d  
	    &&\quad\hbox{in}\quad \BR^N, \\
&-\Di\big( \bS(\bu) -\gamma_2 \eta \bI \big)=\bff
		&&\quad\hbox{in}\quad \BR^N,
	\end{aligned}\right.
\end{equation*}
satisfying the estimates:
\begin{align*}
\|\nabla \eta_0\|_{L_p(\BR^N)} 
+ \sum_{|\alpha|=2} \| \pd^{\alpha}_x\bu_0\|_{L_p(\BR^N)}
\leq C_{L,p,N}  \|(d,\bff)\|_{H^{1,0}_p(\BR^N)},
\end{align*}
\begin{align*}
\sup_{|x|\geq 2L}|x|^{N-1} |\eta_0(x)| 
+ \sum_{|\alpha|=0}^1 \sup_{|x|\geq 2L}|x|^{N-2+|\alpha|} |\pd_x^{\alpha}\bu_0(x)|
\leq C_{L,p,N}  \|(d,\bff)\|_{L_p(\BR^N)},
\end{align*}
for some constant $C_{L,p,N}>0.$ 
Furthermore, we have 
\begin{equation*}
\lim_{\substack{ |\lambda|\rightarrow 0 \\ 
|\arg \lambda| \leq  \pi \slash 4}}
\Big(\|\CM_{\lambda}(d,\bff)-\eta_0\|_{H^{1}_p(B_L)}
+\|\CV_{\lambda}(d,\bff)-\bu_0\|_{H^{2}_p(B_L)} \Big)=0.
\end{equation*}
\end{enumerate}
\end{thm}

\subsection{Construction of the parametrix}
\label{subsec:para}
Without loss of generality, we shall prove Theorem \ref{thm:local_energy} for $L=5R.$ 
We first consider the auxiliary problem:
\begin{equation}\label{resolvent_0_2}
\left\{  \begin{aligned}
&\gamma_1\ol\dv \bu= d
&&\quad &\text{in}& \quad \Omega_{5R}, \\
&-\ol\DV \big( \ol\bS(\bu)
-\gamma_2\eta\bI\big)= \bff
&&\quad &\text{in}& \quad \Omega_{5R}, \\
&\big(\ol\bS(\bu)- \gamma_2\eta\bI\big)
\ol\bn_{\Gamma}= 0
&& \quad &\text{on}&\quad  \Gamma,\\
&\big(\bS(\bu)- \gamma_2\eta\bI\big)
\bn_{_{S_{5R}}}= 0
&& \quad &\text{on}&\quad  S_{5R},
\end{aligned} \right.
\end{equation}
Here, $\bn_{_{S_{5R}}}$ denotes the unit outer normal to 
$S_{5R} = \{x \in \BR^N \mid |x| = 5R\}$. 
Let $3R < b_0 < b_1 < b_2 < b_3 < 4R$ and set 
$$D_{b_1, b_2} = \{x \in \BR^N \mid b_1 < |x| < b_2\}, \quad 
D^+_{b_1, b_2} = \{x \in D_{b_1, b_2} \mid x_j > 0 \enskip(j=1, \ldots, N)\}.
$$
Let $\psi \in  C^\infty_0(\BR^N)$ such that 
${\rm supp}\, \psi \subset D_{b_1, b_2},$ and 
$\psi=1$ on some ball $B \subset D^+_{b_1, b_2}.$
Recall the basis $\{\br_j\}_{j=1}^M$ of the rigid motion in \eqref{eq:rid}. 
We introduce a family of vectors $\fQ_\psi=\{\bq_j\}_{j=1}^M,$  another normalization of $\{\br_j\}_{j=1}^M$ in such a way that 
\begin{equation}\label{norm:1}
(\bq_j, \bq_k)_\psi = (\psi \bq_j, \bq_k)_{\BR^N} =\int_{\BR^N} \psi (x) \,\bq_j(x) \cdot\bq_k(x) \,dx = \delta_{jk}.
\end{equation}
Here and in the following, 
given function $a(x),$ $\bar a(y)$ is defined by $a(x) = \bar a(y)$ in view of \eqref{def:LL_w}. 
Since $x=X_{\bw}(y,T)=y$ for $y \not\in B_{2R},$ we infer from \eqref{norm:1} that 
\begin{equation}\label{norm:2}
\delta_{jk}=(\psi\bq_j, \bq_k)_{\BR^N} 
=(\psi\bq_j, \bq_k)_{\Omega} 
= (\bar\psi\bar\bq_j, J\bar\bq_k)_{\Omega} 
=(\psi\wt\bq_j, J\wt\bq_k)_{\Omega_{5R}},
\end{equation}
for $\wt \bq_{j} \in \{\bq_{j}, \bar{\bq}_{j}\}$
with $j,k=1,\ldots, M.$
Moreover, for simplicity we write
\begin{itemize}
\item $\bff \perp \fQ_R$ if $(\bff, J\bar\bq_j)_{\Omega_{5R}}=0$ for any 
 $\bq_j \in \fQ_\psi$; 
\item $\bff\perp \fQ_\psi$ if $(\bff, \bq_j)_\psi = 0$ for any $\bq_j \in \fQ_\psi$.
\end{itemize}
In fact, $\bff \perp \fQ_R$ is equivalent to $\bff \in \CR_d(\Omega_{5R})^{\perp}$ by \eqref{orth:2} in the previous section.
Thus Theorem \ref{thm:st2} yields the following result for \eqref{resolvent_0_2}.
\begin{thm}\label{resolvent_0_bdd} 
Let $1 < p \leq r$. 
Let $(d, \bff) \in H^{1,0}_p(\Omega_{5R})$ with $\bff \perp\fQ_R$.
Then there exist operators
$$(\CJ, \CW) \in \CL(H^{1,0}_p \big(\Omega_{5R}), 
H^{1, 2}_p(\Omega_{5R})\big)$$
such that $(\eta, \bu) = (\CJ, \CW)(d, \bff)$ is a unique solution of 
\eqref{resolvent_0_2} with $\bu \perp \fQ_R$.
Moreover, the following estimate holds,
\begin{equation*}
\|\eta\|_{H^1_p(\Omega_{5R})} + \|\bu\|_{H^2_p(\Omega_{5R})}  
\leq C\big(\|d\|_{H^1_p(\Omega_{5R})} 
+\|\bff\|_{L_p(\Omega_{5R})}  \big),
\end{equation*}
for some constant $C>0.$
\end{thm}
To prove Theorem \ref{thm:ext_1}, we introduce cut-off functions
$\varphi$, $\psi_0$, and $\psi_\infty$ such that 
$0 \leq \varphi, \psi_0, \psi_\infty \leq 1$, $\varphi$, 
$\psi_0$, $\psi_\infty \in C^\infty(\BR^N),$ and 
\begin{equation}\label{eq:cut-off}
\varphi(x) = \begin{cases} 1 &\quad\text{for $|x| \leq b_1$}, \\ 0&\quad\text{for $|x| \geq b_2$},
\end{cases} \quad 
\psi_0(x) = \begin{cases} 1 &\quad\text{for $|x| \leq b_2$}, \\ 0&\quad\text{for $|x| \geq b_3$},
\end{cases} \quad 
\psi_\infty(x) = \begin{cases} 1 &\quad\text{for $|x| \geq  b_1$}, \\ 0&\quad\text{for $|x| \leq b_0$}.
\end{cases}
\end{equation}
For any $(d, \bff) \in H^{1,0}_{p}(\Omega_{5R})$, we have 
\begin{equation}\label{eq:df_ext}
\|\psi_\infty d\|_{H^1_p(\BR^N)} + \|\psi_\infty\bff\|_{L_p(\BR^N)}
\leq C(\|d\|_{H^1_p(\Omega)} + \|\bff\|_{L_p(\Omega)}).
\end{equation}
Then, by Theorem \ref{thm:wh_2} and \eqref{eq:df_ext}, there exists  a $\lambda_0>0$ such that 
$(\eta_\lambda,\bu_\lambda)=(\CM_\lambda,\CV_\lambda)(\psi_\infty d, \psi_\infty\bff)$
solves the following equations:
\begin{equation}\label{eq:2.1}
	\left\{\begin{aligned}
&\lambda \eta_\lambda  + \gamma_1 \di \bu_\lambda = \psi_\infty d  
	    &&\quad\hbox{in}\quad \BR^N, \\
&\gamma_1\lambda \bu_\lambda-\Di\big( \bS(\bu_\lambda) -\gamma_2 \eta_\lambda \bI \big)=\psi_\infty \bff
		&&\quad\hbox{in}\quad \BR^N,
	\end{aligned}\right.
\end{equation}
and satisfies the estimate:
\begin{equation}\label{es:etau_1} 
\|\eta_\lambda\|_{H^1_p(B_{6R})}
+ \|\bu_\lambda\|_{H^2_p(B_{6R})} 
\leq C(\|d\|_{H^1_p(\Omega)}
+ \|\bff\|_{L_p(\Omega)}).
\end{equation}
Moreover, we set $(\eta_0, \bu_0) = (\CM_0, \CV_0)(\psi_\infty d, \psi_\infty\bff)
\in H^{1,2}_{p, {\rm loc}}(\BR^N)$ fulfilling that
\begin{equation}\label{whole:1}
	\left\{\begin{aligned}
&\gamma_1\di \bu_0 =\psi_\infty d  
	    &&\quad\hbox{in}\quad \BR^N, \\
&-\Di\big( \bS(\bu_0) -\gamma_2 \eta_0 \bI \big)=\psi_\infty \bff
		&&\quad\hbox{in}\quad \BR^N,
	\end{aligned}\right.
\end{equation}
and
\begin{equation}\label{etua_1_lim}
\lim_{\substack{\lambda \in \dot U_{\lambda_0}\\
|\lambda|\to 0}}
(\|\eta_\lambda-\eta_0\|_{H^1_p(B_{6R})} + \|\bu_\lambda-\bu_0\|_{H^2_p(B_{6R})}) = 0.
\end{equation}

On the other hand, let us set
$$\bff_{\CR_d}= \sum_{j=1}^M (\psi_0\bff, J\bar\bq_j)_{\Omega_{5R}} \psi \bq_j, 
\quad \bff_\perp = \psi_0\bff - \bff_{\CR_d}
\in L_p(\Omega_{5R})^N.$$
Obviously, $\bff_{\perp} \perp \fQ_R.$ 
In fact, \eqref{norm:2} implies that
$$(\bff_\perp, J\bar\bq_\ell)_{\Omega_{5R}} 
= (\psi_0\bff, J\bar\bq_\ell)_{\Omega_{5R}} 
- \sum_{j=1}^M (\psi_0\bff, J\bar\bq_j)_{\Omega_{5R}}
(\psi \bq_j, J\bar\bq_\ell)_{\Omega_{5R}}=0,
$$
for any $\ell=1,\ldots,M.$
Then, Theorem \ref{resolvent_0_bdd} yields that there exists a (unique) solution
$(\eta_\sharp, \bu_\sharp)\in H^{1,2}_p(\Omega_{5R})$ 
with $\bu_\sharp \perp \fQ_R$ of the following equations:
\begin{equation}\label{resolvent_0_3} 
\left\{  \begin{aligned}
&\gamma_1\ol\dv \bu_{\sharp}= \psi_0  d
&&\quad &\text{in}& \quad \Omega_{5R}, \\
&-\ol\DV \big( \ol\bS(\bu_{\sharp})
-\gamma_2\eta_{\sharp}\bI\big)= \bff_{\perp}
&&\quad &\text{in}& \quad \Omega_{5R}, \\
&\big(\ol\bS(\bu_{\sharp})- \gamma_2\eta_{\sharp}\bI\big)
\ol\bn_{\Gamma}= 0
&& \quad &\text{on}&\quad  \Gamma,\\
&\big(\bS(\bu_{\sharp})- \gamma_2\eta_{\sharp}\bI\big)
\bn_{_{S_{5R}}}= 0
&& \quad &\text{on}&\quad  S_{5R},
\end{aligned} \right.
\end{equation}
possessing the estimate
\begin{equation}\label{es:etau_2}
\|\eta_\sharp\|_{H^1_p(\Omega_{5R})} 
+\|\bu_\sharp\|_{H^2_p(\Omega_{5R})}
\leq C(\|d\|_{H^1_p(\Omega)} 
+\|\bff\|_{L_p(\Omega)}).
\end{equation}

We now introduce parametrices:
$$\wt\eta_\lambda = \Phi_\lambda(d, \bff) = (1-\varphi)\eta_\lambda + \varphi\eta_\sharp,\quad
\wt\bu_\lambda = \Psi_\lambda(d, \bff) = (1-\varphi)\bu_\lambda + 
\varphi\bu_\sharp $$
for $\lambda \in \dot U_{\lambda_0}\cup\{0\}.$
Notice that
$$\ol\bS(\wt\bu_\lambda) - \gamma_2 \wt\eta_\lambda\bI
= (1-\varphi) \big(\bS(\bu_\lambda)- \gamma_2 \eta_\lambda\bI\big)
+ \varphi(\ol\bS(\bu_\sharp) - \gamma_2 \eta_\sharp\bI) + \bV_\lambda(d, \bff),$$
with
\begin{equation}\label{eq:V_df}
\bV_\lambda(d, \bff)= \mu\big( 
(\bu_{\sharp}-\bu_{\lambda})\otimes\nabla \varphi + \nabla \varphi\otimes(\bu_{\sharp}-\bu_{\lambda}) \big) + (\nu-\mu) \big( (\bu_{\sharp}-\bu_{\lambda}) \cdot \nabla \varphi \big) \bI,
\end{equation}
and ${\rm supp}\, \bV_\lambda(d, \bff) \subset D_{b_1, b_2}.$ 
Then it holds that
\begin{equation}\label{eq:etau_3}
\left\{\begin{aligned}
&\lambda \wt\eta_{\lambda} + \gamma_1\overline{\dv}
\wt\bu_{\lambda} = d +\CD_\lambda(d, \bff) 
&&\quad &\text{in}& \quad \Omega, \\
&\gamma_1\lambda \wt\bu_{\lambda}  
- \ol\DV\big(\ol\bS(\wt\bu_{\lambda} )-\gamma_2\wt\eta_{\lambda}  \bI \big)= \bff +\CF_\lambda(d, \bff)
&&\quad &\text{in}& \quad \Omega, \\
&\big(\ol\bS(\wt\bu_{\lambda} )-\gamma_2\wt\eta_{\lambda}  \bI \big)\ol\bn_{\Gamma}=0 
&&\quad &\text{on}& \quad \Gamma,
\end{aligned}\right.
\end{equation}
where $\CD_\lambda(d, \bff)$ and $\CF_\lambda(d, \bff)$ with $\lambda \in \dot U_{\lambda_0}\cup\{0\}$
are defined as
\begin{align*}
\CD_\lambda(d, \bff) & = \lambda \varphi \eta_\sharp
+\gamma_1 \nabla\varphi\cdot(\bu_\sharp-\bu_\lambda),\\
\CF_{\lambda}(d, \bff) & = 
-\varphi \bff_{\CR_d} +
 \gamma_1 \lambda\varphi\bu_\sharp
+\Big( \big(\bS(\bu_\lambda)-\gamma_2 \eta_\lambda\bI \big) -
\big(\bS(\bu_\sharp)-\gamma_2\eta_\sharp \bI \big)\Big)\nabla\varphi 
 -\DV\bV_\lambda(d, \bff).
\end{align*}
Moreover, by \eqref{es:etau_1}, \eqref{etua_1_lim}, and \eqref{es:etau_2}, we have
\begin{equation}\label{eq:etau_2_lim}
\begin{aligned}
\CQ_\lambda(d, \bff)= &(\CD_\lambda, \CF_\lambda)(d, \bff) 
\in X_{p, 5R}(\Omega) \quad
\text{for any $\lambda \in \dot U_{\lambda_0} \cup \{0\}$}, \\
&\lim_{\substack{\lambda \in \dot U_{\lambda_0} \\|\lambda|\to 0}}
\|\CQ_\lambda - \CQ_0\|_{\CL(X_{p, 5R}(\Omega))} = 0.
\end{aligned}
\end{equation}
In particular, by Rellich's compactness theorem, $\CQ_0$ is a compact operator 
on $X_{p, 5R}(\Omega)$.  Moreover, $\CI + \CQ_0$ is invertible due to the following lemma.
\begin{lem}\label{lem:unique_2} 
Given $\CQ_0=(\CF_0, \CG_0)$ as above, $\CI + \CQ_0$ has 
a bounded inverse in $\CL(X_{p, 5R}(\Omega))$, which is denoted by
$(\CI + \CQ_0)^{-1}$.
\end{lem}
The proof of Lemma \ref{lem:unique_2} is postponed to the next subsection, and we continue the proof of Theorem \ref{thm:ext_1}.
Notice that
\begin{equation*}
\CI+ \CQ_{\lambda} =(\CI+\CQ_0) \big(\CI+(\CI+\CQ_0)^{-1}
 (\CQ_{\lambda}-\CQ_0)\big).
\end{equation*}
Then by \eqref{eq:etau_2_lim} and Lemma \ref{lem:unique_2}, we can choose $\lambda_1 \leq \lambda_0$ such that 
\begin{equation*}
\|(\CI+\CQ_0)^{-1} (\CQ_{\lambda}-\CQ_0)\|_{\CL(X_{p,5R}(\Omega))} \leq 1\slash 2,
\end{equation*}
for any $\lambda \in\dot{U}_{\lambda_1}.$ 
Moreover, the inverse of $\CI+ \CQ_{\lambda}$ exists for any 
$\lambda \in\dot{U}_{\lambda_1}$ formulated by
\begin{equation*}
(\CI +\CQ_{\lambda})^{-1} = \sum_{j=0}^\infty 
\big(-(\CI+\CQ_0)^{-1} (\CQ_{\lambda}-\CQ_0)\big)^j(\CI +\CQ_0)^{-1}.
\end{equation*} 
Thus we can define the operators $(\BM_{\lambda},\BV_{\lambda})$ as
\begin{equation*}
(\eta,\Bu)=(\BM_{\lambda},\BV_{\lambda})(d,\bff)
= (\Phi_{\lambda},\BPsi_{\lambda})\circ(\CI +\CQ_{\lambda})^{-1}(d,\bff),
\end{equation*}
for any $(d,\bff) \in X_{p,5R}(\Omega)$ and $\lambda \in\dot{U}_{\lambda_1}.$
In fact, $(\eta,\Bu)$ clearly satisfies \eqref{resolvent_0} by \eqref{eq:etau_3}.
Furthermore, the existence of operators $\BM_{\lambda}^i$ ($i=1,2$) and $\BV_{\lambda}^j$ ($j=0,1,2$) is immediate from Theorems \ref{thm:wh_2} and \ref{resolvent_0_bdd}. So the details are omitted here.

\subsection{Proof of Lemma \ref{lem:unique_2}} 
Before proving Lemma \ref{lem:unique_2}, we start with the following lemma, which plays an essential role in the later proof.
\begin{lem}\label{lem:unique_1} 
Let $1 < p \leq r$ and $N\geq 3.$  
Assume that $\Omega$ is a $C^3$ exterior domain in $\BR^N.$
Let $(\bu, \eta)\in H^{1,2}_{p, {\rm loc}}(\Omega)$ satisfy the homogeneous equations: 
\begin{equation}\label{eq:hom_1}
\left\{\begin{aligned}
&\overline{\dv}\bu = 0
&&\quad &\text{in}& \quad \Omega, \\
&- \ol\DV\big(\mu \ol\bD(\bu)-\gamma_2\eta \bI \big)= 0
&&\quad &\text{in}& \quad \Omega, \\
&\big(\mu\ol\bD(\bu)-\gamma_2\eta \bI \big)\ol\bn_{\Gamma}=0 
&&\quad &\text{on}& \quad \Gamma.
\end{aligned}\right.
\end{equation}
Assume in addition that $\eta(y) = C_\infty + O(|y|^{-(N-1)})$,
$\bu(y) = O(|y|^{-(N-2)})$, and $\nabla \bu(y) = O(|y|^{-(N-1)})$ as 
$|y| \to \infty$, where $C_\infty$ is some constant.  Then, $u=0$ and $\eta=0$
for almost all $y \in \Omega.$  
In particular, $C_\infty=0.$
\end{lem}
\begin{proof}
{\bf Case $ p\geq 2.$}
It suffices to consider $(\eta, \bu) \in H^{1,2}_{2, {\rm loc}}(\Omega)$ in this situation.
Let $\phi$ be a $C^\infty(\BR^N)$ function which equals $1$ for
$y \in B_1$ and $0$ for $y \not\in B_2$.  Let $\phi_L(y) = \phi(y/L)$
for any large $L > 6R$. Set $\tilde C = \gamma_2 C_\infty$ and 
$\tilde\eta = \gamma_2 (\eta- C_\infty)$.  Then, $(\tilde\eta, \bu) \in H^{1,2}_{2, {\rm loc}}(\Omega)$
solves equations:
\begin{equation}\label{eq:3.2}
\left\{\begin{aligned}
&\overline{\dv}\bu = 0
&&\quad &\text{in}& \quad \Omega, \\
&- \ol\DV\big(\mu \ol\bD(\bu)-\tilde \eta\, \bI \big)= 0
&&\quad &\text{in}& \quad \Omega, \\
&\big(\mu\ol\bD(\bu)-\tilde \eta \,\bI \big)\ol\bn_{\Gamma}
=\tilde C\, \ol\bn_{\Gamma}
&&\quad &\text{on}& \quad \Gamma,
\end{aligned}\right.
\end{equation}
Then, by the divergence theorem
of Gau\ss, \eqref{eq:3.2}, and \eqref{support_0},  we have
\begin{align*}
0 & = -\Big(\ol\DV \big(\mu \ol\bD(\bu) - \tilde\eta\bI\big), J\phi_L\bu\Big)_{\Omega}\\
&= -(\tilde C\, \ol\bn_\Gamma, J\phi_L\bu)_\Gamma 
+\frac{\mu}{2}\int_{\Omega} \ol\bD(\bu): \ol\bD(\phi_L\bu) \, J dy
-\big( J\tilde\eta, \overline{\dv}(\phi_L\bu) \big)_{\Omega}\\
& =  -\tilde C\int_\Omega \overline{\dv}(\phi_L\bu)\,J dy 
+\frac{\mu}{2}\int_{\Omega} \ol\bD(\bu): \phi_L\ol\bD(\bu) \, J dy
+ I_L,
\end{align*}
with
$$I_L =\frac{\mu}{2} \int_{L \leq  |y| \leq 2L} 
\bD(\bu):(\nabla \phi_L \otimes \bu + \bu\otimes \nabla \phi_L )\,dy
+ \int_{L \leq |y| \leq 2L}\tilde \eta (\nabla\phi_L \cdot\bu)\,dy. $$
Using the radition condition, we have
$$|I_L| \leq CL^{-1}\int_{L \leq  |y| \leq 2L}|y|^{-(N-1)}|y|^{-(N-2)}\,dy \leq CL^{-(N-2)}\to0
$$
as $L\to\infty$, because $N \geq 3$.  Moreover, \eqref{support_0} implies that 
$$\int_\Omega \overline{\dv}(\phi_L\bu) J \,dy 
= \int_{L \leq |y| \leq  2L}(\nabla\phi_L)\cdot\bu\,dy.
$$
Since $\dv\bu = 0$ for $y \in \BR^N\setminus B_{2R}$, by \cite[Lemma 6]{Shi2018}
there exists a $\bv \in H^1_2(\BR^N)^N$ supported in $D_{R_2, R_3}$ for
some $2R < R_2 < L/2  <  3L < R_3$ such that $\dv \bv=0$ in $\BR^N$ and 
$(\nabla\phi_L)\cdot\bu = (\nabla\phi_L)\cdot\bv$ in $\BR^N$.   Hence, 
$$\int_{L \leq  |y|\leq 2L}(\nabla\phi_L)\cdot\bu\,dy 
= \int_{|y| < R_3}(\nabla\phi_L)\cdot\bv\,dy
= \int_{|y| < R_3} \dv(\phi_L\bv)\,dy = 0.$$
Thus, letting $L\to\infty$ yields that
$$\int_{\Omega} \ol\bD(\bu):\ol\bD(\bu) \, J dy=0,$$
which implies that $\ol \bD(\bu) = 0$ in $\Omega,$ or equivalently,
$\bu \in \Span \{\bar \bq_1, \ldots,\bar\bq_M\}$ (see subsection \ref{subsec:para}).  But,
$\bu(y) = O(|y|^{-(N-2)}) \to 0$ as $|y| \to \infty$, and so $\bu=0$. 
By \eqref{eq:hom_1}, 
$\nabla\eta =0$ in $\Omega$ and $\eta=0$ on $\Gamma.$
This shows that $\eta=0$.
Since $\eta-C_\infty \to 0$ as $|y|\to\infty,$ $C_\infty=0.$  
\smallbreak

{\bf Case $1 < p < 2.$}
For $(\eta,\bu) \in H^{1,2}_{p, {\rm loc}}(\Omega)$
with $1<p<2,$ we use the hypoellipticity of the Stokes operator.
Let $\omega$ be a $C^\infty(\BR^N)$ function which equals $1$ for $x \in \BR^N\setminus B_{4R}$
and $0$ for $x \in B_{2R}.$  
Let $\bB$ be the Bogovskii operator, and set 
$$(\bv,\zeta) = \big(\omega\bu-\bB[(\nabla\omega)\cdot\bu], 
\gamma_2\omega (\eta-C_{\infty}) \big).$$ 
We see that 
$(\zeta,\bv) \in H^{1,2}_{p, {\rm loc}}(\BR^N)$ satisfies the Stokes equations:
$$\bv-\mu\Delta \bv + \nabla\zeta = \bff, \quad \dv \bv = 0
\quad\text{in $\BR^N$}$$
with 
\begin{align*}
\bff= &-2\mu \big( (\nabla\omega) \cdot \nabla\big)\bu 
-\mu (\Delta\omega)\bu 
+\mu \Delta \bB[(\nabla\omega)\cdot\bu]\\
&+\gamma_2 (\eta-C_\infty)\nabla\omega
+\omega\bu-\bB[(\nabla\omega)\cdot\bu] \in H^1_p(\BR^N)^N.
\end{align*}
Notice that ${\rm supp}\,\bff \subset B_{4R}.$ 
By the Sobolev imbedding theorem, 
$\bff \in L_q(\BR^N)^N$, where $q$ is an exponent such that 
$0<N(1/p-1/q) <1.$
By the standard result for the Stokes equations, 
we have $(\zeta,\bv) \in H^{1,2}_{q,\loc}(\BR^N),$ 
which yields that 
$(\eta,\bu) \in H^{1,2}_{q, {\rm loc}}(\BR^N\setminus B_{4R})$. 
 \smallbreak
 
 We next consider the interior problem:
 \begin{equation}
\left\{  \begin{aligned}
&\overline{\dv }\bu = 0
&&\quad &\text{in}& \quad \Omega_{5R}, \\
&\gamma_1 \lambda \bu
-\ol\DV \big( \mu \ol\bD(\bu)
-\gamma_2\eta\bI\big)= \gamma_1 \lambda \bu
&&\quad &\text{in}& \quad \Omega_{5R}, \\
&\big(\mu \ol\bD(\bu)- \gamma_2\eta\bI\big)
\ol\bn_{\Gamma}= 0
&& \quad &\text{on}&\quad  \Gamma,\\
&\big(\mu \bD(\bu)- \gamma_2\eta\bI\big)
\bn_{_{S_{5R}}}= \bh
&& \quad &\text{on}&\quad  S_{5R},
\end{aligned} \right.
\end{equation}
with $\bh = \big(\mu\bD(\bu)-\gamma_2\eta\bI\big)\bn_{_{S_{5R}}} \in W^{1-1/q}_q(S_{5R})$ by the discussion above.  
Since $\bu \in H^2_p(\Omega_{5R})^N \subset L_q(\Omega_{5R})^N,$ choosing $\lambda>0$ so large if necessary, 
by Corollary \ref{cor:ms1} $(\eta, \bu) \in H^{1,2}_q(\Omega_{5R})$,
and so we have $(\eta, \bu) \in H^{1,2}_{q, {\rm loc}}(\Omega)$.  
If $p < q < 2$, then repeated use of this argument finally implies that 
$(\eta, \bu) \in H^{1,2}_{2, {\rm loc}}(\Omega)$.
This completes the proof of Lemma \ref{lem:unique_1}. 
\end{proof}

\begin{proof}[ Proof of Lemma \ref{lem:unique_2} ]

 According to the compactness of $\CQ_0$ and the 
Fredholm alternative theorem, it is sufficient to verify the injectivity of 
$\CI + \CQ_0$.
Let $(d, \bff) \in {\rm Ker}\,(\CI+ \CQ_0)$, that is,
$(d, \bff) + \CQ_0(d, \bff) = (0, 0)$.  Since ${\rm supp}\,\CQ_0(d, \bff)
\subset D_{b_1, b_2}$, ${\rm supp}\,(d, \bff) \subset D_{b_1, b_2}$. 
As $\psi_0=\psi_\infty=1$ in  $D_{b_1, b_2}$, we have 
\begin{equation}\label{eq:df_0}
(d, \bff)=(\psi_\infty d, \psi_\infty\bff) = (\psi_0 d, \psi_0\bff).
\end{equation}
In particular, $\bff_\perp = \psi_0\bff - \bff_{\CR_d} 
= \bff - \bff_{\CR_d}$ with
\begin{equation}\label{residue:1}
\bff_{\CR_d} = \sum_{j=1}^M (\bff, \bq_j)_{\Omega_{5R}}
\psi \bq_j.
\end{equation}
Then, in view of \eqref{eq:etau_3}, 
$(\wt\eta_0, \wt\bu_0) = (1-\varphi)(\eta_0, \bu_0)
 + \varphi(\eta_\sharp, \bu_\sharp)$
satisfies the homogeneous equations: 
\begin{equation*}
\left\{\begin{aligned}
&  \gamma_1\overline{\dv}\,
\wt\bu_{0} =0
&&\quad &\text{in}& \quad \Omega, \\
&- \ol\DV\big(\mu \ol\bD(\wt\bu_{0} )-\gamma_2\wt\eta_{0}  \bI \big)=0
&&\quad &\text{in}& \quad \Omega, \\
&\big(\mu \ol\bD(\wt\bu_{0} )-\gamma_2\wt\eta_{0}  \bI \big)\ol\bn_{\Gamma}=0 
&&\quad &\text{on}& \quad \Gamma,
\end{aligned}\right.
\end{equation*}
and the radiation condition:
\begin{equation*}
\wt\eta_0(x) = O(|x|^{-(N-1)}), \quad
\wt\bu_0(x) = O(|x|^{-(N-2)}), \quad 
\nabla \wt\bu_0(x) = O(|x|^{-(N-1)})
\end{equation*}
as $|x| \to \infty$ by Theorem \ref{thm:wh_2}.
Thus, by Lemma \ref{lem:unique_1}, 
\begin{equation}\label{eq:vanish_0}
\wt\eta_0=0, \quad \wt\bu_0=0.
\end{equation}
 Then the choice of $\varphi$ yields that
\begin{equation}\label{eq:vanish_1} 
\eta_\sharp(x) =0, \quad 
\bu_\sharp(x) = 0 \quad\text{for $|x| \leq b_1$}, \quad 
\eta_0(x) = 0, \quad \bu_0(x)=0 \quad\text{ for $|x| \geq b_2$}.
\end{equation}
We now introduce the extensions
$$(\wt\eta_\sharp, \wt\bu_\sharp)(x)
=\begin{cases} (\eta_\sharp, \bu_\sharp)(x)&
\quad\text{for $x \in \Omega_{5R}\setminus B_{2R}$},
\\ 0 &\quad\text{for $x \in B_{2R}$}.
\end{cases}$$
By \eqref{eq:vanish_1}, \eqref{resolvent_0_3}, \eqref{eq:df_0}, and \eqref{support_0}, 
$(\wt\eta_\sharp, \wt\bu_\sharp)\in H^{1,2}_p(B_{5R})$
satisfy equations:
\begin{equation*}
\left\{  \begin{aligned}
&\gamma_1\dv \wt\bu_{\sharp}= d
&&\quad &\text{in}& \quad B_{5R}, \\
&-\DV \big( \bS(\wt\bu_{\sharp})
-\gamma_2\wt\eta_{\sharp}\bI\big)= \bff - \bff_{\CR_d}
&&\quad &\text{in}& \quad B_{5R}, \\
&\big(\bS(\wt\bu_{\sharp})- \gamma_2\wt\eta_{\sharp}\bI\big)
\bn_{_{S_{5R}}}= 0
&& \quad &\text{on}&\quad  S_{5R}.
\end{aligned} \right.
\end{equation*}

On the other hand, by \eqref{eq:vanish_1}, and \eqref{eq:df_0},
\begin{equation}\label{eq:lem_unique_3}
\left\{  \begin{aligned}
&\gamma_1\dv \bu_{0}= d
&&\quad &\text{in}& \quad B_{5R}, \\
&-\DV \big( \bS(\bu_{0})
-\gamma_2\eta_{0}\bI\big)= \bff 
&&\quad &\text{in}& \quad B_{5R}, \\
&\big(\bS(\bu_{0})- \gamma_2\eta_{0}\bI\big)
\bn_{_{S_{5R}}}= 0
&& \quad &\text{on}&\quad  S_{5R}.
\end{aligned} \right.
\end{equation}
Setting $(\theta, \bv) = (\wt\eta_\sharp-\eta_0, \wt\bu_\sharp-\bu_0)$, we have
\begin{equation}\label{lem_unique_2}
\left\{  \begin{aligned}
&\dv \bv= 0
&&\quad &\text{in}& \quad B_{5R}, \\
&-\DV \big(\mu \bD(\bv)
-\gamma_2\theta \bI\big)= - \bff_{\CR_d}
&&\quad &\text{in}& \quad B_{5R}, \\
&\big(\mu \bD(\bv)- \gamma_2\theta \bI\big)
\bn_{_{S_{5R}}}= 0
&& \quad &\text{on}&\quad  S_{5R}.
\end{aligned} \right.
\end{equation}
We now take the inner product 
$(\cdot, \bq_j)_{B_{5R}}$ on the both side of 
\eqref{lem_unique_2} and use the divergence theorem of Gau\ss,
\begin{align*}
-(\bff_{\CR_d}, \bq_j)_{B_{5R}} &
= \Big(-\DV \big(\mu \bD(\bv) - \gamma_2\theta\bI\big), \bq_j\Big)_{B_{5R}} \\
&= \frac{\mu }{2}\int_{B_{5R}} \bD(\bv): \bD(\bq_j) \,dy 
-(\gamma_2 \theta,
\dv \bq_j)_{B_{5R}} = 0,
\end{align*}
which yields that $(\bff, \bq_j)_{\Omega_{5R}} = (\bff_{\CR_d}, \bq_j)_{B_{5R}} = 0.$
Thus, $\bff_{\CR_d} = 0$ due to \eqref{residue:1}.
Furthermore, $(\bv,\theta)$ solves the homogeneous
equations:
\begin{equation*}
\left\{  \begin{aligned}
&\dv \bv= 0
&&\quad &\text{in}& \quad B_{5R}, \\
&-\DV \big(\mu \bD(\bv)
-\gamma_2\theta \bI\big)= 0
&&\quad &\text{in}& \quad B_{5R}, \\
&\big(\mu\bD(\bv)- \gamma_2\theta \bI\big)
\bn_{_{S_{5R}}}= 0
&& \quad &\text{on}&\quad  S_{5R}.
\end{aligned} \right.
\end{equation*}
Therefore, $\theta=0$ and $\bD(\bv)=0$ in $B_{5R}.$ In particular, 
$\eta_0=\eta_\sharp$ and $\bD(\bu_0-\bu_\sharp)=0$ in $\Omega_{5R}$.  
According to \eqref{eq:vanish_0} and \eqref{eq:vanish_1}, we have  
 $$0=\eta_0 + \varphi(\eta_\sharp-\eta_0)
 = \eta_0
 \,\,\, \hbox{in}\,\,\,\Omega_{5R}.$$

On the other hand, by  $\bD(\bu_0-\bu_\sharp)=0$ on $\Omega_{5R}$, we assume that 
$ \bu_0= \bu_\sharp + \bq_{\star}$ for some $\bq_{\star} \in \Span\{\fQ_\psi\}$.  
Since $\bu_\sharp \perp \fQ_R$ in
$\Omega_{5R}$ and \eqref{eq:vanish_1}, we have 
$(\bu_\sharp, \bq_j)_{\Omega_{5R}} = 0$ 
for all $j=1,\ldots,M,$
 which yields that  $(\bu_\sharp, \bq_{\star})_{\Omega_{5R}}
 =0$. Thus, by \eqref{eq:vanish_0}, we see that
 $$0 = (\wt\bu_0, \bq_{\star})_{\Omega_{5R}} 
 = \big( \bu_\sharp + (1-\varphi)(\bu_0-\bu_\sharp), \bq_{\star} \big)_{\Omega_{5R}}
 = \big( (1-\varphi)\bq_{\star}, \bq_{\star} \big)_{\Omega_{5R}}=0,
 $$
 which implies that $(1-\varphi)\bq_{\star}\cdot \bq_{\star}=0$ in $\Omega_{5R}.$ 
Then $\bq_{\star}=0$ by the choice of $\varphi.$ 
Again using \eqref{eq:vanish_0}, we have 
 $$0=\bu_0 + \varphi(\bu_\sharp-\bu_0)
 = \bu_0 -\varphi\bq_{\star} = \bu_0 
 \,\,\, \hbox{in}\,\,\,\Omega_{5R}.$$
Therefore, $(\eta_0,\bu_0) = (0, 0)$ in $\Omega_{5R},$ 
which, combined with equations \eqref{eq:lem_unique_3}, yields that
 $(d, \bff)=(0,0).$  This completes the proof of Lemma \ref{lem:unique_2}.
 
\end{proof}


\section{Some auxiliary problem}
\label{sec:aux}

In thin section, we introduce some model problem which will be useful for our later study on the resolvent estimates of the Lam\'e operators.
For any fixed $\lambda \in \Sigma_{\ep}$ and $0<\ep<\pi/2,$  we consider the following system in the uniform $C^2$ domain $G:$
\begin{equation}\label{eq:aux_1}
	\left\{\begin{aligned}
&\zeta \bu-\Di\big( \alpha\bD(\bu) 
+ ( \beta-\alpha+\gamma\lambda^{-1} )\di \bu\bI \big)=\bff
		&&\quad\hbox{in}\quad G,\\
&\big( \alpha\bD(\bu) 
+ ( \beta-\alpha+\gamma\lambda^{-1} )\di \bu\bI \big) \Bn_{\pa G} =\Bg
	    &&\quad\hbox{on}\quad \pa G,
	\end{aligned}\right.
\end{equation}
with the constants $\alpha,\beta,\gamma>0$ 
and the parameter $\zeta \in \BBC.$  
For $\lambda \in \Sigma_{\ep}$ and $\zeta_0>0,$  we introduce that 
\begin{align*}
\Xi_{\ep,\lambda}^1&=\{z\in \Sigma_{\ep}: |\arg z -\arg (\alpha +\beta+\gamma \lambda^{-1})| \leq \pi -\ep \},\\
\Xi_{\ep,\lambda}^2&=\{z\in \Xi_{\ep,\lambda}^1: |\arg z -\arg (2\alpha +\beta+\gamma \lambda^{-1})| \leq \pi -\ep \},\\
\Xi_{\ep,\lambda}&= \begin{cases}
\Xi_{\ep,\lambda}^2 & \text{for} \,\,\, \Im \lambda =0, \\
\{z \in \Xi_{\ep,\lambda}^2 : \Im \big( (\beta+\gamma\lambda^{-1} )^{-1}z\big) \Im \lambda>0\}  & \text{for} \,\,\, \Im \lambda \not=0, 
\end{cases}\\
\Xi_{\ep,\lambda,\zeta_0}&=\{z \in \Xi_{\ep,\lambda} : |z|\geq \zeta_0\}.
\end{align*}
Then the result on \eqref{eq:aux_1} reads as follows:
\begin{thm}\label{thm:aux_1}
Let $0<\ep<\pi \slash 2,$ $\alpha,\beta,\gamma,b>0,$ $1<p<\infty,$ and $\lambda\in \Sigma_{\ep,b}.$ Assume that $G$ is a uniform $C^{2}$ domain in $\BR^N.$
Then for any $(\bff,\Bg)\in L_p(G)^N \times H^{1}_p(G)^N,$
there exists $\zeta_0>0$ such that 
\eqref{eq:aux_1} admits a unique solution $\bu \in H^{2}_p(G)^N$
for any $\zeta \in \Xi_{\ep,\lambda,\zeta_0}.$ Moreover, we have 
\begin{equation*}
\sum_{j=0}^{2}|\zeta|^{(2-j) \slash 2} \|\nabla^{j} \bu\|_{L_p(G)}
\leq C\Big( \|\bff\|_{L_p(G)}
+\sum_{j=0}^{1}|\zeta|^{(1-j) \slash 2} \|\nabla^{j} \Bg\|_{L_p(G)}\Big)
\end{equation*}
for some constant $C$ depending solely on 
$\ep,\alpha,\beta,\gamma,b,p,N.$
\end{thm}

\begin{rmk} \label{rmk:aux_1}
Let us give some comments before the proof of Theorem \ref{thm:aux_1}.
\begin{enumerate}
\item 
It is not hard to see that $\Xi_{\ep,\lambda,\zeta_0} \not= \emptyset$ for any 
$0<\ep<\pi/2,$ $\lambda \in \Sigma_{\ep},$ and $\zeta_0 >0.$
In fact, as $\gamma\lambda^{-1}$ still lies in $\Sigma_{\ep}$ 
with $\gamma>0,$ we see that
$$|\arg (\Fs\alpha +\beta+\gamma\lambda^{-1})| \leq \pi-\ep$$
for any $\Fs,\alpha,\beta>0.$ In particular, 
$\BR_+ \subset \Xi_{\ep,\lambda}^2\subset \Xi_{\ep,\lambda}^1.$ 
On the other hand, for any $\omega>0,$ we have 
\begin{equation*}
 \Im \big( (\beta+\gamma\lambda^{-1} )^{-1}z\big) 
 \Im \lambda = \gamma (|\beta \lambda +\gamma|^{-1}\Im \lambda)^2 >0
 \,\,\, ( \Im \lambda \not=0),
\end{equation*}
which implies $\BR_+ \subset \Xi_{\ep,\lambda}$ as well. 
Thus $\{ \omega \in \BR_+: \omega \geq \zeta_0\} \subset 
\Xi_{\ep,\lambda,\zeta_0}.$  

\item In Theorem \ref{thm:aux_1}, $G$ can be bounded or unbounded. Moreover, Theorem \ref{thm:aux_1} can be extended to the uniform domain $G$ of the class $W^{2-1/r}_r$ with $r>N$ and $r \geq \max\{p,p/(p-1)\}$ in the framework of the $\CR-$boundedness theory. 
We refer to \cite{EvBS2014} for more details.
\end{enumerate}
\end{rmk}

As \eqref{eq:aux_1} is an elliptic system, we study only the model problems in $\BR^N$ and $\BR^N_+$ (i.e. \eqref{eq:aux_1_1} and \eqref{eq:aux_1_2} below) with respect to \eqref{eq:aux_1} here, which essentially give us the interior and boundary estimates of \eqref{eq:aux_1}.

\subsection{The auxiliary problem in $\BR^N$}
To prove Theorem \ref{thm:aux_1}, we first consider the model problem in $\BR^N$ ($N\geq 2$),
\begin{equation}\label{eq:aux_1_1}
\zeta \bu- \alpha\Delta \bu - (\beta+\gamma \lambda^{-1} )\nabla \di \bu=\bff
\quad\hbox{in}\quad \BR^N.
\end{equation}
By the discussion in \cite[Sec. 2]{ShiT2004}, the solution of \eqref{eq:aux_1_1} can be formulated by
\begin{multline} \label{eq:sf_aux_1_1}
\bu(x)=\CA(\zeta,\BR^N) \bff = \frac{1}{\alpha} \CF_{\xi}^{-1} 
\Big[ \frac{ \CF_y[\bff] (\xi) }{\alpha^{-1} \zeta + |\xi|^2}\Big](x)\\
-\frac{\beta+\gamma \lambda^{-1}}{\alpha (\alpha +\beta
+\gamma\lambda^{-1})} \CF_{\xi}^{-1}
\Big[ \frac{ \xi \big(\xi \cdot \CF_y[\bff] (\xi) \big) }{(\alpha^{-1} \zeta + |\xi|^2) \big(  (\alpha +\beta+\gamma \lambda^{-1})^{-1} \zeta + |\xi|^2 \big)}\Big](x), \hspace*{1cm}
\end{multline}
with the Fourier transformation and its inverse in $\BR^N$ defined by
\begin{align*}
 \CF_x[f] (\xi) = \int_{\BR^N} e^{-i x\cdot \xi} f(x) \,dx, \quad
\CF_{\xi}^{-1}[f] (x) = (2\pi)^{-N}\int_{\BR^N} e^{i x\cdot \xi} f(\xi) \,d\xi.
\end{align*}
Then by applying the classical Fourier multiplier theory  and \cite[Lemma 2.1]{ShiT2004}, it is not hard to find that
\begin{thm}\label{thm:aux_1_1}
Let $0<\ep<\pi \slash 2,$ $\alpha,\beta,\gamma,b>0,$ $1<p<\infty$ and $\lambda \in \Sigma_{\ep,b}.$ 
For any $\bff \in L_p(\BR^N)^N$ and $\zeta \in \Xi_{\ep,\lambda}^1,$ 
\eqref{eq:aux_1_1} admits a solution $\bu =\CA(\zeta,\BR^N)\bff \in H^{2}_p(\BR^N)^N$ satisfying
\begin{equation*}
\sum_{j=0}^{2}|\zeta|^{(2-j) \slash 2} \|\nabla^{j} \bu\|_{L_p(\BR^N)}
\leq C\|\bff\|_{L_p(\BR^N)}
\end{equation*}
for some constant $C$ depending solely on 
$\ep,\alpha,\beta,\gamma,b,p,N.$
\end{thm}

\subsection{The auxiliary problem in $\BR^N_+$}
In this subsection,  we consider \eqref{eq:aux_1} whenever $G$ is the 
half-space $\BR^N_+,$ namely,
\begin{align*}
\BR^N_+=\{x=(x',x_N)\in \BR^N:x'=(x_1,\dots, x_{N-1})\in \BR^{N-1},\,\,x_{N}>0\},
\end{align*}
and $\pa G=\BR^N_0=\{x=(x',0)\in \BR^N\}.$
Now, \eqref{eq:aux_1} in $\BR^N_+$ is rewritten by
\begin{equation}\label{eq:aux_1_2}
	\left\{\begin{aligned}
&\zeta \bu- \alpha\Delta \bu 
- (\beta+\gamma\lambda^{-1} )\nabla \di \bu=\bff
		&&\quad\hbox{in}\quad \BR^N_+,\\
& \alpha (\pa_N u_j +\pa_j u_N)=-g_{j} \,\,\, (j=1,\dots, N-1)  
&&\quad\hbox{on}\quad \BR^N_0,\\
& 2\alpha \pa_N u_N +  
( \beta-\alpha+\gamma\lambda^{-1} )\di \bu =-g_N
	    &&\quad\hbox{on}\quad \BR^N_0.
	\end{aligned}\right.
\end{equation}

For \eqref{eq:aux_1_2}, we can establish the following result.
\begin{thm}\label{thm:aux_1_2}
Let $0<\ep<\pi \slash 2,$ $\alpha,\beta,\gamma,b>0,$ $1<p<\infty$ and $\lambda \in \Sigma_{\ep,b}.$ 
For any $\bff \in L_p(\BR^N_+)^N,$ $\Bg \in H^{1}_p(\BR^N_+)^N$ and $\zeta \in \Xi_{\ep,\lambda},$ 
\eqref{eq:aux_1_2} admits a solution 
$\bu =\CA(\zeta,\BR^N_+)\bff \in H^{2}_p(\BR^N_+)^N$ satisfying
\begin{equation*}
\sum_{j=0}^{2}|\zeta|^{(2-j) \slash 2} \|\nabla^{j} \bu\|_{L_p(\BR^N_+)}
\leq C\Big( 
\|\bff\|_{L_p(\BR^N_+)}
+\sum_{j=0}^{1}|\zeta|^{(1-j) \slash 2} \|\nabla^{j} \Bg\|_{L_p(\BR^N_+)}\Big)
\end{equation*}
for some constant $C$ depending solely on 
$\ep,\alpha,\beta,\gamma,b,p,N.$
\end{thm}

To prove Theorem \ref{thm:aux_1_2}, we make some reduction.
As in \cite[Sec. 4]{GS2014}, the odd and even extensions for the function $f$ in $\BR^N_+$ are denoted by
\begin{equation*}
f^{o}(x) = 
\begin{cases}
f(x) & \text{for} \,\,\, x_N>0,\\
-f(x',-x_N) & \text{for} \,\,\, x_N<0,
\end{cases}
\quad 
f^{e}(x) = 
\begin{cases}
f(x) & \text{for} \,\,\, x_N>0,\\
f(x',-x_N) & \text{for} \,\,\, x_N<0.
\end{cases}
\end{equation*}
Set that $\bF=(f^o_1,\dots, f^o_{N-1},f^e_N),$ and 
$\bU=\CA(\zeta,\BR^N) \bF.$ 
Then $\bU$ satisfies 
\begin{equation}\label{eq:aux_1_2_1}
	\left\{\begin{aligned}
&\zeta \bU- \alpha\Delta \bU 
- (\beta+\gamma\lambda^{-1} )\nabla \di \bU=\bff
		&&\quad\hbox{in}\quad \BR^N_+,\\
& \pa_N U_N= 0   
&&\quad\hbox{on}\quad \BR^N_0.
	\end{aligned}\right.
\end{equation}
Let $\bv=\bu-\bU \mathds{1}_{\BR^N_+}.$ Then \eqref{eq:aux_1_2} and \eqref{eq:aux_1_2_1} yield that
\begin{equation}\label{eq:aux_1_2_2}
	\left\{\begin{aligned}
&\zeta \bv- \alpha\Delta \bv 
- (\beta+\gamma\lambda^{-1} )\nabla \di \bv=0
		&&\quad\hbox{in}\quad \BR^N_+,\\
& \alpha (\pa_N v_j +\pa_j v_N)=-h_{j} \,\,\, (j=1,\dots, N-1)  
&&\quad\hbox{on}\quad \BR^N_0,\\
& 2\alpha \pa_N v_N 
+(\beta-\alpha+\gamma\lambda^{-1} )\di \bv =-h_N
	    &&\quad\hbox{on}\quad \BR^N_0,
	\end{aligned}\right.
\end{equation}
where the vector $\bh=(h_1,\dots, h_{N-1}, h_N)$ is given by 
\begin{equation*}
h_j= g_j + \alpha (\pa_N U_j +\pa_j U_N),\quad
h_N= g_N +  ( \beta-\alpha+\gamma\lambda^{-1} )\di \bU.
\end{equation*}

Thanks to Theorem \ref{thm:aux_1_1}, 
it is not hard to see that Theorem \ref{thm:aux_1_2} holds true so long as we can verify the following result.
\begin{thm}\label{thm:aux_1_2_1}
Let $0<\ep<\pi \slash 2,$ $\alpha,\beta,\gamma,b>0,$ $1<p<\infty$ and $\lambda \in \Sigma_{\ep,b}.$ 
For any $\bh \in H^{1}_p(\BR^N_+)^N$ and $\zeta \in \Xi_{\ep,\lambda},$ 
\eqref{eq:aux_1_2_2} admits a solution 
$\bv =\CB(\zeta,\BR^N_+)\bh \in H^{2}_p(\BR^N_+)^N$ satisfying
\begin{equation*}
\sum_{j=0}^{2}|\zeta|^{(2-j) \slash 2} \|\nabla^{j} \bv\|_{L_p(\BR^N_+)}
\leq C \sum_{j=0}^{1}|\zeta|^{(1-j) \slash 2} \|\nabla^{j} \bh\|_{L_p(\BR^N_+)}
\end{equation*}
for some constant $C$ depending solely on 
$\ep,\alpha,\beta,\gamma,b,p,N.$
\end{thm}

The proof of Theorem \ref{thm:aux_1_2_1} is standard. As the first stage, we derive the solution formula of \eqref{eq:aux_1_2_2}.
By applying the partial Fourier transform to \eqref{eq:aux_1_2_2}, we have 
\begin{equation}\label{eq:aux_1_2_3}
	\left\{\begin{aligned}
&(\zeta+\alpha |\xi'|^2) \wh{v}_j- \alpha\pa_N^2 \wh{v}_j 
- (\beta+\gamma\lambda^{-1} )
 i\xi_j (i\xi'\cdot \wh{\bv}' +\pa_N \wh{v}_N)=0
		&&\quad\hbox{for}\quad x_N>0,\\
&(\zeta+\alpha |\xi'|^2) \wh{v}_N- \alpha\pa_N^2 \wh{v}_N 
- (\beta+\gamma\lambda^{-1} )
\pa_N(i\xi'\cdot \wh{\bv}' +\pa_N \wh{v}_N)=0
		&&\quad\hbox{for}\quad x_N>0,\\
& \alpha (\pa_N \wh{v}_j +i\xi_j \wh{v}_N)=-\wh{h}_{j} \,\,\, (j=1,\dots, N-1)  
&&\quad\hbox{for}\quad x_N=0,\\
& 2\alpha \pa_N \wh{v}_N 
+  ( \beta-\alpha+\gamma\lambda^{-1} )(i\xi'\cdot \wh{\bv}' +\pa_N \wh{v}_N)=-\wh{h}_N
	    &&\quad\hbox{for}\quad x_N=0,
	\end{aligned}\right.
\end{equation}
where $\wh{f}(\xi',x_N)=\int_{\BR^{N-1}} e^{-ix'\cdot\xi'} f(x',x_N)\,dx',$
and $\xi'=(\xi_1,\dots,\xi_{N-1}).$
\medskip

Next, we introduce
\begin{equation}\label{eq:AB}
A= \sqrt{(\alpha+\beta+\gamma \lambda^{-1})^{-1}\zeta +|\xi'|^2},
\quad B= \sqrt{\alpha^{-1}\zeta +|\xi'|^2},
\end{equation}
and suppose that
\begin{equation*}
\wh{v}_{\ell}=P_{\ell} (e^{-Bx_N} -e^{-Ax_N}) + Q_{\ell} e^{-Bx_N} , 
\,\,\,\forall\,\,\, \ell =1,\dots,N.
\end{equation*}
Then by inserting the formulas above into \eqref{eq:aux_1_2_3} and equating the coefficients of $e^{-Ax_N}$ and $e^{-Bx_N},$ we see that 
\begin{equation}\label{eq:aux_1_2_4}
	\left\{\begin{aligned}
& \alpha (B^2-A^2) P_j 
- (\beta+\gamma\lambda^{-1} ) i\xi_j (i\xi'\cdot P'-AP_N)=0,\\
& \alpha (B^2-A^2) P_N 
+(\beta+\gamma\lambda^{-1} ) A (i\xi'\cdot P'-AP_N)=0,\\
& i \xi'\cdot (P'+Q') -B (P_N+Q_N)=0,\\
& \alpha\big( (B-A) P_j+BQ_j-i\xi_j Q_N\big)=\wh{h}_{j}(\xi',0) \,\,\, (j=1,\dots, N-1),\\
&   ( \beta+\alpha+\gamma\lambda^{-1} )
\big( (B-A) P_N+BQ_N\big)
- ( \beta-\alpha+\gamma\lambda^{-1} )
i\xi'\cdot Q'=\wh{h}_N(\xi',0).
	\end{aligned}\right.
\end{equation}
From $\eqref{eq:aux_1_2_4}_2$ and $\eqref{eq:aux_1_2_4}_3,$ we infer that
\begin{equation}\label{eq:P_1}
i\xi'\cdot P' =\frac{|\xi'|^2}{AB-|\xi'|^2}(i\xi'\cdot Q' -BQ_N),\quad 
P_N =\frac{A}{AB-|\xi'|^2}(i\xi'\cdot Q' -BQ_N).
\end{equation}
Thus it suffices to solve the following linear system:
\begin{equation*}
\bL \begin{bmatrix}
i \xi'\cdot Q' \\
Q_N
\end{bmatrix}
= \begin{bmatrix}
L_{11} & L_{12}\\
L_{21} & L_{22}\\
\end{bmatrix} \begin{bmatrix}
i \xi'\cdot Q' \\
Q_N
\end{bmatrix}
=\begin{bmatrix}
i \xi'\cdot \wh{\bh'} (\xi',0) \\
\wh{h}_N(\xi',0) 
\end{bmatrix},
\end{equation*}
where the entries of the \emph{Lopatinski matrix} $\bL=[L_{ij}]_{2\times 2}$ are given by
\begin{gather*}
L_{11}=\frac{\alpha A (B^2-|\xi'|^2)}{AB-|\xi'|^2}, \quad 
L_{12}=\frac{\alpha |\xi'|^2(2AB-B^2-|\xi'|^2)}{AB-|\xi'|^2},\\
L_{21}= \frac{\alpha (2AB-A^2-|\xi'|^2)}{AB-|\xi'|^2}-
\frac{(\beta +\gamma\lambda^{-1}) (A^2-|\xi'|^2)}{AB-|\xi'|^2},\\
L_{22}= ( \alpha+\beta+\gamma\lambda^{-1} )\frac{ B(A^2-|\xi'|^2)}{AB-|\xi'|^2}.
\end{gather*}
Assuming that $\det \bL \not=0,$ we obtain 
\begin{align*}
i \xi' \cdot Q' &= (\det \bL)^{-1} 
(L_{22} \, i \xi'\cdot \wh{\bh'} -L_{12} \wh{h}_N)(\xi',0),\\
Q_N &= (\det \bL)^{-1} 
(-L_{21} \, i \xi'\cdot \wh{\bh'} +L_{11} \wh{h}_N)(\xi',0),
\end{align*}
Then \eqref{eq:P_1}, $\eqref{eq:aux_1_2_4}_1,$ and $\eqref{eq:aux_1_2_4}_4$ imply that
\begin{align*}
P_j&=\frac{(\beta+\gamma\lambda^{-1}) i\xi_j (|\xi'|^2-A^2)}{\alpha (B^2-A^2) \det \bL (AB-|\xi'|^2)} \CH,\\
P_N&=-\frac{(\beta+\gamma\lambda^{-1}) A (|\xi'|^2-A^2)}{\alpha (B^2-A^2) \det \bL (AB-|\xi'|^2)} \CH,\\
Q_j &=\frac{1}{\alpha B} \wh{h}_j(\xi',0) 
+\frac{i\xi_j}{B\det \bL } 
(-L_{21} \, i \xi'\cdot \wh{\bh'} +L_{11} \wh{h}_N)(\xi',0)\\
& \quad -\frac{(\beta+\gamma\lambda^{-1}) i\xi_j (|\xi'|^2-A^2)}{\alpha B (B+A) \det \bL (AB-|\xi'|^2)} \CH,
\end{align*}
with 
\begin{equation*}
\CH=\det \bL (  i\xi'\cdot Q' -BQ_N )= 
\big( (L_{22}+BL_{21}) i \xi'\cdot \wh{\bh'} -(L_{12}+B L_{11}) \wh{h}_N \big)(\xi',0).
\end{equation*}
For simplicity, we denote
\begin{equation*}
N(\zeta,\xi';\lambda)=\frac{(\beta+\gamma \lambda^{-1})
(|\xi'|^2-A^2)}{\alpha (AB-|\xi'|^2)} \cdot
\end{equation*}
Then we have for $j=1, \dots, N-1,$
\begin{align*}
\wh{v}_j (\xi',x_N) =& \frac{1}{\alpha B}e^{-B x_N} \wh{h}_j (\xi',0) \\
&- \Big( \frac{i\xi_j}{B} \frac{L_{21}}{\det \bL} e^{-Bx_N}
 -\frac{N(\zeta,\xi';\lambda) }{ B (B+A)}
\frac{i\xi_j}{\det \bL} (L_{22}+BL_{21}) 
\big(BM(x_N)-e^{-Bx_N} \big)
\Big) (i\xi'\cdot \wh{\bh'}) (\xi',0) \\
& + \Big(  \frac{i\xi_j}{B} \frac{L_{11}}{\det \bL} e^{-Bx_N}
-\frac{N(\zeta,\xi';\lambda) }{ B (B+A)}
\frac{i\xi_j}{\det \bL} (L_{12}+BL_{11}) 
\big(BM(x_N)-e^{-Bx_N} \big)\Big) \wh{h}_N (\xi',0),\\
\wh{v}_N (\xi',x_N) =& -\Big(\frac{L_{21}}{\det \bL} e^{-Bx_N}
 +\frac{N(\zeta,\xi';\lambda) }{ B+A}
\frac{A}{\det \bL} (L_{22}+BL_{21}) M(x_N)
\Big) (i\xi'\cdot \wh{\bh'}) (\xi',0) \\
& +\Big( \frac{L_{11}}{\det \bL} e^{-Bx_N}
+\frac{N(\zeta,\xi';\lambda) }{ B+A}
\frac{A}{\det \bL} (L_{12}+BL_{11}) M(x_N)\Big) \wh{h}_N (\xi',0).
\end{align*}
Thus, we define $\CB(\zeta,\BR^N_+) \bh = \CF^{-1}_{\xi'} (\wh{v}_1,\dots,\wh{v}_{N-1},\wh{v}_N).$
Based on the formula of $\CB(\zeta,\BR^N_+)$ and the Fourier analysis (as in \cite[Sec. 4]{GS2014}), Theorem \ref{thm:aux_1_2_1} can be proved with the following two lemmas.
\begin{lem}\label{lem:basic_half}
Let $0<\ep<\pi/2$ and $\xi' \in \BR^{N-1}.$
Then for any $a>0$ and $\lambda \in \Sigma_{\ep},$ we have 
\begin{equation*}
\big| a \lambda +|\xi'|^2  \big| \geq \sin (\ep \slash 2) 
\big(a |\lambda|+|\xi'|^2 \big).
\end{equation*}
\end{lem}
\begin{lem}\label{lem:detL}
There exists some constant $c$ depending solely on $\ep, b,\alpha,\beta,\gamma$ such that 
\begin{equation}\label{es:detL_0}
|\det \bL | \geq c (|\zeta|+|\xi'|^2)
\end{equation}
for all $\zeta\in \Xi_{\ep,\lambda}.$
\end{lem}
Lemma \ref{lem:basic_half} is standard,
while Lemma \ref{lem:detL} is the key point in our results.
The study of $\det \bL$ is the main task of the next subsection.

\subsection{Proof of Lemma \ref{lem:detL}} \label{sec:A}
For $A$ and $B$ defined in \eqref{eq:AB}, we introduce
\begin{equation*}
P=P(\zeta,\xi';\lambda)= \frac{\zeta}{AB-|\xi'|^2}
 = \frac{\alpha(\alpha+\beta+\gamma \lambda^{-1}) (AB+|\xi'|^2)}{\zeta+(2\alpha+\beta+\gamma \lambda^{-1})|\xi'|^2}.
\end{equation*}
Then the entries of $\bL$ can be rewritten as follows
\begin{gather*}
L_{11}= A P, \quad L_{12}=|\xi'|^2(2\alpha-P), \quad 
L_{21}=2\alpha -P,\quad L_{22}=BP.
\end{gather*}
Moreover, $\det \bL = (\zeta +4\alpha |\xi'|^2) P -4\alpha^2|\xi'|^2.$
Based on such observation, we can prove Lemma \ref{lem:detL} in following steps.
\bigskip

{\bf Step 1.} Firstly, we consider the case where $|\zeta| |\xi'|^{-2} =O(r_1)\ll 1.$
Then it is not hard to see that 
$$A,B=|\xi'|\big(1+O(r_1)\big), \quad
P=\frac{2\alpha(\alpha+\beta+\gamma \lambda^{-1})}{2\alpha+\beta+\gamma \lambda^{-1}} +O(r_1).$$
Thus we have 
\begin{equation*}
\det \bL =\frac{4\alpha^2 (\beta+\gamma \lambda^{-1})}{2\alpha+\beta+\gamma \lambda^{-1}}  |\xi'|^2 \big(1+O(r_1)\big).
\end{equation*}
Now we choose $|r_1|$ small enough such that 
\begin{align*}
|\det \bL| \geq \frac{2\alpha^2 |\beta+\gamma \lambda^{-1}|}{|2\alpha+\beta+\gamma \lambda^{-1}|}  |\xi'|^2 
\geq \frac{\sin (\ep\slash 2)\alpha^2 \beta}{(2\alpha+\beta+\gamma b^{-1})} (|\zeta| + |\xi'|^2).
\end{align*}

{\bf Step 2.} Secondly, we consider the case where 
$|\xi'|^2 |\zeta|^{-1}  =O(r_2)\ll 1.$
It is not hard to see that 
\begin{gather*}
A=(\alpha+\beta+\gamma \lambda^{-1})^{-1/2}
\zeta^{1/2}\big(1+O(r_2)\big), 
\quad B=\alpha^{-1/2}\zeta^{1/2}\big(1+O(r_2)\big), \\
P=\alpha^{1/2} (\alpha+\beta+\gamma \lambda^{-1})^{1/2} \big(1+O(r_2)\big).
\end{gather*}
Thus, we have 
\begin{equation*}
\det \bL = \alpha^{1/2} (\alpha+\beta+\gamma \lambda^{-1})^{1/2} \zeta \big(1+O(r_2)\big).
\end{equation*}
Now, we can choose $|r_2|$ small enough such that 
\begin{align*}
|\det \bL| \geq 2^{-1}\alpha^{1/2} |\alpha+\beta+\gamma \lambda^{-1}|^{1/2} |\zeta| 
\geq 4^{-1} \big(\sin (\ep/2)\alpha(\alpha+\beta)\big)^{1/2} (|\zeta| + |\xi'|^2).
\end{align*}

{\bf Step 3.} By previous two steps, we consider
$R_2^{-1} |\zeta|^{1/2} \leq |\xi'| \leq R_1 |\zeta|^{1/2} $
for some fixed $R_1,R_2>0.$
Set that 
\begin{gather*}
\wt{\zeta}= \zeta (|\zeta|^{1/2}+|\xi'|)^{-2}, \quad 
\wt{\xi'}=\xi' (|\zeta|^{1/2}+|\xi'|)^{-1}\\
\wt A= \sqrt{(\alpha+\beta+\gamma \lambda^{-1})^{-1}
\wt \zeta +|\wt \xi'|^2},\quad 
\wt B= \sqrt{\alpha^{-1}\wt \zeta +|\wt \xi'|^2}.
\end{gather*}
For $R_2^{-1} |\zeta|^{1/2} \leq |\xi'| \leq R_1 |\zeta|^{1/2},$ 
$(\wt\zeta,\wt\xi') \in D(R_1,R_2)$ with
\begin{align*}
D(R_1,R_2)=\big\{(z,\eta') \in \Sigma_{\ep} \times \BR^{N-1}:& (1+R_1)^{-2}\leq |z| \leq R_2^2 (1+R_2)^{-2},\\
&(1+R_2)^{-1}\leq |\eta'| \leq R_1 (1+R_1)^{-1}\big\}.
\end{align*}
Note that 
\begin{equation*}
P=\frac{\zeta}{AB-|\xi'|^2} = \frac{\wt \zeta}{\wt A \wt B-|\wt \xi'|^2} ,
\end{equation*}
then we have 
\begin{equation*}
\det \bL = (|\lambda|^{1/2} +|\xi'|)^2 \det \wt \bL
\,\,\,\text{with}\,\,\,
\det \wt \bL=  (\wt \zeta +4\alpha |\wt \xi'|^2) P -4\alpha^2|\wt \xi'|^2.
\end{equation*}
To see \eqref{es:detL_0}, we shall show that $|\det \wt \bL|$ is bounded from below.
Suppose that $\det \wt L=0,$ or equivalently, $\det \bL=0$ for some $(\zeta,\xi').$
Then there exists non-trivial $\Bw$ in the form of 
$$w_{\ell}(x_N)=P_{\ell} (e^{-Bx_N} -e^{-Ax_N}) + Q_{\ell} e^{-Bx_N} , 
\,\,\,\forall\,\,\, \ell =1,\dots,N,$$
satisfying the homogeneous equations
\begin{equation}\label{eq:aux_1_2_5}
	\left\{\begin{aligned}
&(\zeta+\alpha |\xi'|^2) w_j- \alpha\pa_N^2 w_j 
- (\beta+\gamma\lambda^{-1} )
 i\xi_j (i\xi'\cdot \Bw' +\pa_N w_N)=0
		&&\quad\hbox{for}\quad x_N>0,\\
&(\zeta+\alpha |\xi'|^2) w_N- \alpha\pa_N^2 w_N 
- (\beta+\gamma\lambda^{-1} )
\pa_N(i\xi'\cdot \Bw' +\pa_N w_N)=0
		&&\quad\hbox{for}\quad x_N>0,\\
& \alpha (\pa_N w_j +i\xi_j w_N)=0 \,\,\, (j=1,\dots, N-1)  
&&\quad\hbox{for}\quad x_N=0,\\
& 2\alpha \pa_N w_N +  ( \beta-\alpha+\gamma\lambda^{-1} )
(i\xi'\cdot \Bw' +\pa_N w_N)=0
	    &&\quad\hbox{for}\quad x_N=0.
	\end{aligned}\right.
\end{equation}

Denote the inner product $(f,g)=\int_0^{\infty} f(x_N) \overline{g(x_N)} \,dx_N,$ and $\|f\|=(f,f)^{1/2}.$
By taking the inner products between $(\eqref{eq:aux_1_2_5}_1,\eqref{eq:aux_1_2_5}_2)$ and $\Bw,$ 
we have
\begin{align*}
0&=\zeta \|\Bw\|^2 +\alpha \sum_{j,k=1}^{N-1}\|i\xi_k w_j\|^2  
+ \alpha \sum_{j=1}^{N-1} \|i\xi_j w_N\|^2
-\alpha \sum_{j=1}^{N-1} (\pa_N^2 w_j, w_j) - \alpha(\pa_N^2 w_N, w_N) \\
&\quad +\alpha \|i \xi'\cdot \Bw'\|^2 +\alpha (\pa_N w_N, i\xi'\cdot \Bw' )
-(\beta+\gamma\lambda^{-1}) 
\big(\pa_N(i\xi'\cdot \Bw' +\pa_N w_N), w_N \big)\\
& \quad + ( \beta-\alpha+\gamma\lambda^{-1} ) (i\xi'\cdot \Bw' +\pa_N w_N, i\xi'\cdot \Bw')\\
&=\zeta \|\Bw\|^2 +\alpha \sum_{j,k=1}^{N-1}\|i\xi_k w_j\|^2  + 2\alpha \|\pa_N w_N\|^2 + \alpha \|i\xi'\cdot \Bw'\|^2 
+ \alpha \sum_{j=1}^{N-1} \|\pa_N w_j +i\xi_j w_N\|^2\\
&\quad + (\beta-\alpha+\gamma\lambda^{-1})
\|i\xi'\cdot \Bw' +\pa_N w_N\|^2,
\end{align*}
where we have used the integration by parts and the boundary conditions in \eqref{eq:aux_1_2_5} for the last equality.
Then taking the imaginary part and real part respectively, we have
\begin{align}\label{eq:w_Im_1}
0&=\Im \zeta \|\Bw\|^2+\gamma\Im (\lambda^{-1})
\|i\xi'\cdot \Bw' +\pa_N w_N\|^2,\\ \label{eq:w_Re_1}
0&=\Re \zeta \|\Bw\|^2 +\alpha \sum_{j,k=1}^{N-1}\|i\xi_k w_j\|^2  + 2\alpha \|\pa_N w_N\|^2 + \alpha \|i\xi'\cdot \Bw'\|^2 
+ \alpha \sum_{j=1}^{N-1} \|\pa_N w_j +i\xi_j w_N\|^2\\ \nonumber
&\quad + \big(\beta-\alpha+\gamma \Re(\lambda^{-1}) \big)
\|i\xi'\cdot \Bw' +\pa_N w_N\|^2.
\end{align}
From \eqref{eq:w_Re_1}, we obtain that
\begin{align}\label{eq:w_Re_2}
0&\geq \Re \zeta \|\Bw\|^2 
+ \alpha \sum_{j=1}^{N-1} \|\pa_N w_j +i\xi_j w_N\|^2
 + \big(\beta+\gamma \Re(\lambda^{-1}) \big)
\|i\xi'\cdot \Bw' +\pa_N w_N\|^2\\ \nonumber
&\geq \Re \zeta \|\Bw\|^2 
 + \big(\beta+\gamma \Re(\lambda^{-1}) \big)
\|i\xi'\cdot \Bw' +\pa_N w_N\|^2.
\end{align}
In what follows, we shall prove that $\|\Bw\|^2 =0,$ which is a contradiction with our choice of $\Bw.$ In other words, $\det \wt \bL\not=0$ for any $(\wt\zeta,\wt\xi') \in D(R_1,R_2).$ By the continuity of $\det \wt\bL,$ we have $\det \wt \bL \geq c$ for some constant $c>0.$
\medskip

If $\Im (\lambda^{-1}) =0,$ or equivalently $\Im \lambda =0,$ then $\Re (\lambda^{-1}) = |\lambda|^{-2} \Re \lambda>0$ as $\lambda^{-1} \in \Sigma_{\ep}.$ Thus \eqref{eq:w_Im_1} and \eqref{eq:w_Re_2} imply that
\begin{equation*}
\Im \zeta \|\Bw\|^2 =0, \quad \Re \zeta \|\Bw\|^2 \leq 0.
\end{equation*}
Thus  $\|\Bw\|=0$ for $\zeta \in \Sigma_{\ep}.$
\smallbreak

If $\Im (\lambda^{-1}) \not =0,$ then 
$\|i\xi'\cdot \Bw' +\pa_N w_N\|^2 = - \Im\zeta  \|\Bw\|^2 / \big(\gamma\Im (\lambda^{-1}) \big)$ from \eqref{eq:w_Im_1}.
Inserting it into \eqref{eq:w_Re_2}, we have 
\begin{equation}\label{eq:w_1}
\Big( \Re \zeta - \frac{\Im \zeta \big(\beta+\gamma \Re(\lambda^{-1}) \big) }{\gamma\Im (\lambda^{-1}) } \Big) \|\Bw\|^2  \leq 0.
\end{equation}
As $\Im \big( (\beta+\gamma\lambda^{-1} )^{-1}\zeta\big) \Im \lambda>0$ for $\zeta \in \Xi_{\ep,\lambda},$ the coefficient in \eqref{eq:w_1} satisfies
\begin{equation*}
 \Re \zeta-\frac{\Im \zeta \big(\beta+\gamma \Re(\lambda^{-1}) \big) }{\gamma\Im (\lambda^{-1})} = \frac{
| \beta \lambda+\gamma|^2
}{\gamma }
 \frac{\Im \big( (\beta+\gamma\lambda^{-1} )^{-1}\zeta\big)}{\Im \lambda} >0.
\end{equation*}
Hence $\|\Bw\|=0.$


\section{Resolvent problem for $\lambda$ away from zero} 
\label{sec:large}
According to Theorem \ref{thm:ext_1}, we shall study the \eqref{resolvent_0} whenever $\lambda$ is uniformly bounded from below in this section. In next subsection, we first consider the case where $\lambda$ is far away from the origin. Then we study \eqref{resolvent_0} whenever $\lambda $ lies in some ring-shaped region.

\subsection{Resolvent problem for large $\lambda$}
Recall the notion in \eqref{def:domain}. The main result of this subsection reads:
\begin{thm}\label{thm:large_res}
Let $1<p\leq r <\infty,$ and $0<\ep<\pi \slash 2.$
Assume that $\Omega$ is a $C^{2}$ exterior domain in $\BR^N$ for $N\geq 3.$
Then there exist $\lambda_2>0$ and two families of operators 
$$\big(\CP_{\infty}(\lambda),\CV_{\infty}(\lambda)\big) \in 
\Hol\Big(V_{\ep,\lambda_2} ; \CL\big( H^{1,0}_p(\Omega) ;
H^{1,2}_p(\Omega)\big)\Big),$$
such that $(\eta, \bu)= \big(\CP_{\infty}(\lambda),\CV_{\infty}(\lambda)\big)(d,\bff) \in H^{1,2}_p(\Omega)$ is a unique solution of \eqref{resolvent_0}
for any $\lambda \in V_{\ep,\lambda_2}$ 
and any $(d,\bff)\in H^{1,0}_p(\Omega).$
Moreover, we have
\begin{equation}\label{eq:large_res_1}
\|\eta\|_{H^{1}_p(\Omega)} +\|\bu\|_{H^{2}_p(\Omega)} \leq C
\big( \|d\|_{H^{1}_p(\Omega)}+ \|\bff\|_{L_p(\Omega)} \big)
\end{equation}
for some constant $C$ depending solely on 
$\lambda_2,\ep,p,\mu,\nu,\gamma_1,\gamma_2,N.$
\end{thm}
 
The proof of Theorem \ref{thm:large_res} is similar to Theorem \ref{thm:ms1}.
In fact, we can rewrite \eqref{resolvent_0} by 
\begin{equation*}
\left\{\begin{aligned}
&\lambda \eta + \gamma_1\dv\bu = d + \rho(\bu)
&&\quad &\text{in}& \quad \Omega, \\
&\gamma_1\lambda \bu 
- \DV\big(\bS(\bu)-\gamma_2\eta \bI \big)= \bff +\bF(\eta,\bu)
&&\quad &\text{in}& \quad \Omega, \\
&\big(\bS(\bu)-\gamma_2\eta \bI \big)\bn_{\Gamma}
=\bH(\eta,\bu)
&&\quad &\text{on}& \quad \Gamma,
\end{aligned}\right.
\end{equation*}
with 
\begin{align*}
\rho (\Bu)=&-\gamma_1 \bV:\nabla \bu,\\
\bV_1(\bu)=&\mu \big(\bV  \nabla \bu+(\bV  \nabla \bu)^{\top}\big)
+(\nu-\mu) (\bV: \nabla \bu)\bI,\\
\bF(\eta,\bu)=& 
\big(\bV\nabla\mid \ol\bS(\bu) -\gamma_2 \eta\bI \big)
+ \DV\bV_1(\bu),\\
\bH(\eta,\bu)=&-\big(\ol\bS(\bu) -\gamma_2 \eta \bI \big) \bV\bn_{\Gamma} - \bV_1(\bu)\bn_{\Gamma}.
\end{align*}
Moreover, \eqref{assump:1} and \eqref{eq:sobolev} yield that 
\begin{align}\label{eq:large_res_3}
\big\|\big(\rho(\bu), \bH(\eta,\bu) \big)\big\|_{H^1_p(\Omega)}
+\|\bF(\eta,\bu)\|_{L_p(\Omega)} \leq C \sigma \|(\eta,\nabla \Bu)\|_{H^1_p(\Omega)}
\end{align}
for any $1<p\leq r<\infty.$
Then Theorem \ref{thm:large_res} can be proved by taking advantage of \eqref{eq:large_res_3}, the fixed point theorem, and the following result in \cite{EvBS2014}.
So the details of the proof are left to the reader.
\begin{thm}{\cite[Theorem 2.4]{EvBS2014}}
\label{thm:large_res_2}
Assume that $\Omega$ is a $C^{2}$ exterior domain in $\BR^N$ for 
$N\geq 3.$
Let $0<\ep<\pi \slash 2,$ and $1<p<\infty.$  Set
\begin{equation*}
Y_p(\Omega)= H^1_p(\Omega) \times L_p(\Omega)^N \times H^1_p(\Omega)^N.
\end{equation*}
Then there exist $\lambda_2>0$ and two families of operators 
$$\big(\wt\CP_{\infty}(\lambda),\wt\CV_{\infty}(\lambda)\big) \in 
\Hol\Big(V_{\ep,\lambda_2} ; \CL\big( Y_p(\Omega) ;
H^{1,2}_p(\Omega)\big)\Big),$$
such that the following system:
\begin{equation}\label{eq:large_res_2}
\left\{\begin{aligned}
&\lambda \eta + \gamma_1\dv\bu = d 
&&\quad &\text{in}& \quad \Omega, \\
&\gamma_1\lambda \bu 
- \DV\big(\bS(\bu)-\gamma_2\eta \bI \big)= \bff
&&\quad &\text{in}& \quad \Omega, \\
&\big(\bS(\bu)-\gamma_2\eta \bI \big)\bn_{\Gamma}
=\bh
&&\quad &\text{on}& \quad \Gamma,
\end{aligned}\right.
\end{equation}
admits a unique solution
$(\eta, \bu)= \big(\wt\CP_{\infty}(\lambda),\wt\CV_{\infty}(\lambda)\big)(d,\bff,\bh) \in H^{1,2}_p(\Omega)$
for $\lambda \in V_{\ep,\lambda_2}$ 
and for $(d,\bff,\bh)\in Y_p(\Omega).$
Moreover, we have
\begin{equation*}
\|\eta\|_{H^{1}_p(\Omega)} +\|\bu\|_{H^{2}_p(\Omega)} \leq C
\big( \|(d,\bh)\|_{H^{1}_p(\Omega)}+ \|\bff\|_{L_p(\Omega)} \big)
\end{equation*}
for some constant $C$ depending solely on 
$\lambda_2,\ep,p,\mu,\nu,\gamma_1,\gamma_2,N.$
\end{thm}

The existence of the semigroup $\{T(t)\}_{t\geq 0}$ associated to \eqref{LL:Lame_1}
is immediate from Theorem \ref{thm:large_res}.
For $1<p,q<\infty,$ we define 
\begin{align*}
\CD_p(\CA_{\Omega})&=\{(\eta,\bu) \in H^{1,0}_p(\Omega) \mid \bu \in H^2_p(\Omega)^N,\,\,\big(\ol\bS(\bu) - \gamma_2 \eta\bI\big)\ol\bn_{\Gamma} = 0 \},\\
\CD_{p,q}(\Omega) &=\big( H^{1,0}_p(\Omega), \CD_p(\CA_{\Omega}) \big)_{1-1/q,q} \subset H^1_p(\Omega) \times 
B^{2(1-1/q)}_{p,q}(\Omega)^N.
\end{align*}
\begin{thm}\label{thm:sg}
The operator $\CA_{\Omega}$ generates a $C_0$-semigroup $\{T(t)\}_{t\geq 0}$ in 
$H^{1,0}_p(\Omega)$ for any $1<p\leq r<\infty,$ which is analytic as well.  Denote the solution of \eqref{LL:Lame_1} by $(\rho,\bv)(t)=T(t)(\rho_0,\bv_0).$
Then there exists positive constants $\gamma_0$ and $C$ such that the following assertions hold.
\begin{enumerate}
\item For $(\rho_0,\bv_0) 
\in H^{1,0}_{p}(\Omega),$ we have
\begin{equation*}
\|(\rho,\bv)(t)\|_{H^{1,0}_{p}(\Omega)} 
+t \big(\|\pd_t (\rho,\bv)(t)\|_{H^{1,0}_{p}(\Omega)} 
+ \| (\rho,\bv)(t)\|_{\CD_p(\CA_{\Omega})}\big) 
\leq Ce^{\gamma_0 t} \|(\rho_0,\bv_0) \|_{H^{1,0}_{p}(\Omega)}.
\end{equation*}

\item For $(\rho_0,\bv_0) 
\in \CD_{p}(\CA_{\Omega}),$ we have
\begin{equation*}
\|\pd_t (\rho,\bv)(t)\|_{H^{1,0}_{p}(\Omega)} 
+ \| (\rho,\bv)(t)\|_{\CD_p(\CA_{\Omega})}
\leq Ce^{\gamma_0 t} \|(\rho_0,\bv_0) \|_{\CD_p(\CA_{\Omega})}.
\end{equation*}

\item For $(\rho_0,\bv_0) 
\in \CD_{p,q}(\Omega),$ we have
\begin{multline*}
\|e^{-\gamma_0 t} (\pd_t \rho,\rho) \|_{L_q(\BR_+; H^1_p(\Omega))} 
+\|e^{-\gamma_0 t} \pd_t \bv\|_{L_q(\BR_+; L_p(\Omega))} 
+\|e^{-\gamma_0 t} \bv \|_{L_q(\BR_+; H^2_p(\Omega))} \\
\leq C \big( \|\rho_0\|_{H^1_p(\Omega)} + \|\bv_0\|_{B^{2(1-1/q)}_{p,q}(\Omega)} \big).
\end{multline*}
\end{enumerate}
\end{thm}

\subsection{Resolvent problem for $\lambda$ in some compact subset} 
\label{sec:compact}
Thanks to Theorem \ref{thm:large_res}  and Theorem \ref{thm:ext_1}, it remains to study \eqref{resolvent_0} whenever $\lambda$ is uniformly bounded from above and also from below. To this end, let us take some suitable positive constants $\lambda_1'$ and $\lambda_2'$ such that 
\begin{equation*}
0<\lambda_1-\lambda_1' \ll 1, \quad
0<\lambda_2'-\lambda_2 \ll 1, \quad  
\end{equation*}
with $\lambda_1$ and $\lambda_2$ given 
by Theorem \ref{thm:ext_1} and  Theorem \ref{thm:large_res}   respectively. 
For fixed constants $\mu,\nu,\gamma_1,\gamma_2>0,$ we set 
\begin{equation}\label{eq:KD_ep}
\begin{aligned}
K_{\ep}&= \Big\{\lambda \in \BBC\backslash \{0\}: 
\big(\Re \lambda +\frac{\gamma_1\gamma_2}{\mu+\nu} +\ep \big)^2 
+\Im \lambda^2 \geq \big(\frac{\gamma_1\gamma_2}{\mu+\nu} +\ep \big)^2\Big\},\\
D_{\ep}'&=\{\lambda \in \Sigma_{\ep} \cap K_{\ep} : \lambda_1'\leq |\lambda|\leq \lambda_2'\}.
\end{aligned}
\end{equation}
In this section, we address the resolvent problem \eqref{resolvent_0} whenever 
$\lambda$ lies in $D_{\ep}'$ above.
The main theorem of this section reads as follows:
\begin{thm}\label{thm:mid_res}
Suppose that $\Omega$ is a $C^{2}$ exterior domain in $\BR^N$ for $N\geq 3.$ 
Let $0<\ep<\pi \slash 2,$ $N<r<\infty,$ $1<p\leq r,$ and $\lambda\in D_{\ep}'.$ 
Then there exist two families of operators 
$$\big(\CP_{mid}(\lambda),\CV_{mid}(\lambda)\big) \in 
\Hol\Big(D'_{\ep} ; \CL\big( H^{1,0}_p(\Omega) ;
H^{1,2}_p(\Omega)\big)\Big),$$
such that $(\eta, \Bu)= \big(\CP_{mid}(\lambda),\CV_{mid}(\lambda)\big)(d,\bff) \in H^{1,2}_p(\Omega)$ is a unique solution of \eqref{resolvent_0}
for any $\lambda \in D_{\ep}'$ and for any 
$(d,\bff)\in H^{1,0}_p(\Omega).$
Moreover, we have
\begin{equation*}
\|\eta\|_{H^{1}_p(\Omega)} +\|\Bu\|_{H^{2}_p(\Omega)} \leq C
\big( \|d\|_{H^{1}_p(\Omega)}+ \|\bff\|_{L_p(\Omega)} \big)
\end{equation*}
for some constant $C$ depending solely on 
$\lambda_1',\lambda_2',\ep,p,r,\mu,\nu,\gamma_1,\gamma_2,N.$
\end{thm}

To check the solvability of \eqref{resolvent_0} for $\lambda$ in $D_{\ep}',$ we will study some model problem in the bounded domain $\Omega_{5R}$ in Subsection \ref{subsec:model}.
Afterwards, we construct the solution operators by fixing $\lambda,$
and then extend the result to the whole region $D_{\ep}'$ by the compactness of $D_{\ep}'.$

\subsubsection{Some model problem in $\Omega_{5R}$}
\label{subsec:model}
Recall the notion in \eqref{def:domain}. We consider the following system in $\Omega_{5R}$ for any fixed $\lambda \in \Sigma_{\ep} \cap K:$
\begin{equation}\label{eq:mid_res_5}
	\left\{\begin{aligned}
& \lambda \eta+\gamma_1 \ol \di\Bu = d  
	    &&\quad\hbox{in}\quad \Omega_{5R}, \\
&\gamma_1 \lambda \Bu-\ol\Di\big( \ol\bS(\Bu) -\gamma_2 \eta \bI \big)
=\bff
		&&\quad\hbox{in}\quad \Omega_{5R},\\
&\big( \ol\bS(\Bu) -\gamma_2 \eta \bI \big) \ol\Bn_{\Gamma} =0
	    &&\quad\hbox{on}\quad \Gamma,\\
	    &\big( \bS(\Bu) -\gamma_2 \eta \bI \big) \Bn_{_{S_{5R}}} =0
	    &&\quad\hbox{on}\quad S_{5R}.
	\end{aligned}\right.
\end{equation}
The result for \eqref{eq:mid_res_5} can be established as follows.
\begin{thm}\label{thm:mid_res_5}
Let $\Omega,\ep,p,r$ be given as in Theorem \ref{thm:mid_res}.
Then there exist continuous linear operators 
$$\big(\CP_{0}(\lambda),\CV_{0}(\lambda)\big) 
\in \CL\big( H^{1,0}_p(\Omega_{5R}) ;
H^{1,2}_p(\Omega_{5R})\big)$$
for any fixed $\lambda \in \Sigma_{\ep} \cap K$
such that $(\eta, \Bu)= \big(\CP_{0}(\lambda), \CV_{0}(\lambda)\big)(d,\bff) \in H^{1,2}_p(\Omega_{5R})$ is a unique solution of \eqref{eq:mid_res_5} 
for any $(d,\bff)\in H^{1,0}_p(\Omega_{5R}).$
Moreover, we have
\begin{equation*}
\|\eta\|_{H^{1}_p(\Omega_{5R})} +
\|\Bu\|_{H^{2}_p(\Omega_{5R})} \leq C_{\lambda}
\big( \|d\|_{H^{1}_p(\Omega_{5R})} 
+ \|\bff\|_{L_p(\Omega_{5R})} \big)
\end{equation*}
for some constant $C_{\lambda}$ depending solely on 
$\lambda,\ep,p,r,\mu,\nu,\gamma_1,\gamma_2,N.$
\end{thm}

To solve \eqref{eq:mid_res_5}, we consider
\begin{equation}\label{eq:mid_res_1}
	\left\{\begin{aligned}
& \lambda \eta+\gamma_1 \di\Bu = d  
	    &&\quad\hbox{in}\quad \Omega_{5R}, \\
&\gamma_1 \lambda \Bu-\Di\big( \bS(\Bu) -\gamma_2 \eta \bI \big)
=\bff
		&&\quad\hbox{in}\quad \Omega_{5R},\\
&\big( \bS(\Bu) -\gamma_2 \eta \bI \big) 
\Bn_{_{\pa \Omega_{5R}}} =\bg
	    &&\quad\hbox{on}\quad \pa \Omega_{5R}.
	\end{aligned}\right.
\end{equation}
Applying Theorem \ref{thm:aux_1}, we have the following theorem for \eqref{eq:mid_res_1}.
\begin{thm}\label{thm:mid_res_1}
Let $\Omega,\ep$ be given as in Theorem \ref{thm:mid_res} 
and $1<p<\infty.$
Set that 
\begin{equation*}
Y_p(\Omega_{5R})= H^1_p(\Omega_{5R}) \times L_p(\Omega_{5R})^N \times H^1_p(\Omega_{5R})^N.
\end{equation*}
Then there exist continuous linear operators 
$$\big(\wt\CP_{0}(\lambda),\wt\CV_{0}(\lambda)\big) 
\in \CL\big( Y_p(\Omega_{5R}) ;
H^{1,2}_p(\Omega_{5R})\big)$$
for any fixed $\lambda \in \Sigma_{\ep} \cap K$
such that $(\eta, \Bu)= \big(\wt\CP_{0}(\lambda),\wt\CV_{0}(\lambda)\big)(d,\bff,\bg) \in H^{1,2}_p(\Omega_{5R})$ is a (unique) solution of \eqref{eq:mid_res_1} 
for any $(d,\bff,\bg)\in Y_p(\Omega_{5R}).$
Moreover, we have
\begin{equation*}
\|\eta\|_{H^{1}_p(\Omega_{5R})} +
\|\Bu\|_{H^{2}_p(\Omega_{5R})} \leq C_{\lambda}
\big( \|(d,\bg)\|_{H^{1}_p(\Omega_{5R})} 
+ \|\bff\|_{L_p(\Omega_{5R})}\big)
\end{equation*}
for some constant $C_{\lambda}$ depending solely on 
$\lambda,\ep,p,\mu,\nu,\gamma_1,\gamma_2,N.$
\end{thm}
\begin{proof}
Now, let us construct $\wt\CP_0(\lambda)$ and $\wt\CV_0(\lambda).$ We reduce \eqref{eq:mid_res_1}
by inserting $\eta= \lambda^{-1} (d-\gamma_1 \di \Bu)$ into 
$\eqref{eq:mid_res_1}_2,$
\begin{equation}\label{eq:mid_res_2}
	\left\{\begin{aligned}
&\lambda \Bu-\Di\big( \alpha\bD(\Bu) +  (\beta-\alpha
+\gamma_2\lambda^{-1})\di \Bu\,\bI \big)=\bF
		&&\quad\hbox{in}\quad \Omega_{5R},\\
&\big( \alpha\bD(\Bu) + (\beta-\alpha
+\gamma_2\lambda^{-1})\di \Bu\bI \big) 
\Bn_{_{\pa \Omega_{5R}}}=\bG
	    &&\quad\hbox{on}\quad \pa \Omega_{5R},
	\end{aligned}\right.
\end{equation}
with 
\begin{equation}\label{eq:aux_1_FG}
(\alpha,\beta)=\gamma_1^{-1}(\mu,\nu),\,\,\,
\bF=\gamma_1^{-1}\bff-\lambda^{-1}\gamma_1^{-1}
\gamma_2\nabla d,
\,\,\,\text{and}\,\,\,
\bG=\gamma_1^{-1}\bg+\lambda^{-1}\gamma_1^{-1}\gamma_2 d\,
\Bn_{_{\pa \Omega_{5R}}}.
\end{equation}
So we shall study the solvability of \eqref{eq:mid_res_2} in what follows.
\smallbreak

According to Theorem \ref{thm:aux_1} and Remark \ref{rmk:aux_1}, there exists 
$\zeta_{\lambda} \in \Xi_{\ep,\lambda,\zeta_0}$ for some $\zeta_0>0$ such that 
\begin{equation}\label{eq:aux_2}
	\left\{\begin{aligned}
&\zeta_{\lambda} \Bu-\Di\big( \alpha\bD(\Bu) +  (\beta-\alpha+\gamma_2\lambda^{-1})\di \Bu\bI \big)=\bF
		&&\quad\hbox{in}\quad \Omega_{5R},\\
&\big( \alpha\bD(\Bu) + (\beta-\alpha
+\gamma_2 \lambda^{-1})\di \Bu\bI \big) 
\Bn_{_{\pa \Omega_{5R}}} =\bG
	    &&\quad\hbox{on}\quad \pa \Omega_{5R},
	\end{aligned}\right.
\end{equation}
admits a solution 
$\CR_{\zeta_{\lambda}} (\bF,\bG)=\Bu\in H^{2}_p(\Omega_{5R})^N$ for any 
$(\bF,\bG)\in L_p(\Omega_{5R})^N \times H^{1}_p(\Omega_{5R})^N.$
Moreover, we have 
\begin{equation}\label{es:aux_2_1}
\|\Bu\|_{H^{2}_p(\Omega_{5R})} \leq C_{\lambda} 
( \|\bF\|_{L_p(\Omega_{5R})} +\|\bG\|_{H^{1}_p(\Omega_{5R})}).
\end{equation}

Next, we rewrite \eqref{eq:aux_2} by 
\begin{equation*}
	\left\{\begin{aligned}
&\lambda \Bu-\Di\big( \alpha\bD(\Bu) 
+  (\beta-\alpha+\gamma_2\lambda^{-1})\di \Bu\bI \big)=\BF+(\lambda-\zeta_{\lambda})\CR_{\zeta_{\lambda}} (\bF,\bG)
		&&\quad\hbox{in}\quad \Omega_{5R},\\
&\big( \alpha\bD(\Bu) + (\beta-\alpha
+\gamma_2\lambda^{-1})\di \Bu\bI \big) 
\Bn_{_{\pa \Omega_{5R}}} =\bG
	    &&\quad\hbox{on}\quad \pa \Omega_{5R},
	\end{aligned}\right.
\end{equation*}
and then denote 
\begin{equation*}
\CT_{\lambda} (\bF,\bG) = \big( (\lambda-\zeta_{\lambda})\CR_{\zeta_{\lambda}}(\bF,\bG),0\big)
\end{equation*}
for any 
$(\bF,\bG)\in L_p(\Omega_{5R})^N \times H^{1}_p(\Omega_{5R})^N.$
By the Rellich's Theorem, we know that $\CT_{\lambda}$ is compact on 
$L_p(\Omega_{5R})^N \times H^{1}_p(\Omega_{5R})^N.$
Furthermore, we claim that $\CI+\CT_{\lambda}$ is invertible from 
$L_p(\Omega_{5R})^N \times H^{1}_p(\Omega_{5R})^N$ into itself.
By the existence of $(\CI+\CT_{\lambda})^{-1},$
$\Bu= \CR_{\zeta_{\lambda}} \circ (\CI+\CT_{\lambda})^{-1}  (\bF,\bG)$
is a solution of \eqref{eq:mid_res_2}. 
According to \eqref{eq:aux_1_FG} and \eqref{es:aux_2_1},  we obtain 
\begin{equation}\label{es:mid_res_1_1}
\|\Bu\|_{H^{2}_p(\Omega_{5R})} 
\leq C_{\lambda}
\big( \|d\|_{H^{1}_p(\Omega_{5R})} + \|\bff\|_{L_p(\Omega_{5R})} \big).
\end{equation}
Then the bound of $\eta$ follows from $\eqref{eq:mid_res_1}_1$ and \eqref{es:mid_res_1_1}.
\smallbreak
 
At last, we prove the claim above. In fact, the existence of $(\CI+\CT_{\lambda})^{-1}$ is reduced to the injectivity of 
$\CI+\CT_{\lambda} $ due to the Fredholm alternative theorem.
So let us consider the homogeneous equation 
\begin{equation}\label{eq:homo_1}
(\CI+\CT_{\lambda}) (\bF_0,\bG_0) 
= \big( \bF_0 +(\lambda-\zeta_{\lambda})
\CR_{\zeta_{\lambda}}(\bF_0,\bG_0), \bG_0\big) =0
\end{equation}
for some$(\bF_0,\bG_0) \in L_p(\Omega_{5R})^N \times H^{1}_p(\Omega_{5R})^N.$
Immediately, we see from \eqref{eq:homo_1} that $\bG_0 =0,$ and that  
$$\bF_0 +(\lambda-\zeta_{\lambda})\CR_{\zeta_{\lambda}}(\bF_0,0)=0.$$  
Now, set $\Bu_0=\CR_{\zeta_{\lambda}}(\bF_0, 0) \in H^{2}_p(\Omega_{5R})^N,$ which satisfies
\begin{equation*}
	\left\{\begin{aligned}
&\lambda \Bu_0-\Di\big( \alpha\bD(\Bu_0) 
+  (\beta-\alpha+\gamma_2\lambda^{-1})\di \Bu_0\bI \big)
=0
		&&\quad\hbox{in}\quad \Omega_{5R},\\
&\big( \alpha\bD(\Bu_0) + (\beta-\alpha
+\gamma_2 \lambda^{-1})\di \Bu_0\bI \big) 
\Bn_{_{\pa \Omega_{5R}}} =0
	    &&\quad\hbox{on}\quad \pa \Omega_{5R}.
	\end{aligned}\right.
\end{equation*}
By the discussion on \eqref{eq:lem_unique_3_2} in the proof of Lemma \ref{lem:unique_3}, we can conclude that $\Bu_0=\bF_0=0.$ 
This completes our proof for the injectivity of $\CI +\CT_{\lambda}.$ 
\end{proof}

The existence in Theorem \ref{thm:mid_res_5} is immediate from the fixed point argument and Theorem \ref{thm:mid_res_1}. 
To handle the uniqueness issue in Theorem \ref{thm:mid_res_5}, we consider some $(\eta, \Bu) \in H^{1,2}_{p}(G)$ satisfying the following system:
\begin{equation}\label{eq:lem_unique_3_1}
	\left\{\begin{aligned}
& \lambda \eta+\gamma_1 \ol\di\Bu = 0  
	    &&\quad\hbox{in}\quad G, \\
&\lambda \gamma_1 \Bu
-\ol \Di\big(\ol\bS(\Bu) -\gamma_2 \eta \bI \big)=0
		&&\quad\hbox{in}\quad G,\\
&\big( \ol\bS(\Bu) -\gamma_2 \eta \bI \big)\ol\Bn_{\pa G} =0
	    &&\quad\hbox{on}\quad \pa G,
	\end{aligned}\right.
\end{equation}
where $G\in \{\Omega_{5R},B_{5R},\Omega\},$ $\Bn_{\pa G}$ is the unit normal vector on the boundary $\pd G,$ and
$\ol\Bn_{\pa G}= (\bI+\bV)\Bn_{\pa G}.$
Noting that $\eta=-\lambda^{-1}\gamma_1 \ol\di \Bu,$  we obtain from 
$\eqref{eq:lem_unique_3_1}_2,$ 
\begin{equation}\label{eq:lem_unique_3_2}
	\left\{\begin{aligned}
&\lambda \Bu-\ol\Di\big( \alpha \ol \bD(\Bu) 
+(\beta-\alpha +\lambda^{-1}\gamma_2) \ol \di \Bu\bI \big)=0
		&&\quad\hbox{in}\quad G,\\
&\big( \alpha \ol\bD(\Bu) 
+(\beta-\alpha +\lambda^{-1}\gamma_2) \ol\di \Bu\bI \big)
\ol \Bn_{\pa G} =0
	    &&\quad\hbox{on}\quad \pa G,
	\end{aligned}\right.
\end{equation}
with $(\alpha,\beta )= \gamma_1^{-1}(\mu,\nu).$ 
For \eqref{eq:lem_unique_3_2}, we can prove the following lemma, which also yields the uniqueness of \eqref{eq:lem_unique_3_1}.
\begin{lem}\label{lem:unique_3}
Let $0<\ep<\pi\slash 2,$ $1<p<\infty,$ 
and $G\in \{\Omega_{5R},B_{5R},\Omega\}.$
Assume that $\Bu \in H^{1,2}_{p}(G)$ solves \eqref{eq:lem_unique_3_2}
for $\lambda \in \Sigma_{\ep} \cap K.$
Then $\Bu(y)=0$ for any $y \in G.$ 
\end{lem}
\begin{proof}
{\bf Step 1.} We first consider the case $p=2.$
Taking the inner product $(\cdot, J\Bu)_{G}$ for $\eqref{eq:lem_unique_3_2}_1$ and integration by parts yield that
\begin{align} \label{eq:lem_unique_3_3}
0=&\lambda \big\|\sqrt{J}\Bu\big\|_{L_2(G)}^2 
+ \frac{\alpha}{2} \int_{G} 
\ol\bD(\Bu):\ol\bD(\Bu) \, J\, dy
+(\beta-\alpha+\gamma_2\lambda^{-1}) 
\big\|\sqrt{J} \,\ol\di \Bu\big\|_{L_2(G)}^2
\\ \nonumber
 =&\lambda A_0 + \alpha A_1
 +(\alpha+\beta+\gamma_2\lambda^{-1}) A_2,
\end{align}
with 
\begin{equation*}
A_0= \big\|\sqrt{J}\Bu\big\|_{L_2(G)}^2  , \quad 
A_1 = \sum_{1\leq i<j\leq N}  \big\| \sqrt{J} \,\ol \bD(\bu)_{ij}\big\|_{L_2(G)}^2,\quad 
A_2=\big\|\sqrt{J} \,\ol\di \Bu\big\|_{L_2(G)}^2.
\end{equation*}
Now we take the real and imaginary parts of \eqref{eq:lem_unique_3_3},
\begin{gather}\label{eq:A012}
\Re \lambda A_0 + \alpha A_1
+(\alpha+\beta+\gamma_2 |\lambda|^{-2} \Re \lambda ) A_2=0,\quad 
\Im \lambda A_0-\gamma_2 |\lambda|^{-2} \Im \lambda  A_2=0.
\end{gather}

To see $\Bu=0$ from \eqref{eq:A012}, it suffices to show $A_0=0.$
If $\Im \lambda=0,$ then $\Re \lambda >0$ as $\lambda \in \Sigma_{\ep}.$
Thus the first equality in \eqref{eq:A012} gives us $A_0=0.$

On the other hand, for $\Im \lambda \not= 0,$ we insert $A_0 =\gamma_2 |\lambda|^{-2} A_2$ into the first equality in \eqref{eq:A012}, 
\begin{equation*}
\alpha A_1+(\alpha+\beta+2\gamma_2 |\lambda|^{-2} \Re \lambda ) A_2=0.
\end{equation*}
Observe that $A_0=0$ is equivalent to $A_2=0$ as $\lambda \not=0.$
\begin{enumerate}
\item ($\Im \lambda \not= 0$ and $\Re \lambda \geq 0$).
If $\Re \lambda \geq 0,$ we have 
$\alpha A_1+(\alpha+\beta) A_2 \leq 0,$
which yields that $A_1=A_2 = 0.$

\item ($\Im \lambda \not= 0$ and $\Re \lambda<0$).
By the definition of $K,$ $\lambda$ fulfils that
\begin{equation*}
 - \frac{(\alpha +\beta)|\lambda|^2}{2\gamma_2} <\Re \lambda <0.
\end{equation*}
Hence $\alpha+\beta+2\gamma_2 |\lambda|^{-2} \Re \lambda>0,$
and then $A_1=A_2=0.$ 
\end{enumerate}
This completes the proof for $p=2.$
\medskip

{\bf Step 2.} Now we assume that $G=\Omega_{5R}$ or $B_{5R}.$
By Step 1 and the boundedness of $G$, we have the uniqueness for $2\leq p<\infty.$ Suppose that $1<p<2.$
We take advantage of the hypoellipticity of \eqref{eq:lem_unique_3_2}
 and the embedding 
$H^1_p(G) \hookrightarrow L_q(G)$
for some $q$ satisfying $0<N(1/p-1/q)<1.$
Then we can prove that $\bu \in H^{2}_2(G)^N$ by applying this idea in finite times.
\medskip

{\bf Step 3.} At last,  we study \eqref{eq:lem_unique_3_2} in $\Omega.$
We will see that $\bu \in H^2_2(\Omega)$
 by the hypoellipticity of \eqref{eq:lem_unique_3_2}.
Recall the definition of $\varphi$ in \eqref{eq:cut-off}, and set $\bw = \varphi \bu.$ Then we have 
\begin{equation}\label{eq:lem_unique_3_4}
	\left\{\begin{aligned}
&\lambda \bw-\ol\Di\big( \alpha \ol \bD(\bw) 
+(\beta-\alpha +\lambda^{-1}\gamma_2) \ol \di \bw\bI \big)=\bff_u
		&&\quad\hbox{in}\quad \Omega_{5R},\\
&\big( \alpha \ol\bD(\bw) 
+(\beta-\alpha +\lambda^{-1}\gamma_2) \ol\di \bw\bI \big)
\ol \Bn_{_{\pa \Omega_{5R}}} =0
	    &&\quad\hbox{on}\quad \pa \Omega_{5R},
	\end{aligned}\right.
\end{equation}
with
\begin{align*}
\bff_u&= -\DV \bU-\big( \alpha \bD(\bu) 
+(\beta-\alpha +\lambda^{-1}\gamma_2) \di \bu\bI \big) 
\nabla \varphi,\\
\bU&=\alpha(\nabla \varphi \otimes \bu + \bu \otimes \nabla \varphi)
+(\beta-\alpha +\lambda^{-1}\gamma_2)  (\nabla \varphi \cdot \bu) \bI.
\end{align*}
Moreover, we have $\bff_u \in H^1_p(\Omega_{5R})^N.$
Then the discussion in Step 2 implies that $\bw=\varphi \bu \in H^2_2(\Omega_{5R})^N$ for all $1<p<\infty.$
That is, $\bu \in H^2_2(\Omega_{b_1})^N$ for $3R<b_0<b_1<4R.$

Next,  we consider $\bv = \psi_{\infty} \bu$ for $\psi_{\infty}$ in \eqref{eq:cut-off}, which satisfies
\begin{align}\label{eq:lem_unique_3_5}
\lambda \bv-\Di\big( \alpha \bD(\bv) 
+(\beta-\alpha +\lambda^{-1}\gamma_2) \di \bv\bI \big)=\bg_u
		&&\quad\hbox{in}\quad \BR^N,
\end{align}
with 
\begin{align*}
\bg_u&= -\DV \wt\bU-\big( \alpha \bD(\bu) 
+(\beta-\alpha +\lambda^{-1}\gamma_2) \di \bu\bI \big) 
\nabla \psi_{\infty},\\
\wt\bU&=\alpha(\nabla \psi_{\infty} \otimes \bu 
+ \bu \otimes \nabla  \psi_{\infty})
+(\beta-\alpha +\lambda^{-1}\gamma_2)  (\nabla  \psi_{\infty} \cdot \bu) \bI.
\end{align*}
Then we have $\supp \bg_u\subset \Omega_{5R}$ and $\|\bg_u\|_{H^1_p(\Omega_{5R})}$ is finite.
Thus we can use Theorem \ref{thm:wh_1} and the argument in Step 2 to get $\bv=\psi_{\infty} \bu \in H^2_2(\BR^N)^N$ for all $1<p<\infty.$
Then we have $\bu \in H^2_2(\Omega \backslash B_{b_1})^N$ by\eqref{eq:cut-off}.
This completes our proof.
\end{proof}


\subsubsection{Solvability of \eqref{resolvent_0} by fixing $\lambda$}
Here we simplify the strategy in Subsection \ref{subsec:para} to construct the solution operators of \eqref{resolvent_0} for $\lambda \in D_{\ep}'.$
\begin{thm}\label{thm:mid_res_2}
Let $\Omega,\ep,$ $p,r$ be given as in Theorem \ref{thm:mid_res}, and $\lambda\in D_{\ep}'.$ 
Then there exist continuous linear operators 
$$\big(\CP_{mid}(\lambda),\CV_{mid}(\lambda)\big) 
\in \CL\big( H^{1,0}_p(\Omega) ;H^{1,2}_p(\Omega)\big)$$
for any fixed $\lambda \in D_{\ep}'$
such that $(\eta, \Bu)= \big(\CP_{mid}(\lambda),\CV_{mid}(\lambda)\big)(d,\bff) \in H^{1,2}_p(\Omega)$ is a unique solution of \eqref{resolvent_0} 
for any $(d,\bff)\in H^{1,0}_p(\Omega).$
Moreover, we have
\begin{equation*}
\|\eta\|_{H^{1}_p(\Omega)} +
\|\Bu\|_{H^{2}_p(\Omega)} \leq C_{\lambda}
\big( \|d\|_{H^{1}_p(\Omega)} + \|\bff\|_{L_p(\Omega)} \big)
\end{equation*}
for some constant $C_{\lambda}$ depending solely on 
$\lambda,\ep,p,r,\mu,\nu,\gamma_1,\gamma_2,N.$
\end{thm}

Thanks to Lemma \ref{lem:unique_3}, we only focus on the constructing $\CP_{mid}(\lambda)$ and $\CV_{mid}(\lambda).$
By Theorem \ref{thm:wh_1} and Theorem \ref{thm:mid_res_5}, we  introduce that for $\lambda \in D_{\ep}',$
\begin{gather*}
(\eta_\lambda,\Bu_{\lambda})= \big(\CP(\lambda),\CV(\lambda)\big)
(\psi_{\infty} d,\psi_{\infty}\bff), \quad 
(\eta_{\sharp},\Bu_{\sharp})= \big(\CP_0(\lambda)\CV_0(\lambda) \big) (\psi_{0} d,\psi_{0}\bff), \\
\wt\eta_{\lambda} =\Phi_{\lambda}(d,\bff)
=(1-\varphi)\eta_{\lambda} +\varphi \, \eta_\sharp,\quad 
\wt \Bu_{\lambda}=\BPsi_{\lambda}(d,\bff) 
=(1-\varphi)\Bu_{\lambda} +\varphi \, \Bu_\sharp,
\end{gather*}
where $\varphi,\psi_{0}$ and $\psi_{\infty}$ are given in \eqref{eq:cut-off}.
Then $\wt\eta_{\lambda}$ and $\wt\Bu_{\lambda}$ fulfil
\begin{equation*}
	\left\{\begin{aligned}
&\lambda \wt\eta_{\lambda}  
+ \gamma_1 \ol\di \, \wt\Bu_{\lambda} 
= d + \CD_{\lambda}(d,\bff)
	    &&\quad\hbox{in}\quad \Omega, \\
& \gamma_1  \lambda \wt\Bu_{\lambda}
-\ol\DV\big( \ol\bS(\wt\Bu_{\lambda}) -\gamma_2 \wt\eta_{\lambda} \bI \big)
=\bff+\CF_{\lambda}(d,\bff)
		&&\quad\hbox{in}\quad \Omega,\\
&\big( \ol\bS(\wt\Bu_{\lambda}) -\gamma_2\wt\eta_{\lambda} \bI \big) 
\ol\Bn_{\Gamma} =0
	    &&\quad\hbox{on}\quad \Gamma,
	\end{aligned}\right.
\end{equation*}
where $\bV_{\lambda}$ is defined in \eqref{eq:V_df}, and 
\begin{align*}
\CD_{\lambda}(d,\bff)=
& \gamma_1 \nabla \varphi \cdot (\Bu_{\sharp}-\Bu_{\lambda} ),\\
\CF_{\lambda}(d,\bff)=&
\Big( \big(\bS(\Bu_{\lambda}) -\gamma_2\eta_{\lambda}\bI\big)
-\big(\bS(\Bu_{\sharp}) -\gamma_2\eta_{\sharp}\bI\big) \Big)
 \nabla \varphi -\Di\bV_{\lambda}(d,\bff).
\end{align*}
Moreover, we have
\begin{gather*}
\CQ_{\lambda}(d,\bff)=
(\CD_{\lambda},\CF_{\lambda})(d,\bff) \in H^{1,0}_{p}(\Omega), \quad
\forall \,\,\,\lambda \in D_{\ep}', 
\end{gather*}
which is a compact operator on $H^{1,0}_{p}(\Omega)$ in view of the Rellich compactness theorem and 
$$\supp\CQ_{\lambda}(d,\bff) \subset \supp \nabla \varphi \subset 
D_{b_1,b_2}=\{x\in \Omega :b_1\leq |x|\leq  b_2\}.$$ 
In addition, $\CI +\CQ_{\lambda}$ is invertible due to the following lemma.
\begin{lem}\label{lem:unique_5}
For any $\lambda \in D_{\ep}',$ $\CI + \CQ_{\lambda}$ has a bounded inverse $(\CI+ \CQ_\lambda)^{-1}$ in $\CL\big(H^{1,0}_{p}(\Omega)\big).$ 
\end{lem}
\begin{proof}
According to the compactness of $\CQ_{\lambda}$ and the Fredholm's alternative theorem, it is sufficient to verify the injectivity of $\CI + \CQ_{\lambda}.$
For any $(d,\bff)\in \Ker (\CI + \CQ_{\lambda}) 
\subset H_{p}^{1,0}(\Omega),$
$\supp d,\,\supp \bff \subset D_{b_1,b_2}.$
Thus $$(d,\bff)= (\psi_{\infty} d,\psi_{\infty}\bff)= (\psi_{0} d,\psi_{0}\bff).$$
By Theorem \ref{thm:wh_1} and Theorem \ref{thm:mid_res_5},
we denote 
\begin{gather*}
(\eta_{\lambda},\Bu_{\lambda})=\big(\CP(\lambda),\CV(\lambda)\big)(d,\bff),
\,\,\,\text{and}\,\,\,
(\eta_\sharp,\Bu_\sharp)=\big(\CP_0(\lambda),\CV_0(\lambda)\big)(d,\bff)
\end{gather*}
for all $\lambda \in D_{\ep}'.$
Then $(\wt\eta_{\lambda},\wt\bu_{\lambda})
=(1-\varphi)(\eta_{\lambda},\Bu_{\lambda}) 
+\varphi (\eta_\sharp,\Bu_\sharp) $
satisfies 
\begin{equation*}
	\left\{\begin{aligned}
&\lambda \wt\eta_{\lambda}  
+ \gamma_1 \ol\di\, \wt \bu_{\lambda} = 0
	    &&\quad\hbox{in}\quad \Omega, \\
&\gamma_1 \lambda  \wt\bu_{\lambda}
-\ol\Di\big( \,\ol\bS(\wt\bu_{\lambda}) 
-\gamma_2 \wt\eta_{\lambda} \bI \big)=0
		&&\quad\hbox{in}\quad \Omega,\\
&\big( \,\ol\bS(\wt\Bu_{\lambda}) -
\gamma_2 \wt\eta_{\lambda} \bI \big) \ol\Bn_{\Gamma} =0
	    &&\quad\hbox{on}\quad \Gamma.
	\end{aligned}\right.
\end{equation*}
Hence Lemma \ref{lem:unique_3} yields that
\begin{equation}\label{eq:vanish_2}
(\wt{\eta}_{\lambda},\wt{\Bu}_{\lambda})
=(1-\varphi)(\eta_{\lambda},\Bu_{\lambda}) 
+\varphi (\eta_\sharp,\Bu_\sharp) =(0,0)
\quad \hbox{in}\,\,\, \Omega,
\end{equation}
and that
\begin{gather*}
\eta_{\sharp} (x) =0, \,\,\, \Bu_{\sharp}=0\,\,\, \text{for}\,\,\, |x| \leq b_1,\quad 
\eta_{\lambda} (x) =0, \,\,\, \Bu_{\lambda}=0\,\,\, \text{for}\,\,\, |x| \geq b_2.
\end{gather*}

Next, set that
$$(\wt\eta_\sharp, \wt\bu_\sharp)(x)
=\begin{cases} (\eta_\sharp, \bu_\sharp)(x)&
\quad\text{for $x \in \Omega_{5R}\setminus B_{2R}$},
\\ 0 &\quad\text{for $x \in B_{2R}$}.
\end{cases}$$
Then $(\rho,\bv)=(\wt\eta_{\sharp}-\eta_{\lambda},
\wt\bu_{\sharp}-\Bu_{\lambda})$ satisfies
\begin{equation*}
	\left\{\begin{aligned}
& \lambda \rho+\gamma_1 \di\,\bv = 0
	    &&\quad\hbox{in}\quad B_{5R}, \\
&\gamma_1 \lambda \Bv-\Di\big( \bS(\Bv) -\gamma_2 \rho \bI \big)
=0
		&&\quad\hbox{in}\quad B_{5R},\\
	    &\big( \bS(\bv) -\gamma_2 \rho \bI \big) \Bn_{_{S_{5R}}} =0
	    &&\quad\hbox{on}\quad S_{5R}.
	\end{aligned}\right.
\end{equation*}
By modifying the proof of Lemma \ref{lem:unique_3}, we obtain that 
$(\rho,\bv)=(0,0)$ in $B_{5R}.$
In particular, 
$$(\eta_{\sharp},\Bu_{\sharp})
=(\eta_{\lambda},\Bu_{\lambda})
\,\,\, \hbox{in}\,\,\,D_{b_1,b_2}.$$
Then \eqref{eq:vanish_2} implies that 
$(\eta_{\lambda},\Bu_{\lambda})=(0,0)$ in 
$D_{b_1,b_2},$ and $(d,\bff)=(0,0)$ for 
$\supp d,\,\supp \bff \subset D_{b_1,b_2}.$
This completes our proof.
\end{proof}

\begin{proof}[Completing the proof of Theorem \ref{thm:mid_res}]
By Theorem \ref{thm:mid_res_2}, we have 
\begin{equation*}
\big\|\big(\CP_{mid}(\lambda),\CV_{mid}(\lambda)\big)(d,\bff)\big\|_{H^{1,2}_p(\Omega)}
\leq C_{\lambda} \|(d,\bff)\|_{H^{1,0}_p(\Omega)}
\end{equation*}
for any $\lambda \in D_{\ep}'.$ Then 
\begin{equation*}
\big\|\big(\CP_{mid}(\zeta),\CV_{mid}(\zeta)\big)(d,\bff)\big\|_{H^{1,2}_p(\Omega)}\leq C_{\lambda} \|(d,\bff)\|_{H^{1,0}_p(\Omega)}
\end{equation*}
for any $\zeta \in B_{r_{\lambda}}(\lambda)$ for some $r_{\lambda}>0.$
As $D_{\ep}'$ is compact, we can choose finite $\lambda_1,\dots, \lambda_{N_0} \in D_{\ep}' $ such that $D_{\ep}\subset \cup_{k=1}^{N_0} B_{r_{\lambda_k}}(\lambda_{k}).$ Then we have 
\begin{equation*}
\big\|\big(\CP_{mid}(\lambda),\CV_{mid}(\lambda)\big)(d,\bff)\big\|_{H^{1,2}_p(\Omega)}\leq C \|(d,\bff)\|_{H^{1,0}_p(\Omega)}
\end{equation*}
for any $\lambda \in D_{\ep}'.$
\end{proof}


\section{$L_p$-$L_q$ estimates of the linearized problem}
In this section, we will prove $L_p$-$L_q$ estimates of \eqref{LL:Lame_1}.
In next subsection, we show Theorem \ref{thm:local_energy} by using the results in sections \ref{sec:near} and \ref{sec:large}.
Then we review a preliminary result for some model problem in the whole space $\BR^N,$ which plays the key role in the proof of Theorem \ref{thm:main_LpLq}. At last, we derive the $L_p$-$L_q$ decay estimates for \eqref{LL:Lame_1} in subsection \ref{subsec:LpLq}.

\subsection{Local energy estimates}

This short subsection is dedicated to the proof of Theorem \ref{thm:local_energy}. Without loss of generality, we take $L=5R.$
Recall the definition $K_{\ep}$ and $D_{\ep},$ and choose $0 < \ep < \ep'<\pi \slash 2$ such that 
$\lambda_1' e^{i(\pi-\ep')} \in K_{\ep}.$ 
Then we introduce 
$$\Gamma=\Gamma_{+,2} \cup \Gamma_{+,1} \cup \Gamma_{+,0} \cup 
D_+ \cup D_- \cup \Gamma_{-,0} \cup \Gamma_{-,1} \cup \Gamma_{-,2},$$ 
where $\Gamma_{\pm,k},$ $k=0,1,2,$ $D_{\pm}$ are defined by
\begin{gather*}
\Gamma_{+,2} =\{\lambda = r e^{i(\pi-\ep')}, r: \infty \rightarrow \lambda_2'\},
\quad  \Gamma_{+,1}=\{\lambda = r e^{i(\pi-\ep')}, r: \lambda_2' \rightarrow \lambda_1'\},\\
\Gamma_{+,0} =\{\lambda = 
\lambda_1' \cos (\pi-\ep') + i s, s: \lambda_1' \sin(\pi-\ep') \rightarrow 0\},\\
D_{+}=\{\lambda = s e^{i\pi}, s: \lambda_1' \cos \ep' \rightarrow 0\}, \quad 
 D_{-}=\{\lambda = s e^{-i\pi}, s: 0 \rightarrow \lambda_1' \cos \ep' \},\\
 \Gamma_{-,0} =\{\lambda = 
\lambda_1' \cos (-\pi+\ep') + i s, s: 0\rightarrow \lambda_1' \sin(-\pi+\ep') \},\\
\Gamma_{+,1}=\{\lambda = r e^{i(-\pi+\ep')}, r:  \lambda_1'\rightarrow \lambda_2'\},\quad 
\Gamma_{+,2} =\{\lambda = r e^{i(-\pi+\ep')}, r:  \lambda_2'\rightarrow \infty\}.
\end{gather*}
Next, we denote that
\begin{align*}
T(t)(\rho_0,\Bv_0)=\frac{1}{2\pi i} \int_{\Gamma}e^{\lambda t} (\lambda \CI +\CA_{\Omega})^{-1} (\rho_0,\Bv_0) \,d\lambda 
=\sum_{k=0}^2 I_{k}(t) +J(t),
\end{align*}
with 
\begin{align*}
I_k(t)&=\frac{1}{2\pi i} \int_{\Gamma_{+,k} \cup \Gamma_{-,k}}e^{\lambda t} (\lambda \CI +\CA_{\Omega})^{-1} (\rho_0,\Bv_0) \,d\lambda, \,\,\, k=0,1,2,\\
J(t)&=\frac{1}{2\pi i} \int_{D_{+} \cup D_{-}}e^{\lambda t} (\lambda \CI +\CA_{\Omega})^{-1} (\rho_0,\Bv_0) \,d\lambda.
\end{align*}

By Theorem \ref{thm:mid_res} and Theorem \ref{thm:large_res}, 
it is not hard to see that  
\begin{align*}
\|I_{k'}(t)\|_{H^{1,2}_p(\Omega)} 
\leq C_{\ep'} e^{-(\lambda_{k'}'\cos \ep') t} 
\| (\rho_0,\Bv_0)\|_{H^{1,0}_p(\Omega)}, 
\,\,\, \forall\,\, t\geq 1, \,\, k'=1,2.
\end{align*}
Thank to Theorem \ref{thm:ext_1}, $I_0$ is bounded by
\begin{align*}
\|I_0(t)\|_{H^{1,2}_p(\Omega_{5R})} 
\leq C_{\ep'} e^{-(\lambda_1'\cos \ep') t} 
\| (\rho_0,\Bv_0)\|_{H^{1,0}_p(\Omega)}, \,\,\, \forall\,\, t\geq 1.
\end{align*}

Now, let us study $J(t).$ According to Theorem \ref{thm:ext_1}, 
we can rewrite $J(t)=J_0(t)+J_1(t)+J_2(t)$ with
\begin{align*}
J_0(t)&=\frac{1}{2\pi i} \int_{D_{+} \cup D_{-}}e^{\lambda t} 
\lambda^{N\slash 2-1} (\log \lambda)^{\sigma(N)} (0,\BV_{\lambda}^0) (\rho_0,\Bv_0) \,d\lambda,\\
J_1(t)&=\frac{1}{2\pi i} \int_{D_{+} \cup D_{-}}e^{\lambda t} 
\lambda^{N-2} \log \lambda (\BM_{\lambda}^1,\BV_{\lambda}^1) (\rho_0,\Bv_0) \,d\lambda,\\
J_2(t)&=\frac{1}{2\pi i} \int_{D_{+} \cup D_{-}}e^{\lambda t} (\BM_{\lambda}^2,\BV_{\lambda}^2) (\rho_0,\Bv_0) \,d\lambda.
\end{align*}
By the analyticity of operator $(\BM_{\lambda}^2,\BV_{\lambda}^2)$ near 
$D_+ \cup D_-,$ we have $J_2(t)=0.$ Moreover, direct calculations yield 
\begin{align*}
\|J_0(t)\|_{H^{1,2}_p(\Omega_{5R})} 
&\leq C_{\ep'} 
t^{-N/2}
\| (\rho_0,\Bv_0)\|_{H^{1,0}_p(\Omega)},\\
\|J_1(t)\|_{H^{1,2}_p(\Omega_{5R})} 
&\leq C_{\ep'} t^{-(N-1)}
\| (\rho_0,\Bv_0)\|_{H^{1,0}_p(\Omega)}, 
\,\,\, \forall\,\, t\geq 1.
\end{align*}
Summing up all the bounds, we obtain that 
\begin{equation*}
\|T(t)(\rho_0,\Bv_0)\|_{H^{1,2}_p(\Omega_{5R})} 
\leq C_{\ep'} t^{-N/2}
\| (\rho_0,\Bv_0)\|_{H^{1,0}_p(\Omega)}, 
\,\,\, \forall\,\, t\geq 1\,\,\,\text{and}\,\,\, N\geq 3.
\end{equation*}

At last, note that 
\begin{align*}
\pa_t^k T(t)(\rho_0,\Bv_0)=\frac{1}{2\pi i} 
\int_{\Gamma}e^{\lambda t} \lambda^k (\lambda \CI +\CA_{\Omega})^{-1} (\rho_0,\Bv_0) \,d\lambda.
\end{align*}
By the similar calculations, we can obtain the desired estimates for the higher derivatives in Theorem \ref{thm:local_energy}.

\subsection{Some result for the problem in $\BR^N$}
\begin{equation}\label{eq:GLame_whole}
	\left\{\begin{aligned}
&\pa_t  \eta  + \gamma_1 \di \Bu = 0  
	    &&\quad\hbox{in}\quad  \BR^N \times \BR_+, \\
&\gamma_1 \pa_t \Bu-\Di\big( \bS(\Bu) -\gamma_2 \eta \bI \big)=0
		&&\quad\hbox{in}\quad  \BR^N \times \BR_+,\\
		&(\eta, \Bu)|_{t=0} = (\eta_0, \Bu_0) 
        &&\quad\hbox{in}\quad \BR^N.
	\end{aligned}\right.
\end{equation}
For \eqref{eq:GLame_whole}, we recall the following result proved in \cite{KS2002}.
\begin{thm}{\cite[Theorems 2.3 and 2.4]{KS2002}}
\label{thm:GLame_whole}
There exist $\eta^0, \eta^{\infty},$ $\Bu^0,$ and $\Bu^{\infty}$ such that 
$(\eta,\Bu)= (\eta^0+\eta^{\infty},\Bu^0+\Bu^{\infty})$ 
solves \eqref{eq:GLame_whole}, and the following assertions hold for non-negative integers $\ell,m$ and $n.$
\begin{enumerate}
\item For all $t\geq 1,$ there exists a $C=C_{m,\ell,p,q}>0$ such that 
\begin{equation*}
\sum_{|\alpha|=\ell} \|\pa_t^m \pa_x^{\alpha}(\eta^0,\Bu^0)(t)\|_{L_p(\BR^N)}
\leq C t^{-N(1/q-1/p)/2-(m+\ell)/2} \|(\eta_0,\Bu_0)\|_{L_q(\BR^N)}
\end{equation*}
with $1\leq q \leq 2 \leq p \leq \infty.$
Moreover, 
\begin{equation*}
\sum_{|\alpha|=\ell} \|\pa_t^m \pa_x^{\alpha}(\eta^0,\Bu^0)(t)\|_{L_p(\BR^N)}
\leq C \|(\eta_0,\Bu_0)\|_{L_q(\BR^N)}, \,\,\, \forall\,\, 0<t\leq 2,
\end{equation*}
for $1\leq q \leq p \leq \infty$ and $(p,q) \not=(1,1), (\infty,\infty).$
\item Let $(\ell)^{+}=\ell$ if $\ell \geq 0$ and $(\ell)^+=0$ if $\ell<0.$
Let $1<p<\infty.$ Then, there exist positive constants $C$ and $c$ such that
for any $t>0,$
\begin{align*}
\sum_{|\alpha|=\ell} \|\pa_t^m \pa_x^{\alpha}\eta^\infty(t)\|_{L_p(\BR^N)}
\leq C e^{-ct}
&\Big( t^{-n/2}\|\eta_0\|_{H_p^{(2m+\ell-n-2)^+}(\BR^N)}
+\|\eta_0\|_{H^{\ell}_p(\BR^N)}\\
&  +t^{-n/2}\|\Bu_0\|_{H_p^{(2m+\ell-n-1)^+}(\BR^N)}
+\|\Bu_0\|_{H^{(\ell-1)^+}_p(\BR^N)} \Big),
\end{align*}
\begin{align*}
\sum_{|\alpha|=\ell} \|\pa_t^m \pa_x^{\alpha}\Bu^\infty(t)\|_{L_p(\BR^N)}
\leq C e^{-ct}
&\Big( t^{-n/2}\|\eta_0\|_{H_p^{(2m+\ell-n-1)^+}(\BR^N)}
+\|\eta_0\|_{H^{(\ell-1)^+}_p(\BR^N)}\\
&  +t^{-n/2}\|\Bu_0\|_{H_p^{(2m+\ell-n)^+}(\BR^N)}
+\|\Bu_0\|_{H^{(\ell-2)^+}_p(\BR^N)} \Big).
\end{align*}
\end{enumerate}
\end{thm}

\subsection{Proof of Theorem \ref{thm:main_LpLq}} \label{subsec:LpLq}
To obtain the decay estimates of $(\rho,\Bv)=T(t) (\rho_0,\Bv_0),$ 
we borrow the cut-off functions $\varphi,$ $\psi_0$ and $\psi_{\infty}$ in \eqref{eq:cut-off}. The proof is mainly divided into two steps, 
where we treat the bound of $(\rho,\Bv)$ in $\Omega_{5R}$ and $\Omega \backslash B_{5R}$ respectively. 
For simplicity, we only check the $L_p$-$L_q$ decay estimates for $t\geq 3$ with the large time behaviour of the solutions involved. Afterwards, the bound for all $t\geq 1$ can be obtained by refining the constant $C$ in Theorem \ref{thm:main_LpLq}. 
\medskip

{\bf Step 1.} In this part, we would like to verify the following estimate in $\Omega_{5R},$
\begin{equation}\label{es:GL_le_1}
\|T(t) (\rho_0,\Bv_0)\|_{H^{1,2}_p(\Omega_{5R})} 
\leq C t^{-N/(2q)} \vertiii{(\rho_0, \Bv_0)}_{p,q}, 
\quad \forall\,\, t\geq 3.
\end{equation}
To prove \eqref{es:GL_le_1}, we introduce $(\eta,\Bu)$ by solving the following system:
\begin{equation*}
	\left\{\begin{aligned}
&\pa_t  \eta  + \gamma_1 \di \Bu = 0  
	    &&\quad\hbox{in}\quad  \BR^N \times \BR_+, \\
&\gamma_1 \pa_t \Bu-\Di\big( \bS(\Bu) -\gamma_2 \eta \bI \big)=0
		&&\quad\hbox{in}\quad  \BR^N \times \BR_+,\\
		&(\eta, \Bu)|_{t=0} = (\psi_{\infty}\rho_0, \psi_{\infty}\Bv_0) 
        &&\quad\hbox{in}\quad \BR^N.
	\end{aligned}\right.
\end{equation*}
In addition, Theorem \ref{thm:GLame_whole} and the definition of $\psi_{\infty}$ in \eqref{eq:cut-off} imply that
\begin{equation}\label{es:GL_etau_1}
\begin{aligned}
\|(\eta,\Bu)(t)\|_{H^{1,2}_p(\Omega_{5R})} &\leq C_{p,R} 
\big( \|(\eta^0,\Bu^0)(t)\|_{H^{1,2}_\infty(\BR^N)} +
\|(\eta^{\infty},\Bu^{\infty})(t)\|_{H^{1,2}_p(\BR^N)} \big)\\ 
& \leq C_{p,q,R}\, t^{-N/(2q)} \vertiii{(\rho_0, \Bv_0)}_{p,q}, 
\quad \forall\,\, t\geq 1.
\end{aligned}
\end{equation}
\smallbreak

Next, $(\theta,\Bw)=(\rho,\Bv)-(1-\varphi)(\eta,\Bu)$ fulfils that
\begin{equation}\label{eq:GLame_2}
	\left\{\begin{aligned}
&\pa_t \theta + \gamma_1 \ol\di\, \Bw =d_u
	    &&\quad\hbox{in}\quad \Omega \times \BR_+, \\
&\gamma_1 \pa_t \Bw-\ol\Di \big(\ol\bS(\Bw)-\gamma_2
\theta \bI\big)=\bff_{\eta,u}
		&&\quad\hbox{in}\quad \Omega \times \BR_+,\\
&\big(\ol\bS(\Bw)-\gamma_2\theta \bI\big) \ol\Bn_{\Gamma} = 0
	    &&\quad\hbox{on}\quad \Gamma \times \BR_+, \\
&(\theta, \Bw)|_{t=0} = (\varphi\rho_0, \varphi\Bv_0) 
        &&\quad\hbox{in}\quad \Omega,
	\end{aligned}\right.
\end{equation}
with
\begin{align*}
d_u&=\gamma_1 \nabla \varphi \cdot \bu, \quad
\bff_{\eta,u} =-\big( \bS(\Bu)-\gamma_2 \eta \bI \big) \nabla \varphi
+\Di \bU,\\
\bU&=-\mu (\Bu \otimes \nabla \varphi
+\nabla \varphi\otimes \Bu) 
-(\nu-\mu) (\nabla \varphi \cdot \Bu) \bI,
\end{align*} 
where  we have used the fact that $(1-\varphi) \bV =0$ from \eqref{support_0} and \eqref{eq:cut-off}.
Clearly, 
$$\supp (d_u, \bff_{\eta,u},\bU) \subset D_{b_1,b_2}\subset \Omega_{5R}.$$ 
Moreover,  \eqref{es:GL_etau_1} yields that 
$ (d_u, \bff_{u,\eta}) \in H^{1,1}_p(\Omega_{5R})$ with
\begin{equation}\label{es:GL_df_1}
\|(d_u,\bff_{u,\eta}) (t)\|_{H^{1,k}_p(\Omega_{5R})} 
\leq C \|(\eta,\Bu) (t)\|_{H^{k,1+k}_p(\Omega_{5R})} 
\leq  C_{p,q,R}\, t^{-N/(2q)} \vertiii{(\rho_0, \Bv_0)}_{p,q}, 
\end{equation}
for $k=0,1,$ and $t\geq 1.$
\medskip

By the definition of $(\theta,\Bw)$ and \eqref{es:GL_etau_1}, 
\eqref{es:GL_le_1} holds true if one can show
\begin{equation}\label{es:GL_le_2}
\|(\theta, \Bw)(t)\|_{H^{1,2}_p(\Omega_{5R})} 
\leq C t^{-N/(2q)} \vertiii{(\rho_0, \Bv_0)}_{p,q}, 
\quad \forall\,\, t\geq 3.
\end{equation}
To obtain \eqref{es:GL_le_2}, we use the Duhamel principle,
\begin{equation}\label{eq:thetaw_1}
(\theta,\Bw)(t)= T(t)(\varphi\rho_0, \varphi\Bv_0) 
+ \int_0^t T(t-s)(d_u,\bff_{\eta,u})(s) \,ds.
\end{equation}
As $(\varphi\rho_0, \varphi\Bv_0) \in X_{p,5R}(\Omega),$ Theorem \ref{thm:local_energy} yields that 
\begin{align}\label{es:GL_I_0}
\|T(t)(\varphi\rho_0, \varphi\Bv_0)\|_{H^{1,2}_p(\Omega_{5R})} 
\leq C t^{-N/2} \|(\varphi\rho_0, \varphi\Bv_0)\|_{H^{1,0}_p(\Omega)}
\leq C t^{-N/2} \|(\rho_0,\Bv_0)\|_{H^{1,0}_p(\Omega)}
\end{align}
for all $t\geq 1.$
To bound the second term on the right-hand side of \eqref{eq:thetaw_1}, we write 
\begin{align*}
\int_0^t T(t-s)(d_u,\bff_{\eta,u})(s) \,ds
&= \Big(\int_0^1 +\int_1^{t-2} +\int_{t-2}^t\Big)
T(t-s)(d_u,\bff_{\eta,u})(s) \,ds\\
&=I_1(t) +I_2(t)+I_3(t),
\quad \forall\,\, t\geq 3,
\end{align*}
and then we study $I_k(t)$ $k=1,2,3,$ respectively in what follows.
\medskip

\underline{Bound of $I_1(t)$.}
Firstly, Theorem \ref{thm:GLame_whole} yields that 
\begin{align}\label{es:GL_etau_2}
\|(\eta,\Bu)(t)\|_{H^{0,1}_p(\BR^N)}
&\leq C( t^{-1/2}\| (\psi_{\infty}\rho_0, \psi_{\infty}\Bv_0) \|_{L_p(\BR^N)}
+ \| (\psi_{\infty}\rho_0, \psi_{\infty}\Bv_0) \|_{L_q(\BR^N)})
\\ \nonumber
&\leq C t^{-1/2}\vertiii{(\rho_0, \Bv_0)}_{p,q}\quad \text{for}\,\,0<t<1.
\end{align}
Then, by Theorem \ref{thm:local_energy} and \eqref{es:GL_etau_2}, 
we have 
\begin{align} \label{es:GL_I_1}
\|I_1(t)\|_{H^{1,2}_p(\Omega_{5R})} 
&\leq C \int_0^1 (t-s)^{-N/2}\|(d_u,\bff_{\eta,u})(s)\|_{H^{1,0}_p(\Omega_{5R})} \,ds\\ \nonumber
&\leq C \int_0^1 (t-s)^{-N/2}\|(\eta,\Bu)(s)\|_{H^{0,1}_p(\Omega_{5R})} \,ds\\ \nonumber
& \leq C t^{-N/2} \vertiii{(\rho_0, \Bv_0)}_{p,q}, 
\quad \forall\,\, t\geq 3.
\end{align}
\smallbreak

\underline{Bound of $I_2(t)$.}
According to Theorem \ref{thm:local_energy} and \eqref{es:GL_df_1}, we have 
\begin{align}\label{es:GL_I_2_1}
\|I_2(t)\|_{H^{1,2}_p(\Omega_{5R})} 
\leq C \int_1^{t-2} (t-s)^{-N/2} s^{-N/(2q)} \,ds
\vertiii{(\rho_0, \Bv_0)}_{p,q}.
\end{align}
Note that 
\begin{align*}
\int_1^{t-2} (t-s)^{-N/2} s^{-N/(2q)} \,ds 
&\leq  \Big(\int_1^{t/3} +\int_{t/3}^{t-2}\Big)
 (t-s)^{-N/2} s^{-N/(2q)} \,ds \\
& \leq  \int_1^{t/3} \big((t+s)/2\big)^{-N/2} s^{-N/(2q)} \,ds\\
&\quad + \int_{t/3}^{t-2} (t-s)^{-N/2} \big((t+s)/4\big)^{-N/(2q)} \,ds
=A_1(t) +A_2 (t).
\end{align*}

For the bound of  $A_1(t),$
we first assume that $N\slash (2q)>1,$ and obtain
\begin{equation*}
A_1(t) \leq C_N t^{-N/2} \int_1^{t/3}s^{-N/(2q)} \,ds \leq C_{q,N} t^{-N/2}.
\end{equation*}
Otherwise, if $0<N/(2q)\leq 1$ for some $N\geq 3,$ then we have 
\begin{align*}
A_1(t)\leq C_{N} \int_1^{t/3} \frac{s^{1+\sigma-N/(2q)}}{(t+s)^{N/2}} s^{-(1+\sigma)} \,ds \leq C_{\sigma,q,N}  t^{-(N/2-1-\sigma+N/(2q))},
\end{align*}
with $0<\sigma<1/2.$ Thus we have 
\begin{equation*}
A_1(t)  \leq C_{q,N} t^{-N/(2q)} \,\,\,
\text{for all $1\leq q\leq 2$ and $N\geq 3.$} 
\end{equation*}

On the other hand, 
\begin{equation*}
A_2(t) \leq C_{q,N}t^{-N/(2q)} \int_{t/3}^{t-2} (t-s)^{-N/2}  \,ds \leq  C_{q,N}t^{-N/(2q)}.
\end{equation*}
Now, combining the bounds of $A_1(t)$ and $A_2(t)$ above, we obtain
\begin{equation*}
\int_1^{t-2} (t-s)^{-N/2} s^{-N/(2q)} \,ds  \leq  C_{q,N}t^{-N/(2q)},
\end{equation*}
which, together with \eqref{es:GL_I_2_1}, implies that
\begin{equation}\label{es:GL_I_2_2}
\|I_2(t)\|_{H^{1,2}_p(\Omega_{5R})}  \leq C t^{-N/(2q)} \vertiii{(\rho_0, \Bv_0)}_{p,q}, 
\quad \forall\,\, t\geq 3.
\end{equation}

\underline{Bound of $I_3(t)$.}
By Theorem \ref{thm:sg}, we have
\begin{equation*}
\|T(t)(d,\bff)\|_{H^{1,2}_p(\Omega)}  \leq 
\begin{cases}
C t^{-1} \|(d,\bff)\|_{H^{1,0}_p(\Omega)},
& \forall \,\, (d,\bff) \in  H^{1,0}_{p}(\Omega),\\
C \|(d,\bff)\|_{H^{1,2}_{p}(\Omega)},
&\forall \,\, (d,\bff) \in  \CD_{p}(\CA_\Omega),
\end{cases}
\end{equation*}
for $0<t<2.$ Then by the interpolation, we see that 
\begin{equation}\label{es:GL_sg_1}
\|T(t)(d,\bff)\|_{H^{1,2}_p(\Omega)} \leq C t^{-1/2} \|(d,\bff)\|_{H^{1,1}_p(\Omega)}.
\end{equation}
Now, we infer from \eqref{es:GL_sg_1} and \eqref{es:GL_df_1},
\begin{equation}\label{es:GL_I_3}
\begin{aligned}
\|I_3(t)\|_{H^{1,2}_p(\Omega_{5R})} 
&\leq C \int_{t-2}^{t} (t-s)^{-1/2}\|(d_u,\bff_{\eta,u})(s)\|_{H^{1,1}_p(\Omega_{5R})} \,ds\\ 
&\leq C \int_{t-2}^{t} (t-s)^{-1/2} s^{-N/(2q)} \,ds
\,\vertiii{(\rho_0, \Bv_0)}_{p,q} \\ 
& \leq C t^{-N/(2q)} \vertiii{(\rho_0, \Bv_0)}_{p,q}, 
\quad \forall\,\, t\geq 3.
\end{aligned}
\end{equation}

Finally, we can conclude the claim \eqref{es:GL_le_2} by summing up \eqref{es:GL_I_0}, \eqref{es:GL_I_1}, \eqref{es:GL_I_2_2} and \eqref{es:GL_I_3}.  This completes the proof of \eqref{es:GL_le_1}.
\medskip

{\bf Step 2.} Let us study the bound of $(\rho,\Bv)=T(t) (\rho_0,\Bv_0)$ near the area $\Omega\backslash B_{5R}$ 
by holding \eqref{es:GL_etau_1}.
Recall the definition of $\psi_\infty$ in \eqref{eq:cut-off}. Then $(\wt\rho,\wt\Bv)=(\psi_{\infty} \rho, \psi_{\infty}\Bv)$ can be regarded as a solution of 
\begin{equation}\label{eq:GLame_3}
	\left\{\begin{aligned}
&\pa_t  \wt \rho  + \gamma_1 \di \wt\Bv = \wt d  
	    &&\quad\hbox{in}\quad  \BR^N \times \BR_+, \\
&\gamma_1 \pa_t \wt\Bv-\Di\big( \bS(\wt\Bv) -\gamma_2 \wt\rho \bI \big)=\wt\bff
		&&\quad\hbox{in}\quad  \BR^N \times \BR_+,\\
		&(\wt\rho, \wt\Bv)|_{t=0} = (\psi_{\infty}\rho_0, \psi_{\infty}\Bv_0) 
        &&\quad\hbox{in}\quad \BR^N,
	\end{aligned}\right.
\end{equation}
with 
\begin{gather*}
\wt d = \gamma_1 \nabla \psi_{\infty} \cdot \Bv, \quad 
\wt \bff =-\big( \bS(\Bv) -\gamma_2 \rho \bI \big) \nabla \psi_{\infty}
-\Di \wt\bV\\
\wt \bV=\mu (\Bv \otimes \nabla \psi_{\infty} 
+\nabla \psi_{\infty} \otimes \Bv)
+(\nu-\mu) (\nabla \psi_{\infty} \cdot \Bv) \, \bI.
\end{gather*}
Moreover, we have $\supp (\wt d, \wt \bff) \subset \Omega_{5R}$ and
\begin{align*}
\vertiii{(\wt d,\wt \bff)(t)}_{p,q}  \leq C 
\|(\rho,\Bv)(t)\|_{H^{0,1}_p(\Omega_{5R}) \cap H^{0,1}_q(\Omega_{5R})} 
\leq C_{p,q,R}\|(\rho,\Bv)(t)\|_{H^{0,1}_p(\Omega_{5R})} 
\end{align*}
for any $t>0$ and $1\leq q \leq 2 \leq p < \infty.$ 
Similar to \eqref{es:GL_sg_1}, we can conclude 
\begin{equation}\label{es:GL_sg_2}
\|T(t)(\rho_0,\Bv_0)\|_{H^{1,1}_p(\Omega)} \leq 
C t^{-1/2} \|(\rho_0,\Bv_0)\|_{H^{1,0}_p(\Omega)} 
\,\,\, (0<t<3)
\end{equation}
from the following inequalities in Theorem \ref{thm:sg}:
\begin{equation*}
\|T(t)(\rho_0,\Bv_0)\|_{H^{1,0}_p(\Omega)} \leq C \|(\rho_0,\Bv_0)\|_{H^{1,0}_p(\Omega)},\quad 
\|T(t)(\rho_0,\Bv_0)\|_{H^{1,2}_p(\Omega)} \leq C  t^{-1} \|(\rho_0,\Bv_0)\|_{H^{1,0}_p(\Omega)}.
\end{equation*}
Thus \eqref{es:GL_sg_2} and \eqref{es:GL_le_1} yield that
\footnote{Without loss of generality, we can divide the time interval $\BR_+$ by $t=2,$ as the bound of $(\wt d,\wt \bff)(t)$ for $t\in [2,3]$ is uniform up to some constant C.} 
\begin{equation}\label{es:wt_df_1}
\vertiii{(\wt d,\wt \bff)(t)}_{p,q}  \leq 
\begin{cases}
C  t^{-1/2} \|(\rho_0,\Bv_0)\|_{H^{1,0}_p(\Omega)}, 
& \text{for}\,\, 0<t<2,\\
C t^{-N/(2q)} \vertiii{(\rho_0, \Bv_0)}_{p,q}, 
&\text{for}\,\, t\geq 2.
\end{cases}
\end{equation}

According to Theorem \ref{thm:wh_1}, we denote $\{\wt T(t)\}_{t\geq 0}$ for the (analytic) semigroup associated to \eqref{eq:GLame_3}. 
Then the Duhamel principle furnishes  
\begin{equation}\label{eq:wt_rhov_1}
(\wt \rho,\wt\Bv)(t)= \wt T(t)(\psi_{\infty}\rho_0, \psi_{\infty}\Bv_0) 
+ \int_0^t \wt T(t-s)(\wt d,\wt\bff)(s) \,ds.
\end{equation}
Applying Theorem \ref{thm:GLame_whole}, we obtain for $t\geq 1,$
\begin{align}
\|\wt T(t)(\psi_{\infty}\rho_0, \psi_{\infty}\Bv_0) \|_{L_p(\BR^N)} 
&\leq C_{p,q} t^{-(N/q-N/p)/2}\vertiii{(\rho_0, \Bv_0)}_{p,q}, \nonumber \\
\|\nabla \wt T(t)(\psi_{\infty}\rho_0, \psi_{\infty}\Bv_0) \|_{L_p(\BR^N)} 
&\leq C_{p,q} t^{-(N/q-N/p)/2-1/2}\vertiii{(\rho_0, \Bv_0)}_{p,q},  
\label{es:GL_wt_T_1}\\ \nonumber
\|\nabla^2 \CP_v \wt T(t)(\psi_{\infty}\rho_0, \psi_{\infty}\Bv_0) \|_{L_p(\BR^N)} 
&\leq C_{p,q} t^{-(N/q-N/p)/2-1}\vertiii{(\rho_0, \Bv_0)}_{p,q}.
\end{align}
\medskip

Next, to study the second term on the right-hand side of \eqref{eq:wt_rhov_1},
we write
\begin{align*}
\int_0^t \wt T(t-s)(\wt d,\wt\bff)(s) \,ds
&= \Big(\int_0^2 +\int_2^{t-1} +\int_{t-1}^t\Big)\wt T (t-s)(\wt d,\wt \bff)(s) \,ds\\
&=\wt I_1(t) +\wt I_2(t)+\wt I_3(t),
\quad \forall\,\, t\geq 3.
\end{align*}

\underline{Bound of $\wt I_1(t)$.}
For $0\leq s\leq 2$ and $t\geq 3,$ we have $(t-s) \geq (t+s)/5.$
Then \eqref{es:GL_wt_T_1} and \eqref{es:wt_df_1} yield for $t\geq 3,$
\begin{align}
\|\wt I_1(t) \|_{L_p(\BR^N)} &\leq C \int_0^2 (t-s)^{-(N/q-N/p)/2} 
s^{-1/2} \,ds\, \|(\rho_0,\Bv_0)\|_{H^{1,0}_p(\Omega)} \nonumber \\
&\leq C t^{-(N/q-N/p)/2} \|(\rho_0,\Bv_0)\|_{H^{1,0}_p(\Omega)}, \nonumber \\
\|\nabla \wt I_1(t)\|_{L_p(\BR^N)} 
&\leq C\int_0^2 (t-s)^{-(N/q-N/p)/2-1/2}
s^{-1/2} \,ds \,\|(\rho_0,\Bv_0)\|_{H^{1,0}_p(\Omega)}
\label{es:wt_I_1}\\ 
&\leq C t^{-(N/q-N/p)/2-1/2} 
\|(\rho_0,\Bv_0)\|_{H^{1,0}_p(\Omega)}, \nonumber \\
\|\nabla^2 \CP_v \wt I_1(t)\|_{L_p(\BR^N)} 
&\leq C \int_0^2  (t-s)^{-(N/q-N/p)/2-1}
s^{-1/2} \,ds \,\|(\rho_0,\Bv_0)\|_{H^{1,0}_p(\Omega)} \nonumber\\
&\leq C t^{-(N/q-N/p)/2-1} 
\|(\rho_0,\Bv_0)\|_{H^{1,0}_p(\Omega)}.\nonumber
\end{align}
\smallbreak

\underline{Bound of $\wt I_2(t)$.} 
As $\supp (\wt d, \wt \bff) \subset \Omega_{5R},$ 
observe that 
$$\vertiii{(\wt d,\wt \bff)(t)}_{p,1} \leq C_{q,R} \vertiii{(\wt d,\wt \bff)(t)}_{p,q}, 
\quad \forall\,\,\, t>0.$$
Then \eqref{es:GL_wt_T_1} and \eqref{es:wt_df_1} imply for $t\geq 3,$
\begin{align}
\|\wt I_2(t) \|_{L_p(\BR^N)} &\leq C \int_2^{t-1} (t-s)^{-(N-N/p)/2} 
 s^{-N/(2q)}  \,ds \,\vertiii{(\rho_0, \Bv_0)}_{p,q} ,\nonumber \\
\|\nabla \wt I_2(t)\|_{L_p(\BR^N)} 
&\leq C\int_2^{t-1} (t-s)^{-(N-N/p)/2-1/2}
 s^{-N/(2q)}  \,ds \, \vertiii{(\rho_0, \Bv_0)}_{p,q},
\label{es:wt_I_2_1}\\ 
\|\nabla^2 \CP_v \wt I_2(t)\|_{L_p(\BR^N)} 
&\leq C \int_2^{t-1}  (t-s)^{-(N-N/p)/2-1}
s^{-N/(2q)}\,ds\, \vertiii{(\rho_0, \Bv_0)}_{p,q}. \nonumber
\end{align}
To continue the discussions, we write for $\kappa >0,$
\begin{align*}
\int_2^{t-1} (t-s)^{-\kappa} s^{-N/(2q)}  \,ds
& = \Big(\int_2^{2t/3} +\int_{2t/3}^{t-1}\Big)
 (t-s)^{-\kappa} s^{-N/(2q)} \,ds \\
&\leq \int_2^{2t/3} \big((t+s)/5\big)^{-\kappa} s^{-N/(2q)} \,ds\\
& \quad +\int_{2t/3}^{t-1} (t-s)^{-\kappa} \big(2(t+s)/5\big)^{-N/(2q)} \,ds.
\end{align*}
Then the similar argument to the bound \eqref{es:GL_I_2_1} gives us
\begin{align}
\|\wt I_2(t) \|_{L_p(\BR^N)} 
&\leq C t^{-(N/q-N/p)/2} \vertiii{(\rho_0, \Bv_0)}_{p,q}, \nonumber \\
\|\nabla \wt I_2(t)\|_{L_p(\BR^N)} 
&\leq C t^{-\sigma_1(p,q,N)} 
\vertiii{(\rho_0, \Bv_0)}_{p,q}, 
\label{es:wt_I_2_2}\\ 
\|\nabla^2 \CP_v \wt I_2(t)\|_{L_p(\BR^N)} 
&\leq C t^{-\sigma_2(p,q,N)} 
\vertiii{(\rho_0, \Bv_0)}_{p,q}, 
\quad \forall\,\, t\geq 3.\nonumber
\end{align}
\smallbreak

\underline{Bound of $\wt I_3(t)$.}
According to Theorem \ref{thm:wh_1} and the standard semigroup theory, we have 
\begin{align}
\|\wt T(t) (d,\bff)\|_{H^{1,0}_p(\BR^N)} &\leq C \|(d,\bff)\|_{H^{1,0}_p(\BR^N)},
\nonumber \\
\|\nabla \wt T(t) (d,\bff)\|_{L_p(\BR^N)} &\leq C t^{-1/2}
\|(d,\bff)\|_{L_p(\BR^N)}, 
\label{es:GL_wt_T_2}\\ \nonumber
\|\nabla^2 \CP_v \wt T(t) (d,\bff)\|_{L_p(\BR^N)} 
&\leq C t^{-1/2} \|(d,\bff)\|_{H^{1,1}_p(\BR^N)},
\end{align}
for $0<t<2$ and $(d,\bff) \in H^{1,1}_p(\BR^N).$ 
In addition, \eqref{es:wt_df_1} and \eqref{es:GL_le_1} yield that
\begin{equation}\label{es:wt_df_2}
 \|(\wt d,\wt\bff)(s)\|_{H^{1,1}_p(\BR^N)} 
 \leq C \big( \|(\wt d,\wt\bff)(s)\|_{H^{1,0}_p(\Omega_{5R})} +\|\Bv(t)\|_{H^2_p(\Omega_{5R})}\big)
\leq  C t^{-N/(2	q)} \vertiii{(\rho_0, \Bv_0)}_{p,q}
\end{equation}
for all $s\geq 2.$ Thus we infer from \eqref{es:GL_wt_T_2} and \eqref{es:wt_df_2} that
\begin{align}
\|\wt I_3(t) \|_{L_p(\BR^N)} &\leq C \int_{t-1}^{t} 
 \|(\wt d,\wt\bff)(s)\|_{H^{1,0}_p(\BR^N)}\,ds 
\leq C t^{-N/(2	q)} \vertiii{(\rho_0, \Bv_0)}_{p,q} , \nonumber \\
\|\nabla \wt I_3(t) \|_{L_p(\BR^N)} 
&\leq C \int_{t-1}^{t} (t-s)^{-1/2}
\|(\wt d,\wt\bff)(s)\|_{L_p(\BR^N)}\,ds
\label{es:wt_I_3}  \\ 
&\leq C t^{-N/(2	q)} \vertiii{(\rho_0, \Bv_0)}_{p,q}, \nonumber\\ 
\|\nabla^2 \CP_v \wt I_3(t)\|_{L_p(\BR^N)} 
&\leq C \int_{t-1}^{t} (t-s)^{-1/2}
\|(\wt d,\wt\bff)(s)\|_{H^{1,1}_p(\BR^N)} \,ds \nonumber\\
& \leq C  t^{-N/(2	q)} \vertiii{(\rho_0, \Bv_0)}_{p,q},
\quad \forall\,\, t\geq 3 .\nonumber
\end{align}

At last, we combine \eqref{es:GL_wt_T_1}, \eqref{es:wt_I_1}, \eqref{es:wt_I_2_2} and \eqref{es:wt_I_3},
\begin{align}
\|(\wt \rho,\wt\Bv)(t) \|_{L_p(\BR^N)} 
&\leq C t^{-(N/q-N/p)/2} \vertiii{(\rho_0, \Bv_0)}_{p,q}, \nonumber \\
\|\nabla(\wt \rho,\wt\Bv)(t)\|_{L_p(\BR^N)} 
&\leq C t^{-\sigma_1(p,q,N)} 
\vertiii{(\rho_0, \Bv_0)}_{p,q}, 
\label{es:wt_rhov}\\ 
\|\nabla^2\wt\Bv (t)\|_{L_p(\BR^N)} 
&\leq C t^{-\sigma_2(p,q,N)} 
\vertiii{(\rho_0, \Bv_0)}_{p,q}, 
\quad \forall\,\, t\geq 3.\nonumber
\end{align}
Furthermore, we immediately infer from \eqref{es:wt_rhov} that  
\begin{align}
\|(\rho,\Bv)(t) \|_{L_p(\Omega\backslash B_{5R})} 
&\leq C t^{-(N/q-N/p)/2} \vertiii{(\rho_0, \Bv_0)}_{p,q}, \nonumber \\
\|\nabla (\rho,\Bv)(t) \|_{L_p(\Omega\backslash B_{5R})} 
&\leq C t^{-\sigma_1(p,q,N)} 
\vertiii{(\rho_0, \Bv_0)}_{p,q}, 
\label{es:GL_rhov_1} \\  \nonumber
\|\nabla^2 \Bv (t)\|_{L_p(\Omega\backslash B_{5R})}  
&\leq C  t^{-\sigma_2(p,q,N)}
\vertiii{(\rho_0, \Bv_0)}_{p,q}
\end{align}
Then \eqref{es:GL_le_1} and \eqref{es:wt_rhov} yield the desired bound for 
$t\geq 3.$
This completes the proof of Theorem \ref{thm:main_LpLq}.

\appendix

\section{Principal of the linearization}
\label{sec:PL}
\subsection{Formulation via partial Lagrange coordinates}
Let us define
\begin{equation}\label{def:LL}
x =X_{\bu}(y,t)= y + \int^t_0 \kappa(y) \bu(y, s)\,ds \in \Omega_t \cup \Gamma_t, 
\,\,\,\forall \,\, y\in \Omega \cup \Gamma,
\end{equation}
for some vector $\bu=\bu(\cdot,s)$ defined in $\Omega\cup \Gamma,$ and $0\leq t\leq T.$
By assuming the condition
\begin{equation}\label{1.4}
\int^T_0\|\kappa(\cdot)\bu(\cdot, s)\|_{H^1_\infty(\Omega)}\,ds \leq \delta<1/2
\end{equation}
for small constant $\delta>0.$
we denote $X_{\bu}^{-1}(\cdot,t)$ for the inverse of $X_{\bu}(\cdot,t)$ in 
\eqref{def:LL}. Suppose that 
\begin{equation*}
\rho(x, t) = \eta \big( X^{-1}_\bu(x, t), t \big), \quad 
\bv(x, t) = \bu \big( X^{-1}_\bu(x, t), t \big), \quad 
\Omega_t = \big\{x = X_{\bu}(y,t) \mid y \in \Omega \big\}, 
\end{equation*}
solve \eqref{1.1} for some function $\eta$ defined in $\Omega.$
We will derive the equations formally satisfied by $(\rho,\bu)$ in $\Omega$ in what follows.

Assume that $\Gamma$ is a compact hypersurface of $C^2$ class. 
The kinematic (non-slip) condition  
$V_{\Gamma_t}=\bv\cdot\bn_t$ is automatically satisfied under
the transformation $X_{\bu},$ because $\kappa = 1$ near the boundary $\Gamma.$  
The calculations for the rest equations in \eqref{1.1} are very similar to the full Lagrangian transformation case. 
Denote that
$$\frac{\pd x}{\pd y} =\nabla_y X_{\bu}= \bI + \int^t_0\nabla_y(\kappa(y)\bu(y, s))\,ds,$$
and $J_{\bu}=\det (\nabla_y X_{\bu}).$
Then by the assumption \eqref{1.4}, there exists the inverse of $\nabla_y X_{\bu},$ that is, 
$$\frac{\pd y}{\pd x} =\big(\nabla_y X_{\bu}\big)^{-1}
= \bI + \bV_0(\bk), \quad 
\bk =\int^t_0\nabla_y(\kappa(y)\bu(y, s))\,ds,
 $$
where $\bV_0(\bk)=[V_{0ij}(\bk)]_{N\times N}$ is a matrix-valued function given by
$$\bV_0(\bk)= \sum_{j=1}^{\infty} (-\bk)^j.$$
In particular, $\bV_0(0) = 0.$
By the chain  rule, we introduce the gradient, divergence and stress tensor operators with respect to the transformation \eqref{def:LL},
\begin{gather}
\nabla_{\bu} = \big(\bI + \bV_0(\bk) \big)\nabla_y,\quad
\dv_{\bu} \bu = \big(\bI+\bV_0(\bk)\big):\nabla_y \bu
= J^{-1}\dv_y \Big( J \big(\bI+ \bV_{0}(\bk) \big)^\top 
\bu\Big), \nonumber\\
\bD_{\bu} (\bu) =\big(\bI + \bV_0(\bk) \big)\nabla \bu 
+(\nabla \bu)^{\top}\big(\bI + \bV_0(\bk) \big)^{\top}
= \bD(\bu)+\bV_0(\bk) \nabla \bu
+\big(\bV_0(\bk) \nabla \bu\big)^{\top},
\label{sym:LL_1} \\
\bS_{\bu}(\bu)=\mu\bD_{\bu} (\bu)
 +(\nu-\mu) (\dv_{\bu} \bu) \bI,\quad 
\DV_{\bu} \bA= J_{\bu}^{-1} \DV_{y} \Big(J_{\bu}\bA \big(\bI + \bV_0(\bk) \big) \Big). \nonumber
\end{gather}
In addition, the $i$th component of $\DV_{\bu} \bA$ can be also written via
\begin{equation}\label{eq:LL_DV_1}
(\DV_{\bu} \bA)_i = 
\sum_{j,k=1}^N \big[\bI + \bV_0(\bk) \big]_{jk}\pd_k A_{ij},\quad 
\forall\,\,i=1,\dots,N.
\end{equation}
In particular, $\DV_{\bu} \bA=0$ if $\bA$ is a constant matrix.
Then according to \eqref{sym:LL_1}, $(\rho,\bu)$ fulfils
\begin{equation}\label{eq:LL_CNS_1}
	\left\{ \begin{aligned}
&\pd_t\eta + (1-\kappa) \bu \cdot \nabla_{\bu}\eta
+ (\rho_e+\eta)\, \dv_{\bu} \bu=0
&&\quad&\text{in} &\quad \Omega \times (0,T), \\
&(\rho_e+\eta)\big(\pd_t\bu +(1-\kappa) \bu\cdot\nabla_{\bu}\bu \big) 
-\DV_{\bu}\big(\bS_{\bu}(\bu) - P(\rho_e+\eta)\bI\big) = 0 
&&\quad&\text{in}& \quad \Omega \times (0,T), \\
&\big(\bS_{\bu}(\bu) - P(\rho_e+\eta)\bI\big)\bn_{\bu} = P(\rho_e)\bn_{\bu}
&&\quad&\text{on}& \quad \Gamma \times (0,T), \\
&(\eta, \bu)|_{t=0} = (\rho_0, \bv_0)
&&\quad&\text{in}& \quad \Omega, \\
\end{aligned}
\right.
\end{equation}
where $\bn_{\Gamma}$ 
denotes for the unit normal vector to $\Gamma,$ 
and $\bn_{\bu}$ is defined by
$$\bn_{\bu} = \frac{\big(\bI + \bV_0(\bk)\big)\bn_{\Gamma}}{\big|\big(\bI + \bV_0(\bk)\big)\bn_{\Gamma} \big|} \cdot 
$$
It is clear that the boundary condition in \eqref{eq:LL_CNS_1} is equivalent to
\begin{equation}\label{1.7}
\Big(\bS_{\bu}(\bu) - \big(P(\rho_e+\eta)-P(\rho_e)\big)\bI\Big)\big(\bI + \bV_0(\bk)\big)\bn_{\Gamma} = 0.
\end{equation}

\subsection{Modified equations from \eqref{eq:LL_CNS_1}}
In this subsection, we derive another linearized form of \eqref{1.1} from \eqref{eq:LL_CNS_1}.
We assume that $$\bk =\int^t_0\nabla_y(\kappa(y)\bu(y, s))\,ds$$ satisfies \eqref{1.4} .
Given a smooth function 
$$G=G(\bk):\BR^{N\times N} \rightarrow \BR,$$ 
we use Taylor's theorem 
\begin{align*}
G(\bk)
& = G\Big(\int^T_0\nabla_y\big(\kappa(y)\bu(y, s)\big)\,ds 
- \int^T_t\nabla_y\big(\kappa(y)\bu(y, s)\big)\,ds\Big)\\
& =  G(\bK) + \sum_{|\ell|=1,\,\beta\in \BN_0^{N^2}} R_{\ell} \Big(\int^T_t\nabla_y \big(\kappa(y)\bu(y,s)\big)\,ds\Big)^{\ell}\\
&=G(\bK) +\wt G,
\end{align*}
with $\bK = \int^T_0\nabla_y\big(\kappa(y)\bu(y, s)\big)\,ds,$ $\wt G=G(\bk)-G(\bK),$ and $R_{\ell}$ given by
$$\quad 
R_{\ell} = -\int^1_0 (\pd_{\ell}G)\Big(\bK
- \theta\int^T_t\nabla_y\big(\kappa(y)\bu(y, s)\big)\,ds\Big)\,d\theta.$$
With this notation $\bK,$ we introduce that
\begin{gather}
J(\bK)= \Big(\det \big(\bI+\bV_0(\bK)\big) \Big)^{-1}
=\det \big(\nabla_y X_{\bu}(T,y)\big),\quad \dv_{T} \bu = \big(\bI+\bV_0(\bK)\big):\nabla_y \bu,\nonumber\\
\bD_{T} (\bu) =\big(\bI + \bV_0(\bK) \big)\nabla \bu 
+(\nabla \bu)^{\top}\big(\bI + \bV_0(\bK) \big)^{\top},
\label{sym:LL_2} \\
\bS_{T}(\bu)=\mu \bD_{T} (\bu)
 +(\nu-\mu) (\dv_{T} \bu) \bI,\quad 
\DV_{T} \bA= J(\bK)^{-1} \DV_{y} \Big(J(\bK)\bA \big(\bI + \bV_0(\bK) \big) \Big), \nonumber
\end{gather}
Then the symbols in \eqref{sym:LL_1} are rewritten by noting \eqref{eq:LL_DV_1} as follows,
\begin{gather}
\dv_{\bu} \bu = \dv_{T} \bu + \wt\bV_0:\nabla_y \bu,\quad
\bD_{\bu} (\bu)= \bD_T(\bu)+\wt\bV_0 \nabla \bu+(\wt\bV_0 \nabla \bu)^{\top}, \nonumber\\
\bS_{\bu} (\bu)=\bS_{T} (\bu)
+\mu\big( \wt\bV_0 \nabla \bu+(\wt\bV_0 \nabla \bu)^{\top} \big)
+ (\nu-\mu)(\wt\bV_0:\nabla_y \bu)\bI, \label{sym:LL_3}\\
\DV_{\bu} \bA=\DV_T \bA + (\wt \bV_0 \nabla|\bA), \nonumber
\end{gather}
where the $i$th component of $(\bB \nabla|\bA)$ equals to $\sum_{j,k} B_{jk} \pd_k A_{ij}.$
\smallbreak

On the other hand, as $P(\cdot)$ is smooth, we infer from Taylor's theorem that
\begin{equation}\label{LL:P_1}
P(\rho_e+\eta)-P(\rho_e)=P'(\rho_e) \eta 
+ \frac{\eta^2}{2}\int_0^1 P''(\rho_e +\theta \eta)(1-\theta)\,d\theta.
\end{equation}
Then we note from \eqref{LL:P_1} that 
\begin{equation}\label{LL:stress_1}
\bS_{\bu}(\bu) - \big(P(\rho_e+\eta)-P(\rho_e)\big)\bI 
=\bS_{T}(\bu) - P'(\rho_e)\eta\bI \, + \bV_1(\eta,\bu),
\end{equation}
with $\bV_1(\eta,\bu)$ given by
\begin{align*}
\bV_1(\eta,\bu)=&\mu\big( \wt\bV_0 \nabla \bu
+(\wt\bV_0 \nabla \bu)^{\top} \big)
+ (\nu-\mu)(\wt\bV_0:\nabla_y \bu)\bI\\
&+ \frac{\eta^2}{2}\int_0^1 P''(\rho_e +\theta \eta)(1-\theta)\,d\theta \bI.
\end{align*}

Then thanks to \eqref{sym:LL_3} and \eqref{LL:stress_1},
\eqref{eq:LL_CNS_1} with \eqref{1.7} can be reformulated by
\begin{equation}\label{eq:LL_CNS_2}
	\left\{ \begin{aligned}
&\pd_t\eta 
+ \rho_e\, \dv_{T} \bu  =f(\eta,\bu)
&&\quad&\text{in} \quad \Omega \times (0,T), \\
&\rho_e\pd_t\bu 
-\DV_{T}\big(\bS_{T}(\bu) - P'(\rho_e)\eta\bI\big) = \bg(\eta,\bu)
&&\quad&\text{in} \quad \Omega \times (0,T), \\
&\big(\bS_{T}(\bu) - P'(\rho_e)\eta\bI\big)
\big(\bI + \bV_0(\bK)\big)\bn = \bh(\eta,\bu)
&&\quad&\text{in} \quad \Gamma \times (0,T), \\
&(\eta, \bu)|_{t=0} = (\rho_0, \bv_0),
\end{aligned}
\right.
\end{equation}
where the nonlinear terms on the r.h.s. of  \eqref{eq:LL_CNS_2} are given  by
\begin{align*}
f(\eta,\bu)&=-(1-\kappa) \bu \cdot \nabla_{\bu}\eta-\eta\, \dv_{\bu} \bu-\rho_e \wt\bV_0:\nabla_y \bu,\\
\bg(\eta,\bu)&=-\eta\pd_t \bu-
(\rho_e+\eta)(1-\kappa) \bu\cdot\nabla_{\bu}\bu \\
&\quad +\DV_{\bu}\bV_1(\eta,\bu)
+\big(\wt \bV_0 \nabla|\bS_{T}(\bu) - P'(\rho_e)\eta\bI \big),\\
\bh(\eta,\bu)&=-\Big(\bS_{\bu}(\bu) - 
\big(  P'(\rho_e) \eta 
+ \frac{\eta^2}{2}\int_0^1 P''(\rho_e +\theta \eta)(1-\theta)\,d\theta       \big)\bI\Big)\wt\bV_0 \bn\\
&\quad -\bV_1(\eta,\bu)\big(\bI + \bV_0(\bK)\big)\bn.
\end{align*}
At last, it is easy to see \eqref{LL:Lame_1} from \eqref{eq:LL_CNS_2}.


\end{document}